\documentclass[a4paper,11pt]{report}

\usepackage[frenchb,francais]{babel}
\usepackage[T1]{fontenc}
\usepackage{amsmath}
\usepackage{amssymb}
\usepackage{oldgerm}
\usepackage{graphicx}
\usepackage{eufrak}
\usepackage{euscript}
\usepackage{epsfig}
\begin{document}
\def\gO{\mathcal{O}}
\renewcommand{\thefootnote}{\fnsymbol{footnote}}
\newcommand{\RR}{\mathbb{R}}
\newcommand{\NN}{\mathbb{N}}
\newcommand{\Zed}{\mathbb{Z}}
\newcommand{\rto}{\rightarrow}
\newcommand{\f}[2]{\ensuremath{\mathchoice%
        {\dfrac{#1}{#2}}
        {\dfrac{#1}{#2}}
        {\frac{#1}{#2}}
        {\frac{#1}{#2}}
        }}
\newcommand{\Dp}[2]{\ensuremath{\f{\partial#1}{\partial#2}}}
\newcommand{\CC}{\mathbb{C}}
\newtheorem{Th}{Théorème}
\newtheorem{Lemme}[Th]{Lemme}
\newtheorem{Cor}[Th]{Corollaire}
\newtheorem{Lemmetech}{Lemme}
\newtheorem{Def}{Définition}
\def\Dom{{\EuScript D}}
\def\HH{\hbox{\large{$\EuFrak H$}}}
\def\ZZ{\hbox{\large{$\EuFrak Z$}}}
\def\HHH{{\mathcal H}}
\def\ZZZ{{\mathcal Z}}
\def\Hh{\hbox{\large{$\EuScript H$}}}
\def\Zz{\hbox{\large{$\EuScript Z$}}}
\def\gY{{\mathcal Y}}
\def\eps{\varepsilon}
\def\Ups{{\boldsymbol\Upsilon}}
\def\gR{{\mathsf R}}
\def\gF{{\mathcal F}}
\def\gG{{\mathcal G}}
\def\DD{\mathcal{D}}
\def\Ai{\mathsf{Ai}}
\def\HB{H$\bullet$}
\def\EuH{{\EuFrak H}}
\def\EuY{{\EuFrak Y}}
\def\Euy{{\EuFrak y}}
\def\EuO{{\EuFrak o}}

\def\og
{\leavevmode\raise.3ex\hbox{$\scriptscriptstyle\langle\!\langle$\hglue 
1pt\penalty 10000}}

\def\fg{\hglue 1pt\penalty 10000
\leavevmode\raise.3ex\hbox{$\scriptscriptstyle\,\rangle\!\rangle$\hskip 2pt
\penalty -1500}}

\thispagestyle{empty}
{\large
\vskip 12.5cm
\centerline{\bf ÉTUDE D'ÉQUATIONS DIFFÉRENTIELLES ORDINAIRES} 
\centerline{\bf SINGULIÈREMENT PERTURBÉES,} 
\centerline{\bf AU VOISINAGE D'UN POINT TOURNANT}

\vfill\eject

\hbox{}
\newpage

\tableofcontents
\newpage
\listoffigures
\newpage

\hbox{}
\vskip 70pt plus 5pt

\section*{Introduction}

\bigskip
\bigskip
Soit l'équation de Van der Pol forcée
\begin{equation*}
\eps \ddot u+(u^2-1)\dot u +u=\alpha\,,
\end{equation*}
où $\eps$ est un paramètre que l'on fait tendre vers $0$, et où $\dot
u$ désigne la dérivée de $u$ par rapport au temps $t$.

Cette équation est équivalente au système différentiel
\[ \left\{ \begin{array}{rl}
   \eps \dot u &=y-\frac{u^3}{3}+u\\
      \dot y   &=\alpha -u
           \end{array}
\right. \]
dont on a découvert il y a plus de vingt ans certaines solutions
particulières, appelées grands canards: pour certaines valeurs de
$\alpha$, dépendant de $\eps$, ces solutions restent au voisinage de
la courbe lente $y=\f{u^3}{3}-u$ (obtenue en posant $\eps=0$) y
compris le long de la partie instable de cette courbe, entre les deux
sommets $y=-1$ et $y=1$ (cf. chapitre~\ref{chap:vdp}). Des études avec
les outils de l'analyse non standard, sur la droite réelle, ont permis
de montrer que l'existence de canards est souvent inéluctable, mais
qu'ils ont une durée de vie extr\^emement courte: les valeurs du
paramètre $\alpha$ menant à des canards sont toutes dans un intervalle
de taille $\exp\left(-M/\eps\right)$ (quand $\eps\rightarrow 0^+$, et
pour un certain $M>0$).

Ce système étant lui-même équivalent à l'équation
\begin{equation}\tag{\ref{eq:vdpdebase}}
\eps v\frac{dv}{du}=(1-u^2)v+\alpha -u\,,
\end{equation} 
l'étude de ces solutions canards a pu progresser par le passage d'une
variable $u$ réelle à une variable dans le champ
complexe~\cite{FS}. La simplification de l'étude vient essentiellement
du fait que les fonctions solutions, au lieu de n'être que de classe
${\cal C}^1$, sont alors toutes holomorphes. Par ailleurs, le lien
avec les points tournants apparaît alors clairement.

On constate d'abord que cette équation admet une solution formelle
$$\hat v(u,\eps)=\sum_{n=0}^\infty v_n(u) \eps^n\,,\ \ 
\hat\alpha(\eps)=\sum_{n=0}^\infty a_n \eps^n\,,$$
dont la propriété (exceptionnelle) est que tous les termes $v_n(u)$
sont analytiques en $u=1$. Ensuite, on peut chercher à faire le
lien entre cette solution formelle et des solutions holomorphes,
bornées en $u=1$ quand $\eps$ tend vers $0$. Si on impose des domaines
d'existence les plus grands possibles, on trouve exactement deux
solutions $v^+(u,\eps)$ et $v^-(u,\eps)$ pour leurs paramètres
respectifs $\alpha^+(\eps)$ et $\alpha^-(\eps)$. Avec ces solutions qui
sont définies toutes deux sur le segment $u\in [-1+\rho,1]$, pour tout
$\rho>0$ fixé, on est capable de donner une valeur à $M$; on pourra
choisir
$$M=R(-1+\rho)-R(1) =
\Re\left(\int_1^{-1+\rho} (t-1)(t+1)^2dt\right)\,.$$ 
Outre le premier point tournant $u=1$, on voit ici que l'autre point
tournant $u=-1$ joue aussi un rôle important. Et qu'il pourrait être
intéressant de s'approcher de ce second point, pour remplacer
$-1+\rho$ par $u_1=-1+o_\eps(1)$, $u_1>-1$. D'où l'idée d'étudier le
voisinage des points tournants, pour montrer qu'on peut toujours
contrôler la croissance, quand $\eps\rightarrow 0$, des solutions
$v^+$ et $v^-$ jusqu'en un tel point $u_1$.

\bigskip
Considérons le cas général d'une équation différentielle du premier
ordre non linéaire
$$\eps y'=y f(x,\eps)+h(x,\eps)+\eps y^2 P(x,\eps,y)\,,$$ 
où $\eps$ est un petit paramètre complexe, où $y'$ désigne la dérivée
$\f{dy}{dx}$, et où les fonctions $f$, $h$ et $P$ sont holomorphes en
$x$ et $y$ dans certains domaines et admettent un développement en
série (éventuellement non convergent) en $\eps$.

On peut calculer une solution formelle du type
$$\hat y(x,\eps)=\sum_{n=0}^\infty y_n(x) \eps ^n$$ 
pour cette équation, à l'aide d'une récurrence sur les coefficients
$y_n(x)$. On peut ensuite au moins espérer~\cite{Wasow1} que cette
série en $\eps$, le plus souvent divergente, corresponde à une
solution holomorphe $y(x,\eps)$ quand $\eps$ tend vers $0$ dans des
secteurs $\arg\eps\in\ ]\theta_1,\theta_2[\,$. On sait d'ailleurs 
qu'on a effectivement une solution, sur un voisinage de tout point
$x$ tels que $f(x,0)\neq 0$, qui reste bornée quand $\eps$ tend vers 
$0$~\cite{Sib}.

En revanche, les points $x_0$ tels que $f(x_0,0)=0$ font partie le
plus souvent~\cite{Wasow2} des points tournants pour cette équation:
on voit bien que les coefficients $y_n$ calculés sont alors en général
singuliers en $x=x_0$. Il para\^it alors difficile que la solution
formelle $\hat y$ puisse être le développement asymptotique d'une
solution holomorphe dans tout un voisinage de $x_0$ (la non-existence
d'une telle solution est précisément la définition d'un point
tournant). Cependant, dans certains cas particuliers, il existe tout
de m\^eme une série formelle dont tous les termes sont analytiques en
$x_0$; ce genre de résultat est obtenu par exemple en introduisant un
(multi-)paramètre sous la forme d'une série formelle. Si une telle
série formelle existe, cela permet d'envisager l'existence d'une
solution à l'équation différentielle, qui reste analytique en
$x_0$. Ces solutions particulières, holomorphes en un point
généralement tournant pour les solutions d'une équation différentielle
donnée, sont celles qu'on a appelées, par généralisation, solutions
canards holomorphes, ou solutions surstables. L'existence de ces
canards dépend souvent du choix d'un (multi-)paramètre $\alpha(\eps)$,
présent dans les fonctions $f$, $h$ et/ou $P$, qui s'écrit lui aussi,
bien sûr, au moins comme une série formelle en $\eps$.

On a ainsi pu démontrer \cite{BFSW} dans un cas très général que les
valeurs du paramètre $a$ donnant des canards sont exponentiellement
proches en $\eps$, et qu'on peut trouver un bon majorant de la
constante $M$ généralisée de la manière suivante: si $\alpha^+(\eps)$
et $\alpha^-(\eps)$ sont deux valeurs à canards,
$$\|\alpha^+-\alpha^-\| < C \exp\left(-\f{M}{\eps}\right)\,,
\text{ pour tout $\eps>0$ assez petit,} $$ 
et pour un certain $C$ assez grand dépendant de la valeur de $M$
choisie; comme dans le cas de l'équation de Van der Pol, $M$ peut
s'écrire $M=R(x_1)-R(x_0)$, où $R$ désigne l'application
$R(x)=\Re\left(\int_0^x f(t,0) dt\right)$; $x_1$ est à nouveau un
point dans le domaine d'existence $\Dom$ de plusieurs solutions
canards \og maximales\fg, point qu'on a intérêt à prendre le plus
proche possible d'un autre point tournant de l'équation pour obtenir
la meilleure inégalité.

L'idée des recherches entreprises ici a donc été d'étudier les
solutions au voisinage des points tournants pour voir si cela ne
permettait pas de s'approcher de ces points (i.e. de prendre un $x_1$
qui tend vers une singularité quand $\eps$ tend vers $0$), et donc
de préciser l'inégalité ci-dessus. L'intérêt principal d'un tel
résultat étant qu'il permet ensuite de donner une estimation ou au
moins une bonne majoration des coefficients des séries formelles
correspondant à $\alpha$.

\bigskip
Dans le cas de l'équation de Van der Pol, au chapitre~\ref{chap:vdp},
on montrera que les solutions canards $v^+$ et $v^-$ existent
effectivement dans des domaines contenant le segment $[-1+X_l
\eps^{1/3}, 1]$ (au lieu de $[-1+\rho,1]$ pour tout $\rho$ strictement
positif), ceci pour une constante $X_l$ suffisamment grande, mais qui
peut être choisie indépendante de $\eps$; la démonstration utilisera
une équation intérieure équivalente à \eqref{eq:vdpdebase}, mais non
singulièrement perturbée. La propriété démontrée permet ensuite
effectivement de donner un équivalent pour les coefficients $a_n$ des
paramètres $\alpha$ à canards:
$$a_n\underset{n\infty}\sim \frac{-2\sqrt 6}{\pi\sqrt\pi 
e^{4/3}}n^{-1/2}\left(\frac{3}{4}\right)^n n!$$

\medskip
Le résultat de l'étude des solutions au voisinage des points tournants
sera généralisé et complété au §~\ref{sec:pres} pour une équation
générale  que l'on réécrit sous la forme
\begin{equation*}
\eps y'=\bigl(x^p f(x)+\eps g(x,\eps,\alpha)\bigr) y +h(x,\eps,\alpha)
+\eps y^2 P(x,\eps,\alpha, \eps y)\,,
\end{equation*}
dont on étudie le point tournant $0$ (pour $p>0$).  

Cette équation singulièrement perturbée admet une solution (dite:
extérieure) à distance finie du point tournant à l'intérieur de
certains secteurs centrés en ce point tournant; la démonstration de
l'existence de cette fonction est rappelée au paragraphe
§~\ref{sec:loin}. Cette première équation est équivalente à une
équation intérieure; celle-ci, si elle est régulière en $\eps$,
possède au moins une solution particulière, qui peut \^etre connectée
avec la solution extérieure. Le domaine en $x$ de toute solution
holomorphe bornée en $\eps$ de l'équation de départ peut donc être
étendu jusqu'à une distance $X_l\eps^{1/(p+1)}$ des points tournants:
\[\Dom=\left\{x\in\CC\,/\,|x|>X_l \eps^{1/(p+1)}\ \text{ et } 
\arg(x)\in\,\, ]\theta_1,\theta_2[\,\right\}\, ,\text{ pour des $X_l$
assez grands}\,.\] 
Dans ce nouveau domaine, cette solution n'est pas
nécessairement bornée quand $\eps\rightarrow 0$, mais sa croissance
reste cependant contr\^olée et ne devient pas exponentielle.

\smallskip
Au §~\ref{sec:alpha}, à l'aide de ce résultat, on démontre
l'existence, dans certains cas, de solutions surstables dans un
voisinage complet d'un point tournant, si on choisit correctement le
(p-multi-)paramètre $\alpha$ en fonction de $\eps$.

\medskip
Au chapitre~\ref{chap:Brusselator}, on s'intéresse alors à l'équation
du Brusselator, transfomée en l'équation différentielle suivante, avec
$a$ pour paramètre,
$$\eps z(x)\frac{dz(x)}{dx}=\frac{2x}{(1+x)^2}
\Bigl(z-\frac{1}{2(1+x)^3}\Bigr)-\frac{a-1}{(1+x)^2}z
-z\Bigl(z-\frac{1}{2(1+x)^3}\Bigr)\frac{2\eps}{1+x}\,.$$ 
Le système du Brusselator équivalent est connu pour avoir des
solutions canards pour certaines valeurs d'un paramètre, mais l'étude
du domaine de ces solutions dans le plan complexe n'avait pas encore
été faite. Gr\^ace à la généralisation du chapitre précédent, il
devient facile de redémontrer l'existence des canards holomorphes et
de décrire leur domaine d'existence étendu. Comme pour le cas de
l'équation de Van der Pol, on peut ensuite trouver un équivalent
(quand $n\rightarrow \infty$) pour les coefficients des paramètres à
canards.

\medskip
À noter que, pour l'équation du Brusselator comme pour l'équation de
Van der Pol, il faut commencer par réussir à calculer des coefficients
de Stokes pour une équation différentielle non linéaire avant de
pouvoir donner les équivalents cités. Les multiplicateurs de Stokes,
qui apparaissent autour d'un point tournant ne sont pas exactement
calculables dans le cas non linéaire, mais dans les deux exemples
cités, on réussit à en déterminer un équivalent complet.

\chapter{Définitions, propriétés élémentaires}
\section{Notion de relief}

Dans tout ce qui suit, il sera souvent fait mention de «relief» et de
chemins descendants ce relief, notions introduites par
J.~Callot~\cite{Callot}
voici quelques années.

Pour faire comprendre son intérêt, regardons l'équation différentielle
singulièrement perturbée (c'est-à-dire avec $\eps$ qui tend vers $0$)
linéarisée générale
\begin{equation}\label{ed}\tag{E.diff.}
\eps y'=y f(x,\eps)+h(x,\eps)+\eps y^2 P(x,\eps,y)\,,
\end{equation}
où toutes les fonctions considérées sont holomorphes dans certains
domaines en toutes leurs variables. On souhaite montrer que cette
équation admet  une solution holomorphe en ses deux variables
$y(x,\eps)$. Si une telle solution existe, en utilisant la formule de
variation de la constante, on voit que nécessairement $y$ vérifie
l'égalité suivante, dans  laquelle $\gamma^s_x$ désigne un chemin
différentiable de $\CC$  allant de $s$ à $x$:
$$y(x,\eps)=y(s,\eps)
+\frac{1}{\eps}\int_{\gamma^s_x}e^\frac{F(x,\eps)-F(t,\eps)}{\eps}
\Bigl(h(t,\eps) +\eps y^2 P(t,\eps,y)\Bigr)dt,$$
où $F(x,\eps)=\int^x f(t,\eps) dt$. Dans tous les cas envisagés 
ci-après, on supposera que $F(x,\eps)$ est peu différent de $F(x,0)$: 
on aura au pire $F(x,\eps)-F(x,0)=\gO(\eps)$ (ce sera bien s\^ur le 
cas si $f$ est holomorphe en $\eps=0$).

Pour voir s'il existe une solution $y$ bornée quand $\eps$ tend vers
$0$, on s'intéresse particulièrement au terme exponentiel. Pour $\eps$
réel positif, on voit bien intuitivement que si la partie réelle de
$F(x,\eps)-F(t,\eps)$ est positive pour un certain $t$, l'intégrale
est en général exponentiellement grande quand $\eps$ tend vers
$0$. C'est à cause de cela que l'on considère l'application relief qui
à tout $x$ de $\CC$ associe $R(x,\eps)=\Re(F(x,\eps))$ dans $\RR$,
dont les lignes de niveau ont pour équation
$\Re(F(x,\eps))=\mathrm{constante}$. Il paraît alors naturel de
supposer l'existence de solutions dans des domaines $\{x /\ \exists\,
\gamma^s_x$, chemin descendant le relief en chacun de ses points$\}$. 
En fait, on en demande un peu plus:
\begin{Def}\label{def:acc}
On dira d'un domaine $\Dom$
qu'il est {\bf accessible} si:
\begin{enumerate}
\item pour tout $x\in\Dom$ on peut trouver un
chemin différentiable  $\gamma_x(u)$ allant d'un point fixe $s\in\Dom$
(appelé sommet de $\Dom$) à $x$; 
\item $\gamma_x$ est tel que $\frac{dR(\gamma_x(u))}{du}
\leq -C_{\gamma_x} |F'(\gamma_x(u)) \gamma'_x(u)|$ pour tout $u$;  
\item les $C_{\gamma_x}$, tous positifs, admettent une borne inférieure
strictement positive: $\displaystyle\inf_{x\in\Dom}C_{\gamma_x}$ 
existe et est égal à $C>0$. 
\end{enumerate}
\end{Def}
Si on a une telle propriété, le terme exponentiel de l'intégrale sera
exponentiellement petit quand $\eps$ tend vers $0$, et il est alors
raisonnable d'espérer pouvoir contrôler la valeur de l'intégrale pour
tous les $x\in\Dom$.

Si $\eps$ tend vers $0$ dans une autre direction ($\arg
(\eps)=\theta\neq 0$), on change de relief, et donc de domaine
accessible, en optant cette fois pour $\Re
\frac{F(x,\eps)-F(t,\eps)}{e^{i\theta}}$, la condition d'accessibilité
restant la même avec le nouveau relief. On préfère en général
effectuer les études avec un relief fixe (qui donne des chemins
$\gamma_x$ fixe); c'est la raison pour laquelle on travaillera
toujours avec $\arg(\eps)$ constant et $|\eps|\rightarrow 0$.

\medskip Pour l'intégration numérique aussi, la propriété de relief
est importante:  il n'est pas difficile de vérifier que pour deux
solutions de l'équation différentielle $y_a$ et $y_b$ telles que
$y_a(x_0)\neq y_b(x_0)$, la différence entre les deux fonctions se
comportera de manière totalement différente (pour $\eps$
arbitrairement petit) suivant qu'on part de $x_0$ le long d'un chemin
qui monte, ou qui descend du relief.  En effet, on aura
$$(y_a-y_b)(x,\eps)=(y_a-y_b)(x_0,\eps)
+\frac{1}{\eps}\int_{\gamma^{x_0}_x}e^\frac{F(x,\eps)-F(t,\eps)}{\eps}
\Bigl(\eps dP(t,\eps,y_a,y_b)\Bigr)dt\,.$$ L'intégrale dans cette
égalité sera, si $\gamma_x$ descend le relief, exponentiellement
petite en $\eps$; ce qui signifie que les deux solutions resteront
voisines (du moins si le chemin entre $x_0$ et $x$ est de longueur
finie).  Mais l'intégrale deviendra exponentiellement grande si ce
chemin monte le relief; et donc dans ce cas $y_a$ et $y_b$
s'éloigneront très rapidement l'une de l'autre.

Lors d'une intégration numérique d'une équation différentielle du type
cité, en prenant un petit $\eps$, il sera donc très important de
partir d'un sommet $x_0$, \og haut\fg sur le relief, et d'en rester
aux points accessibles à partir de ce point.

\bigskip
Cette constatation nous amène à deux nouvelles définitions:
\begin{Def}
On appellera {\bf montagne}, pour le relief $R(x)$, par rapport au
point $x_0\in\CC$ (dans les cas cités, $x_0$ sera un point col pour le
relief $R$), un domaine de $\CC$ connexe et simplement connexe tel que
$x_0$ est accessible avec ce relief à partir de tout élément du
domaine: $R(x)>R(x_0)$.

De m\^eme, on appellera {\bf vallée} un domaine, toujours simplement
connexe, accessible à partir de $x_0$; les points $x$ du domaine
vérifient donc: $R(x)<R(x_0)$.
\end{Def}

\section{Série formelle Gevrey}

Quand on essaye de donner une solution sous forme de série à une
équation différentielle singulièrement perturbée, on arrive à des
séries divergentes, mais dont les termes vérifient cependant des
propriétés de croissance. On est amené naturellement à introduire les
séries Gevrey (cf.~\cite{Wasow1},\cite{Wasow2}).

\begin{Def}
On notera dans toute la suite le {\bf secteur} ouvert centré en $x_0$,
de rayon $\rho$ et d'angle $]\theta_1,\theta_2[$ sous la forme
$$S(x_0,\rho,\theta_1,\theta_2),\text{ soit } \left\{x\in\CC\,\Bigl/ 
0<|x-x_0|<\rho, \text{ et } \arg(x-x_0)\in\,
]\theta_1,\theta_2[\right\}$$
Si le rayon $\rho$ n'est pas précisé, il sera supposé infini.
\end{Def}
Par extension, on parlera aussi de secteurs centrés en l'infini:
$$S(\infty,\varrho,\theta_1,\theta_2)=\left\{x\in\CC\,\Bigl/ 
|x|>\varrho, \text{ et } \arg(x)\in\, ]\theta_1,\theta_2[\right\}$$

\begin{Def}
On dit d'une série formelle $$\sum_{n=0}^\infty a_n \eps^n$$ qu'elle
est {\bf Gevrey} d'ordre $1$ si la série a un rayon de convergence
nul, mais qu'il existe des constantes positives $\mathrm A$ et
$\mathrm B$ telles que, pour tout $n$
$$|a_n|\leq \mathrm{A B}^n n!$$
\end{Def}
Cette définition s'étend bien s\^ur aux séries de fonctions 
$$\sum_{n=0}^\infty y_n(x)\eps^n,$$ en remplaçant $|a_n|$ par
$\|y_n(x)\|$ où $\|\cdot\|$ désigne la norme que l'on choisit pour les
fonctions $y_n$.

\begin{Def}
On dit qu'une fonction $a(\eps)$ admet la série $\sum a_n \eps^n$
comme {\bf développement asymptotique Gevrey} dans un certain secteur
$S=S(0,\rho,\theta_1,\theta_2)$ si $a(\eps)$ est holomorphe sur $S$ et
s'il existe $\mathrm{A}$ et $\mathrm{B}$, constantes positives, telles
que:
$$\text{pour tout } N\text{ assez grand et tout $\eps$ dans $S$ }, 
\left| a(\eps)-\sum_{n=0}^{N-1} a_n \eps^n\right|
<\mathrm{A}\mathrm{B}^N N! |\eps|^N.$$ 
Dans le cas d'une série de fonctions (de la forme ci-dessus), si la
majoration est indépendante de $x$, pour tout $x$ dans un certain
domaine $\Dom$, on parle de série asymptotique {\bf uniforme} sur
$\Dom$.

\end{Def}

On sait qu'il existe une infinité de fonctions holomorphes $a(\eps)$
admettant $\sum a_n \eps^n$ comme développement asymptotique en $0$
(on écrit $a\sim \sum a_n \eps^n$).

Si $a\sim \sum a_n \eps^n$, on sait par ailleurs qu'on obtient une
 approximation exponentiellement bonne en $\eps$ de la valeur de
 $a(\eps)$ en sommant la série des $a_n\eps^n$ jusqu'à ce que, à
 $\eps$ fixé, le terme $|a_n\eps^n|$ soit minimum (on parle alors de
 sommation \og au plus petit terme\fg).

\section{Série formelle et vraie solution}

On se place dans le cas où une équation différentielle (singulièrement
perturbée) du type~\eqref{ed}
$$\eps y'=y f(x,\eps)+h(x,\eps)+\eps y^2 P(x,\eps,y)\,$$
admet une solution formelle 
$$\hat y(x,\eps) =\sum_{n=0}^\infty y_n(x)\eps^n.$$

\begin{Def}
Si la série $\sum y_n(x)\eps^n$ est solution formelle de~\eqref{ed},
on appelle classiquement {\bf«courbe lente»} la fonction $y_0(x)$,
obtenue pour la valeur du paramètre $\eps=0$.
\end{Def}

Le problème est de montrer qu'il existe une solution holomorphe
$y(x,\eps)$ (existant dans certains secteurs centrés en $0$ pour la
variable $\eps$, et dans certains domaines pour la variable $x$), qui
admette effectivement $\hat y$ comme développement asymptotique quand
$\eps$ tend vers $0$.

La réponse générale est loin d'être évidente. Un certain nombre
d'articles (\cite{KS}, \cite{CHAP}) ont été récemment publiés à ce
sujet, partant de l'équation de croissance cristalline
$$\eps^2 \Phi'''+\Phi' = \cos \Phi$$ avec les valeurs aux bornes
$$\Phi(s) \rightarrow-\pi/2\ \text{ quand } s\rightarrow-\infty,\ \ \
\Phi(s) \rightarrow \pi/2\ \text{ quand } s\rightarrow +\infty.$$
Ils démontrent que s'il existe bien une solution formelle $\hat \Phi=
\sum_{n=0}^\infty \eps^{2n}\Phi_n(s)$, où $\Phi_0$ vérifie bien les
valeurs aux bornes et $\lim_\infty \Phi_n =0$ pour tout $n$, il
n'existe pas de solution holomorphe à l'équation vérifiant les limites
voulues en l'infini (le problème venant d'un terme exponentiellement
petit en $\eps$, qui n'appara\^it donc pas dans la série). Ils donnent
aussi des ébauches de procédures pour déterminer, pour des équations
plus générales, si on est en présence de tels cas.

\bigskip

Le travail principal dans cette thèse correspond au problème inverse de celui
envisagé pour l'équation de croissance cristalline. Sachant qu'on a
une solution formelle d'une équation différentielle qui vérifie une
certaine propriété (dans notre cas, la continuité en un point qui est
en général singulier), il s'agit de montrer qu'il existe effectivement
une vraie solution à l'équation admettant cette série formelle comme
développement asymptotique. Cette solution vérifie alors la même
propriété (la continuité en ce point, ou plus exactement le fait que
cette fonction reste bornée au voisinage de ce point quand $\eps$ tend
vers $0$).

\medskip
M\^eme pour ce problème particulier, on connaît des exemples où
l'existence de solutions formelles à fonctions coefficients continues
en $x=0$ n'implique pas l'existence de vraies solutions bornées, quand
$\eps$ tend vers $0$, au voisinage de $x=0$.  En effet, pour
l'équation
$$\eps u'=xu+\exp\left(\frac{-1}{\sqrt\eps}\right),$$ 
la solution nulle est solution formelle de l'équation, qui est
évidemment analytique en $0$. Cependant, les vraies solutions s'écrivent
(en utilisant la formule de variation de la constante, qui met à
nouveau le relief en évidence, fig.~\ref{fig:relintro})
$$u(x)=u(x_0)\exp\left(\frac{x^2-x_0^2}{2\eps}\right)
+\int^x_{x_0}\frac{1}{\eps} \exp\left(\frac{-1}{\sqrt\eps}\right)
\exp\left(\frac{x^2-t^2}{2\eps}\right)dt $$
ne peuvent pas être bornées en $x$, indépendamment de $\eps$, dans un
intervalle ouvert autour de $x=0$. Ceci étant vrai quel que soit les
$x_0$ et $u(x_0)$ choisis: pour $x_0>0$, la solution sera bornée en
$\eps$ pour tout $x\geq 0$, mais ne le sera pas pour un $x<0$
quelconque fixé. On remarquera cependant que, dans cet exemple, les
coefficients de l'équation ne sont pas tous holomorphes en $\eps$.

\begin{figure}[ht!]
\centerline{\rotatebox{-90}{\includegraphics[scale=0.74]{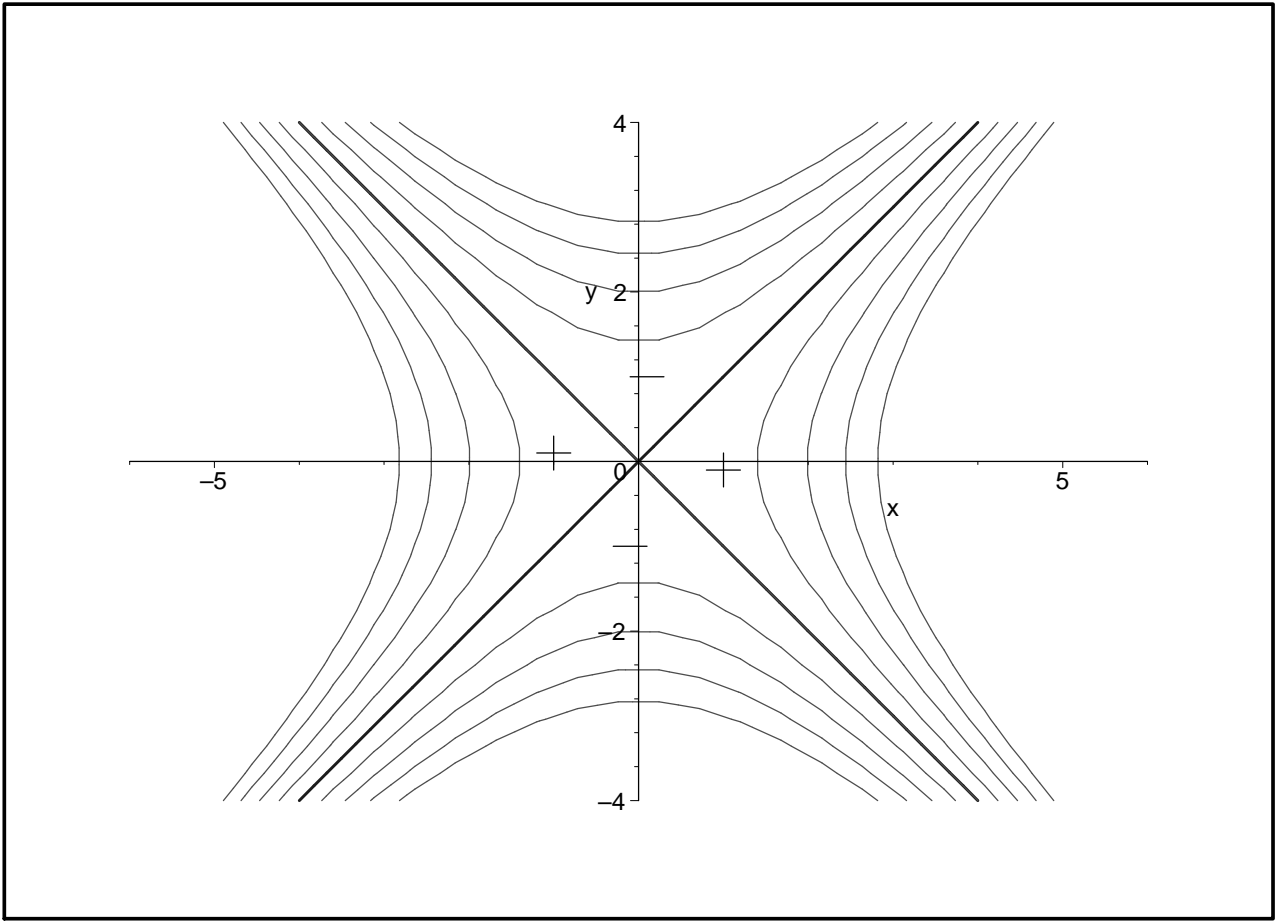}}}
\caption{Lignes de relief avec $R(x)=\Re\left(\int^x t dt\right)$; les
signes $-$ et $+$ désignent respectivement les deux vallées et
montagnes par rapport au point col $0$.}
\label{fig:relintro}
\end{figure}

\section{Points tournants d'une équation différentielle}

On trouve diverses définitions des points tournants pour une équation
différentielle dans la littérature. La plus classique~\cite{Wasow1}
pour une équation du type \eqref{ed} est la définition «par défaut»
suivante:
\begin{Def}
$x_0$ est appelé point tournant, pour une équation \eqref{ed}
admettant une solution formelle $\hat y=\sum y_n \eps^n$, s'il
n'existe pas de vraie solution de l'équation admettant comme
développement asymptotique la série formelle $\hat y$, uniformément
pour tout $x$ dans un voisinage de $x_0$ (quand $\eps$ tend vers $0$).
\end{Def}
On remarquera que les points tournants sont bien des points du plan
complexe en $x$: les séries formelles $\hat y$, sauf si elles
convergent, ne sont des développements asymptotiques de solutions de
l'équation que pour des $\eps$ tendant vers $0$ dans certains secteurs
en $\eps$ centrés en $0$.

D'après ce qui a été dit à propos du relief, on constate que les seuls
points susceptibles d'\^etre des points tournants sont les points cols
du relief donné par $\Re\Bigl(\int^x f(t)dt\Bigr)$, autrement dit là
où $f$ s'annule. Ce qui est le cas dans l'exemple précédent
(fig.~\ref{fig:relintro}). Ou bien évidemment des points singuliers de
ce m\^eme relief. On peut voir en effet que pour un point normal $x$,
il existe toujours un point «au-dessus» de $x$ à partir duquel un
voisinage complet de $x$ est accessible; et il existe donc une
solution bornée dans ce voisinage. En revanche, pour un point $x_0$
qui est un point col, ou singulier, le développement asymptotique sera
valable en général au mieux sur une montagne naissant en $x_0$ et
l'essentiel du domaine accessible à partir du sommet de cette
montagne. Pas forcément tout le domaine accessible, car la validité du
développement jusqu'en $x_0$ lui m\^eme ne soit pas forcément
assurée. Dans ce qui suit, nous étudierons comment se comportent les
solutions d'une équation différentielle au voisinage des points col
pour le relief correspondant. Nous verrons que certains de ces points
ne seront pas, exceptionnellement, des points tournants suivant la
définition donnée ci-dessus.

\section{Phénomène de Stokes}

Autour d'un point tournant, correspondant à une solution formelle
unique, on trouve en général différentes solutions existant dans des
(sous-ensemble de) secteurs ouverts centrés en ce point; on retrouvera
cette propriété classique au théorème~\ref{th:loin}. Ces domaines
d'existence sont tous différents, mais l'intersection de deux de ces
domaines distincts n'est pas forcément vide. L'explication de cette
différence de comportement de deux solutions, très voisines pour
certains $x$, quand on s'approche des limites d'un domaine (l'une des
solutions reste bornée, alors que l'autre y devient singulière) est le
phénomène de Stokes.

Dans le cas linéaire, où les domaines d'existence sont exactement des
secteurs ouverts, on sait~\cite{Wasow1} que la différence entre deux
solutions est un terme du genre $\exp\left(-{P(x)}/{\eps}\right)$. Ce
terme est exponentiellement petit en $\eps$ dans l'intersection des
domaines d'existence (ce qui explique que le développement
asymptotique reste le m\^eme), mais il devient exponentiellement grand
quand $x$ atteint les limites des secteurs: pour $\eps>0$, quand la
partie réelle de $P(x)$ devient négative.

\smallskip
Dans le cas non linéaire, le phénomène est certes plus complexe, mais
reste analogue. Quand on parlera ici de calculer un coefficient de
Stokes, il s'agira de déterminer la différence de deux solutions à une
équation donnée, au moins pour certains $x$ où ces deux solutions sont
toutes deux définies. Cependant alors que dans le cas linéaire, on
peut trouver le terme exact, on n'aura dans le cas général qu'un
équivalent (dans le meilleur des cas).

\smallskip
Le phénomène de Stokes est indissociable des points tournants: s'il
n'est pas présent, c'est qu'il existe une solution holomorphe
univalente tout autour du point tournant, correspondant à une unique
série formelle. Dans ce cas, si la série formelle est bien définie en
ce point, la solution est holomorphe jusque là. On n'a plus alors de
point tournant, par définition.

\section{Solutions {\protect «}canards{\protect »}}

C'est dans ce cas particulier où il existe une solution holomorphe
jusqu'en un point col du relief, qui est un «candidat point tournant»,
que l'on parlera de solution canard:
\begin{Def}
La (vraie) solution d'une équation différentielle du type~\eqref{ed}
$$\eps y'=y f(x,\eps)+h(x,\eps)+\eps y^2 P(x,\eps,y)\,,$$ 
dont le relief correspondant $\Re\Bigl(\int^x f(t,\eps)dt\Bigr)$ admet
un point col est appelée solution canard si ce point col n'est pas un
point tournant pour elle.
\end{Def}

On considèrera, dans la suite, uniquement les canards dans des
domaines ouverts du plan complexe. Ils seront bien entendu holomorphes
dans ces domaines. Sur la droite réelle, il est en revanche possible
de trouver des canards moins réguliers, dont seules les premières
dérivées restent bornées dans un intervalle autour du point col.

L'existence de solutions canards a déjà été démontrée dans des cas
très généraux. Parmi les résultats récents, on peut citer les deux
théorèmes suivants:
\begin{Th}\cite{FSs&r}
Soit l'équation
$$\eps u' = f(x) u +\eps P(x,u,\eps)$$ 
où les fonctions $f$ et $P$ vérifient les propriétés suivantes: 
\begin{enumerate} 
\item La fonction $f$ est analytique dans un voisinage complexe 
d'un intervalle réel $[a,b]$ autour de $0$, et est à valeurs réelles
sur cet intervalle.  
\item On a $xf(x)>0$ pour tout $x$ réel non nul de ce voisinage; de 
plus, il existe $\lambda>0$ et $p$ entier impair tels que
$f(x)=\lambda x^p \left(1+\gO(x)\right)$ au voisinage de $0$.
\item Le fonction $P$ est analytique dans un voisinage de $[a,b] 
\times \{0\} \times \{0\}$ dans $\CC^3$.
\end{enumerate} 
Alors on a l'équivalence suivante: il existe une solution formelle
$\hat u=\sum u_n(x)\eps^n$ dont tous les coefficients $u_n$ sont
analytiques au voisinage de $0$ si et seulement si il existe une
solution canard, tendant vers $u_0$ quand $\eps\rightarrow 0$,
uniformément dans un voisinage complexe de $x=0$.
\end{Th}
Ce premier théorème montre que dans le cas des solutions canards, il y
a en fait une relation presque automatique entre l'existence d'une
solution formelle et l'existence d'une vraie solution. Le cas des
équations avec paramètre est un peu différent, puisqu'on ne trouve,
formellement, un paramètre que sous la forme d'une série
formelle. Dans le cas de ces équations, on conna\^it le résultat
suivant:
\begin{Th}\cite{BFSW}
On considère l'équation avec le multiparamètre $a$ suivante, définie
pour $x$ dans un domaine $\Dom$: 
$$\eps u' =\Psi(x,u,a,\eps)\,,$$ 
où $\Psi$ est une fonction analytique dans un ouvert connexe $D_\Psi$ 
de $\CC\times\CC\times\CC^p\times\CC$. On suppose que 
\begin{enumerate} 
\item Il existe $a_0\in \CC^p$ et une fonction analytique $u_0$ 
tels que $\Psi(x,u_0(x),a_0,0)=0$ pour tout $x\in \Dom$ (il faut bien
s\^ur que si $x\in\Dom$, $(x,u_0(x),a_0,0)$ soit dans $D_\Psi$).
\item La fonction $f(x)=\Dp{\Psi}{u}\Bigl(x,u_0(x),a_0,0\Bigr)$
admet un unique zéro, $x_0$ dans $\Dom$, qui est d'ordre $p$. Cette
fonction $f$ est celle donnant le relief pour l'équation, avec la 
relation $R(x)=\Re\left(\int^x f(t) dt\right)$; $x_0$ est donc le seul 
point col pour ce relief dans $\Dom$.
\item Le domaine $\Dom\setminus\{x_0\}$ est la réunion de $p+1$
domaines accessibles, ces domaines étant éventuellement non bornés.
\item On note $\Psi_a:x\rightarrow \Psi(x,u_0(x),a,0)$. L'application
\begin{alignat*}{2} 
J:\ & \CC^p&\longrightarrow &\,\CC^p\\
    & a &\longrightarrow
    &\left(\Psi_a(x_0),\Psi'_a(x_0),\ldots,\Psi^{(p-1)}_a(x_0)\right)
\end{alignat*}
est un difféomorphisme local au voisinage de $a_0$.
\end{enumerate}
Alors l'équation différentielle admet une solution canard $u(x)$,
voisine de $u_0$, pour une valeur du paramètre $a$ tendant vers $a_0$
avec $\eps$.
\end{Th}

Ce dernier théorème d'existence sera redémontré au §~\ref{sec:alpha} 
(d'une autre manière que dans l'article cité).

\newpage
\begin{figure}[hp!]
\center{\includegraphics*[40mm,40mm][330mm,230mm]{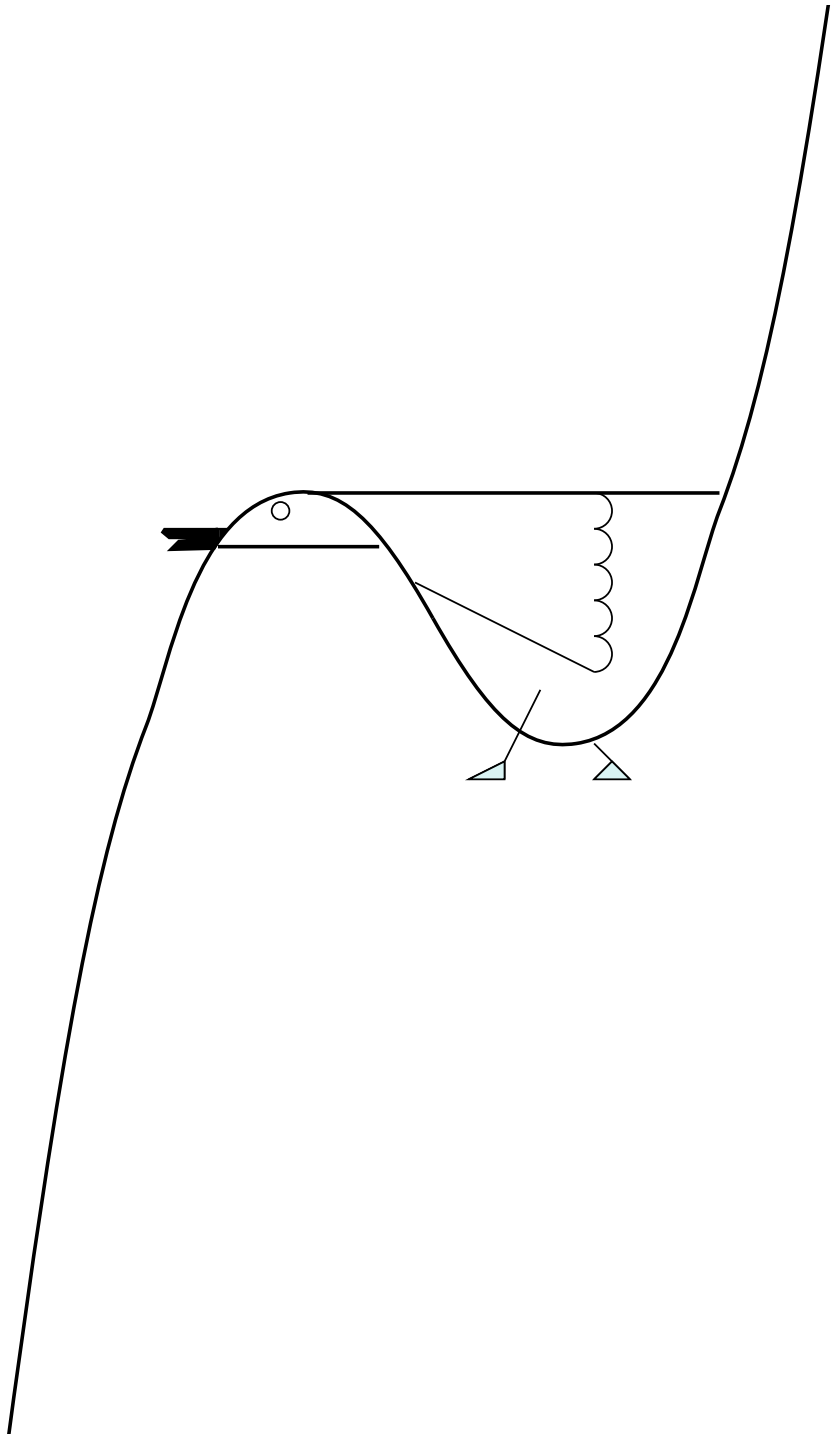}}
\caption{Canard, dit {\protect «}de Van der Pol\protect »}\label{fig:COIN}
\end{figure}

\chapter{Les canards de Van der Pol}

\label{chap:vdp}

\section{Introduction}

\begin{figure}[ht!]
\center{\includegraphics[scale=0.8]{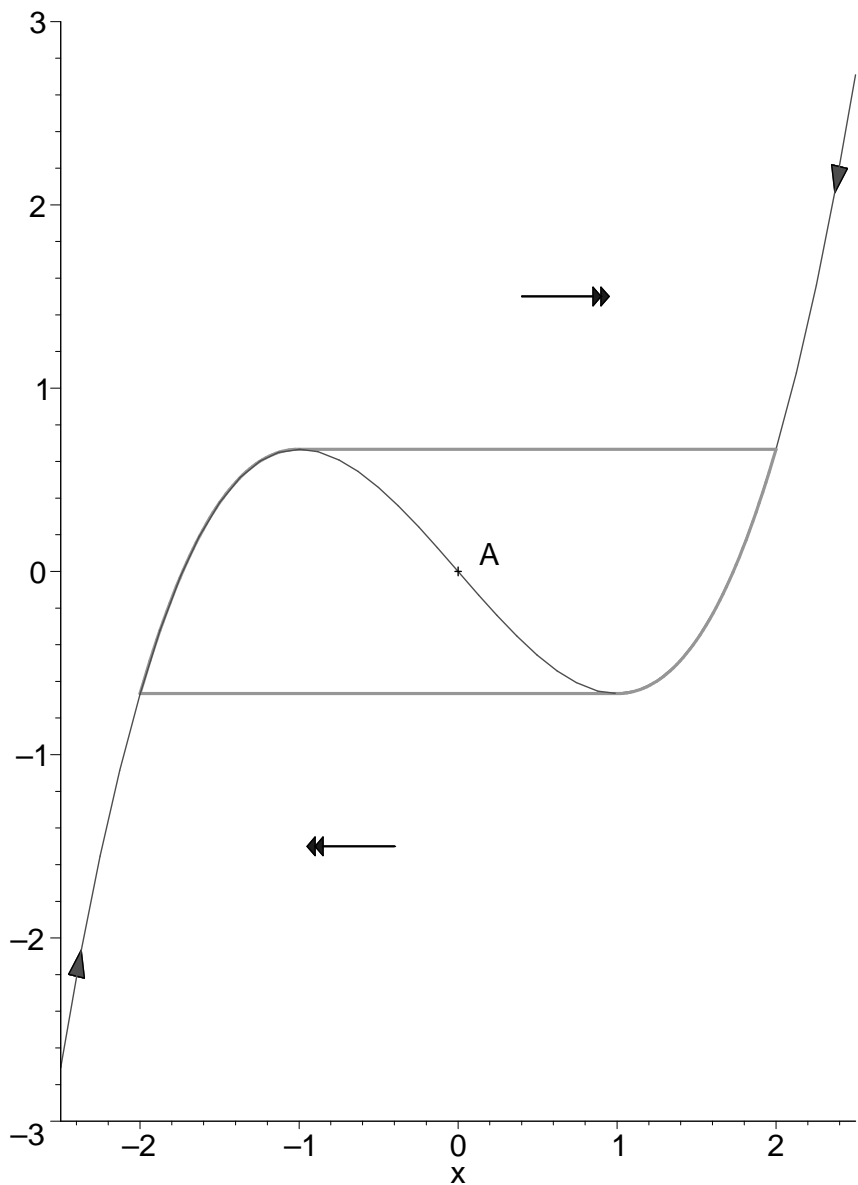}}
\caption{Cycle pour $\alpha=0$}
\label{fig:canard}
\end{figure}

L'équation de van der Pol \cite{vdP}
$$\eps \ddot{u} +(u^2-1)\dot u +u=0$$
est un exemple connu d'équation différentielle dont les solutions
présentent des oscillations rapides entre deux périodes de relaxation.

On considèrera ici l'équation de van der Pol \og forcée\fg
\par\penalty 500
$$\eps \ddot{u} +(u^2-1)\dot u +u=\alpha$$
\par\penalty 500
où $\dot u$ désigne la dérivée $\frac{du}{dt}$, où $\eps$ est un
paramètre arbitrairement petit et $\alpha$ un paramètre réel.

Dans toute la suite, on se limitera au cas $\alpha\geq 0$ (le cas
$\alpha$ négatif lui étant symétrique). 

Pour étudier les solutions de cette équation, on se place dans le plan
des phases (ou plan de Liénard,~\cite{Lie}), en prenant le système 
équivalent suivant:
\[ \left\{ \begin{array}{rl}
   \eps \dot u &=y-\frac{u^3}{3}+u\\
      \dot y   &=\alpha -u
           \end{array}
\right. \]

Une étude sommaire à partir de ce système montre que, pour $\alpha>1$,
la trajectoire rejoint rapidement le voisinage de la cubique
d'équation $y=u^3/3-u$, puis joint le point $A(\alpha,
\alpha^3/3-\alpha)$ qui est stationnaire.  Pour $\alpha<1$,
$\alpha$ pas trop proche de $1$, la trajectoire finit par contre en un
cycle comprenant deux parties stables de la cubique et deux segments
horizontaux (cf. figure~\ref{fig:canard}). Une étude locale au
voisinage de $A$ montre d'ailleurs qu'il y a une bifurcation de Hopf
au voisinage de $\alpha=1$: le point stationnaire $A$ est attractif si
$\alpha<1$, et répulsif si $\alpha>1$ et on a dans ce cas un cycle
limite dont la taille est d'ordre $\sqrt{1-\alpha}$ (pour les $\alpha$
très peu différents de $1$).

\medskip

Le phénomène dont il est question ici est cependant distinct, quoique
très voisin quand on fait varier $\alpha$, de celui de la
bifurcation. Il a été découvert quand E.~Beno{\^\i}t, J.L.~Callot,
F. et M.~ Diener ont cherché à voir, d'abord numériquement, comment on
passait (pour $\alpha<1$) d'un cycle de petite taille
$\sqrt{1-\alpha}$ au cycle beaucoup plus grand de la
figure~\ref{fig:canard}; d'après la propriété de dépendance continue
des solutions (et donc du cycle) par rapport au paramètre $\alpha$
(qui intervient de manière non singulière dans l'équation), on savait
qu'il devait exister un régime intermédiaire entre ces deux cycles.

\begin{figure}[!hp]
\center{\includegraphics[scale=0.5]{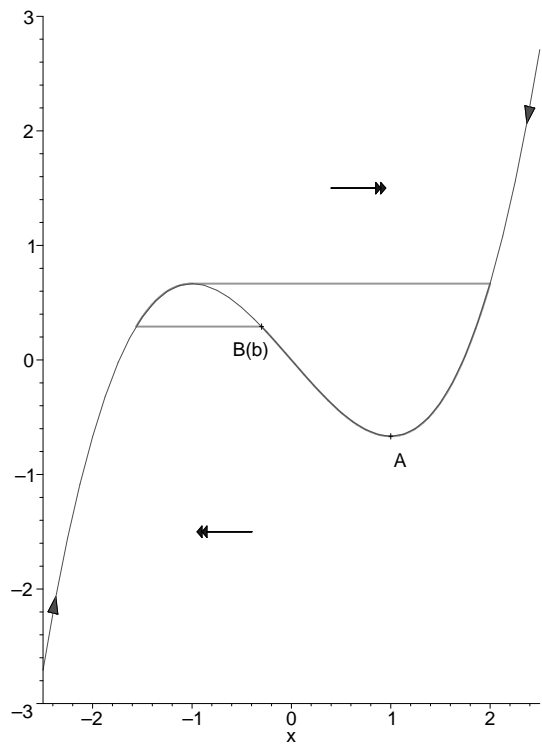}}
\caption{Canard intermédiaire; B, d'abscisse $b$, est appelé 
\protect «col\protect » du canard}
\label{fig:candebase}
\end{figure}
\begin{figure}[!hp]
\center{\includegraphics[scale=0.5]{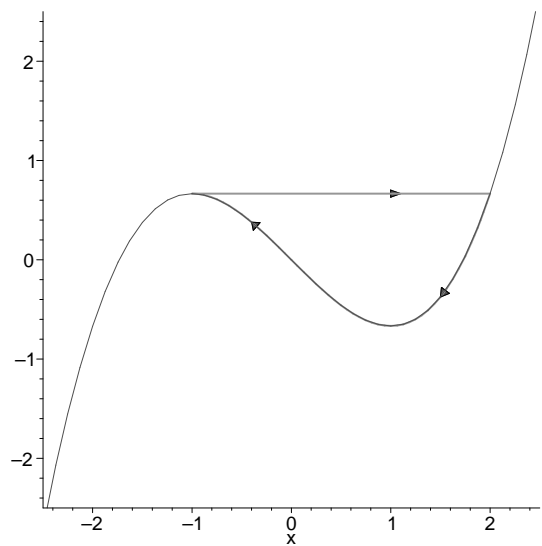}}
\caption{Grand canard}\label{fig:gdcanard}
\end{figure}

Les premiers résultats obtenus sur ce problème sont les suivants: 
pour $\eps$ arbitrairement petit fixé, on montre \cite{BCDD} que quand
$\alpha$ varie au voisinage de $1$, on peut obtenir un régime
intermédiaire entre les deux types de trajectoires cités, dans lequel
la trajectoire suit la cubique sur un morceau (aussi grand que l'on
veut) de la partie instable de cette courbe; on appelle les solutions
correspondantes des \og canards\fg. Plus précisément, 
\begin{description}
\item{1.} pour des valeurs $\breve\alpha$  voisines de 1 bien 
  choisies, la trajectoire suit la cubique entre ses deux sommets 
  $\bigl(1,-\frac{2}{3}\bigr)$, $\bigl(-1,\frac{2}{3}\bigr)$: on a 
  alors ce qu'on appelle des \og grands canards\fg 
  (cf. figure~\ref{fig:gdcanard});
\item{2.} la trajectoire suit la cubique jusqu'au point d'abscisse 
  $b\in\ ]-1,1[$ si et seulement si
  $|\alpha-\breve\alpha|= 
  \exp\bigl[-\frac{1}{\eps}\bigl(R(b)+o_\eps(1)\bigr)\bigr]$,
  où  $\breve\alpha$ est une des valeurs à grand canard, et 
  $R(b)=\int^b_1 (\xi-1)(\xi+1)^2 d\xi$ (voir
  figure~\ref{fig:candebase}, et la figure~\ref{fig:COIN} pour
  l'origine de l'appellation de \og canard\fg); 
\item{3.} si $\breve\alpha$ est une valeur à grand canard, la trajectoire
  passe aussi dans un voisinage de taille $\eps$ du point $(u=1,
  y=-2/3)$ et la solution du système correspondante est donc
  nécessairement continue en ce point (ou plut\^ot: bornée en $\eps$
  dans un voisinage de ce point).
\end{description}

\bigskip

Ce sont ces solutions particulières qui vont \^etre étudiées. Dans ce
but, on effectue le changement de variables suivant, pour se placer au
voisinage de la cubique: $ y=\frac{u^3}{3}-u+\eps v$, ce qui nous mène au
système 
\[ \left\{ \begin{array}{rl}
   \dot u      &=v\\
   \eps \dot v &=v(1-u^2)+\alpha -u
           \end{array}
\right. \]
On élimine la variable temporelle:
\begin{equation}
\eps v\frac{dv}{du}=(1-u^2)v+\alpha -u\,.\tag{E}\label{eq:vdpdebase}
\end{equation}

Comme courbe lente pour cette équation ({\it i.e.} pour le paramètre
$\eps=0$), on trouve la 
fonction $u\mapsto-\frac{\alpha-u}{1-u^2}$, qui n'est continue en $1$
que si $\alpha=1$. Dans ce cas, on trouve facilement la solution
formelle 
\begin{gather}
\hat \alpha =\sum_{n=0}^{\infty} a_n \eps^n \quad  \quad
\hat v(u) =\sum_{n=0}^{\infty} v_n(u) \eps^n\notag\\
\intertext{avec les formules de récurrence suivantes}
\left. \begin{array}{rl} 
\displaystyle{a_0 =1}     \quad & 
\quad\displaystyle{a_{n+1} =\sum_{j=0}^n v_j(1)v'_{n-j}(1)}\\
\displaystyle{v_0(u) =\frac{-1}{u+1}}  \quad & 
\quad\displaystyle{v_{n+1}(u) =v_0(u)\frac{\sum_{j=0}^n
  v_j(u)v'_{n-j}(u) -a_{n+1}}{u-1}} 
       \end{array}
\right\}\label{eq:serie}
\end{gather}

Une étude de cette série, réalisée par M.~Canalis-Durand avec les
outils de l'analyse Gevrey, a montré
qu'elle était divergente, mais Gevrey d'ordre 1 \cite{CD}.

\medskip
 
Dans \cite{FS}, en étudiant l'équation dans le champ complexe,
il est démontré que:
\begin{description}
\item{1.} pour $\eps\in\RR^+$, il existe une unique solution 
$(\alpha^+,v^+)$ tel que la solution $v^+$ existe en particulier, pour
n'importe quel $\rho$ réel positif, dans l'ensemble du secteur infini  
$S\bigl(-1+\rho,-\frac{\pi}{6},\frac{5\pi}{6}\bigr)= 
\left\{ z\in\CC\, /\ \arg(z+1-\rho)\in\ ]-\frac{\pi}{6},\frac{5\pi}{6} [
\,\right\}$; 
$v^+$ correspond dans son domaine de définition à la solution formelle
$\hat v$;
\item{2.} pour $\alpha^-=\overline{\alpha^+}$, on obtient aussi une vraie
solution $v^-$ existant dans un domaine symétrique de celui de $v^+$;
\item{3.} il existe une solution ($\alpha(\eps)$, $v(u,\eps)$) 
holomorphe en $\eps$ avec une singularité en $\eps=0$ qui prolonge la
solution $\bigl(\alpha^+,v^+\bigr)$ dans un large secteur centré en
$0$ (pour la variable $\eps$); cette solution vérifie en plus la
propriété: $\alpha^-(\eps)= \alpha\left(\eps e^{-2i\pi}\right)$,
$v^-(u,\eps)= v\left(\overline u, \eps e^{-2i\pi}\right)$.
\item{4.} pour tout $\eps \in \CC$ assez petit tel que
$|\arg(\eps)|<\frac{\pi}{2} +\delta$, on a
$$|\alpha(\eps e^{-2i\pi})-\alpha(\eps)|\leq
\exp\biggl(\Re\Bigl(\frac{R(-1+\delta')}{\eps}\Bigr)\biggr)$$ 
où $\delta$ et $\delta'$ sont des constantes positives, et
$R(x)=\int_1^x (t-1)(t+1)^2 dt$.
\end{description}
\bigskip

L'analyse à partir de cette série formelle ne permet pas cependant de
voir ce qui se passe au voisinage de $(-1)$, puisque la courbe lente y
a une singularité; alors que le dernier résultat ci-dessus laisse
penser qu'on gagnerait en précision si on pouvait remplacer $\delta'$
par quelquechose qui tend vers $0$ avec $\eps$.

Nous allons regarder spécifiquement ce qui se passe en ce point en
introduisant un nouveau changement de variable qui crée une loupe près
de $(-1)$:
$$u=-1+\eps^{1/3} X\quad,\qquad v=\eps^{-1/3} Y$$
d'où
\begin{equation}
\eqref{eq:vdpdebase}\Longleftrightarrow \qquad 
Y\frac{dY}{dX} = 2XY \bigl(1-\frac{\eps^{1/3}}{2}X\bigr) 
+\alpha+1-\eps^{1/3}X\tag{E'}\label{eq:vdpnonsing}
\end{equation}
La courbe lente $(\eps\rightarrow 0\Rightarrow\alpha\rightarrow 1)$
suivra l'équation
\begin{equation}
Y_0 \frac{dY_0}{dX} =2XY_0+2\label{eq:Y0}
\end{equation}

Nous allons étudier les solutions de \eqref{eq:Y0} bornées à l'infini
et montrer qu'il existe une vraie solution de l'équation
différentielle~\eqref{eq:vdpnonsing} qui admette un développement
asymptotique du type $Y(X,\eps)=\sum_{n=0}^\infty
Y_n(X)\eps^{n/3}$. Cette solution correspondra donc aussi à la série
formelle obtenue de \eqref{eq:serie} par le changement de variable;
l'intersection de son domaine d'existence avec le domaine de $z^+$
étant non nulle et comme elles sont égales dans cette intersection,
par unicité, $Y(X,\eps)$ sera le prolongement analytique de $z^+$
jusqu'au point $-1$.

L'existence de cette solution permettra de trouver un équivalent exact
de la différence $\alpha^+-\alpha^-$ (paragraphe \ref{sec:calculvdp}). À
l'aide de cet équivalent, nous en déduirons ensuite un équivalent pour
les coefficients $a_n$ de $\alpha$ (dans le paragraphe \ref{sec:an}).

\section{Étude de la courbe lente}

Il est possible de donner, d'une certaine manière, une solution exacte
de l'équation de la courbe lente $\eqref{eq:Y0}$ qui vérifie les
conditions aux limites voulues. Plus exactement, nous allons démontrer
le théorème suivant:
\begin{Th}
Il existe des solutions de l'équation différentielle
$$\eqref{eq:Y0}\ \ \qquad Y_0 \frac{dY_0}{dX} =2XY_0+2$$ qui tendent
vers $0$ en l'infini et dont les domaines d'existence contiennent des
secteurs ouverts d'ouverture au moins $\frac{2\pi}{3}$ en l'infini. En
particulier, il y en a une qui est holomorphe dans un secteur
$S\left(\infty,\varrho,-\f{\pi}{6},\f{5\pi}{6}\right)$, et une autre
dans un domaine symétrique contenant
$S\left(\infty,\varrho,-\f{5\pi}{6},\f{\pi}{6}\right)$. Ces deux
solutions peuvent \^etre prolongées jusqu'en $0$.
\end{Th}
\bigskip

Pour démontrer ce théorème, on commence par effectuer le changement de
variable suivant:
\begin{alignat}{2}
Y_0(X) &= X^2+z(X)\ \Longrightarrow\quad \eqref{eq:Y0}
&z'(X) &=\frac{2}{X^2+z(X)} \notag\\
\intertext{On regarde ensuite la fonction réciproque de la fonction
  $X\mapsto z(X)$ : on considère en fait $z\mapsto X(z)$; on sait, 
  d'après les théorèmes sur les
  fonctions implicites que si $z'$ ne s'annule pas, la dérivée de $X$ 
  s'écrira:}
&  & \frac{dX(z)}{dz} &=\frac{X^2(z)+z}{2}\label{eq:vdpx}\,.\\
\intertext{On pose ensuite}
X &=-2\frac{u'}{u}  \hbox{, \ \ d'où}       &u''+\frac{1}{4}uz&=0,\, 
\text{  qui est une équation d'Airy.} \notag 
\end{alignat}
Les solutions de l'équation d'Airy peuvent toutes s'écrire (par
exemple)  à partir des fonctions d'Airy $\Ai$ et ${\mathsf{Bi}}$.

On note 
$$\mu=\sqrt[3]{\frac{1}{4}},\qquad j=e^{2i\pi/3}.$$
Une des solutions correspondra par exemple à $u=\Ai(-\mu j^2 z)$:
$$x=2\mu j^2 \frac{\Ai'(-\mu j^2 z)}{\Ai(-\mu j^2 z)}$$

Or, les résultats connus (voir \cite{Ol}) sur $\Ai$ sont:
\begin{align}
 \begin{split}
\Ai(z) &=\frac{e^{-\xi}}{2\sqrt{\pi}z^{1/4}} \biggl( 1+ \gO\Bigl(
\frac{1}{\xi}\Bigr) \biggr)\\
\Ai'(z) &=\frac{-e^{-\xi}z^{1/4}}{2\sqrt{\pi}} \biggl( 1+ \gO\Bigl(
\frac{1}{\xi}\Bigr) \biggr)
 \end{split}\label{eq:Ai}
\intertext{où $\xi=\frac{2}{3}z^{3/2}$ (en prenant la valeur principale de la
racine), ceci pour} 
-\pi+\delta \leq& \arg z \leq \pi-\delta\quad(\delta>0)\,.\notag
\end{align}
Comme $\arg(-\mu j^2 z)=\arg(z)+\pi/3$, on s'intéresse aux $z$ tels
que $-\frac{4\pi}{3}+\delta \leq \arg z \leq \frac{2\pi}{3}+\delta$,
et on voit que dans ce cas
\begin{align*}
X^2(z) &= 4\frac{{u'}^2(z)}{{u}^2(z)} =4\mu^2 j^4 (-\mu j^2 z) \biggl( 1+ 
\gO\Bigl(\frac{1}{\xi}\Bigr) \biggr)\\
    &=-z+\gO\biggl(\frac{1}{\sqrt{|z|}}\biggr)\,.
\end{align*}
Cette dernière égalité permet d'affirmer, en utilisant le théorème
d'inversion locale, que la fonction $z\mapsto X(z)$ a un inverse quand
$z$ est dans un voisinage de l'infini dans le secteur proposé, soit la
fonction $X\mapsto z(X)$. 

Et alors $Y_0(X)=X^2+z(X)=\gO({1}/{\sqrt{|z|}})$ tend vers $0$ quand
$X\rightarrow\infty$ pour
$$-\frac{\pi}{6} +\delta' \leq \arg X \leq \frac{ 5\pi}{6} -\delta'$$

On obtient bien une solution $Y_0$ bornée sur un large secteur de
sommet $O$. D'après l'équation différentielle \eqref{eq:Y0}, $Y_0$ ne
pourrait tendre vers l'infini que quand $X$ tend vers l'infini. En
effet, supposons qu'il existe une suite de points $(X_n)$ convergeant
vers un point fini $X_\infty$, telle que $Y_0(X_n)$ soit une suite non
bornée. Alors la fonction $u_0=1/Y_0$ vérifie l'équation 
$$\f{du_0}{dX}=-2Xu_0^2-2u_0^3\,,$$ 
et elle est telle que $u_0(X_n)\rto 0$. Comme $u_0$ n'a pas de
singularité en $X_\infty$, et que d'après le théorème de Painlevé elle
vérifie $u_0(X_\infty)=0$, il s'agit de la l'unique solution de son
équation qui s'annule en $X_\infty$, donc de la fonction nulle. Ce qui
est absurde. Les seules singularités de $Y_0$ sont donc des zéros de 
$Y_0$.

L'équation différentielle nous indique aussi que
$$Y_0\underset{X\rightarrow\infty}{\sim} -\frac{1}{X}\,:$$ comme $Y_0$
est bornée dans un domaine comprenant un secteur en l'infini, d'après
la formule de Cauchy, $Y'_0$ est bornée elle aussi à l'intérieur de ce
domaine. Donc quand $X$ devient grand, on a nécessairement que
$2XY_0+2 \rto 0$. D'où le résultat.

Si on veut aller plus loin dans le développement de $Y_0$, il suffit
de décomposer la fonction $Y_0=-{1}/{X}+\tilde Y$. La fonction $\tilde
Y$ est solution (bornée quand $X$ devient grand dans le bon secteur)
de
$$\tilde {Y'}\tilde Y =2X  \tilde Y+\f{1}{X^3}-\f{\tilde Y}{X^2}+
\f{\tilde {Y'}}{X}\,.$$
$\tilde Y=o(1/X)$; avec Cauchy, on en déduit donc que $\tilde {Y'}$
est au plus un $o(1/X^2)$.  Dans l'équation différentielle, les deux
termes les plus importants sont cette fois $2X \tilde Y$ et
${1}/{X^3}$, donc $\tilde Y\sim {-1}/{2X^4}$.

Au total,
$$Y_0=-\frac{1}{X}-\f{1}{2X^4}+o\left(\f{1}{X^4}\right)\,.$$

\vskip 10pt

On peut montrer que la solution proposée  existe aussi sur l'ensemble
du demi-axe $\{ z\in \CC\ /\ \arg z=-\frac{\pi}{3}\}$.

En effet, pour $z=-jt\quad (t\in {\RR}^+)$,
\begin{gather*}
X(z)=2\mu j^2 \frac{\Ai'(\mu t)}{\Ai(\mu t)} \qquad 
\arg X(z)=\arg j^2+\pi-0=\frac{\pi}{3}\\
\frac{dX}{dt}=\frac{dX}{dz} \frac{dz}{dt}=\frac{X^2(z)+z}{2}(-j) =
-\frac{j^2}{2}\biggl(4\mu^2 \frac{{\Ai'}^2(\mu t)}{\Ai^2(\mu t)} -t
\biggr)\,. \\
\end{gather*}
Le calcul de la dérivée de $4\mu^2 {\Ai'}^2(\mu t)-t\Ai^2(\mu t)$, 
montre qu'elle est nulle; cette fonction est par conséquent
une fonction constante, et réelle, donc $\frac{-1}{j^2}
 \frac{dx}{dt}$ est une fonction réelle qui garde un signe constant
 (puisqu'on l'obtient de la précédente en divisant par $\Ai^2(\mu t)\in
 \RR^+$). 

D'où l'on peut déduire que $t\mapsto\frac{x(t)}{-j^2}$ est réelle et
uniforme sur ${\RR}^+$, donc inversible.

\vskip 10pt

On peut bien s\^ur prolonger analytiquement cette fonction sur tout
domaine simplement connexe ne contenant pas de singularité de $Y_0$
(les seules singularités possibles étant les points où $Y_0$
s'annule). 
Cela semble \^etre le cas par exemple (quand on fait des simulations
numériques) du quart de plan $(\Re z\geq 0,\,\Im
z\geq 0)$. On sait en tous cas que $Y_0$ existe sur une partie de
l'axe réel, à partir d'un réel $X_l\geq 0$ jusqu'au voisinage de l'infini.

\section{Existence d'une solution de \eqref{eq:vdpnonsing}, 
$Y(X,\eps)$, qui est voisine de $Y_0$}

Nous utiliserons un théorème de point fixe pour montrer que certaines
solutions de~\ref{eq:vdpnonsing} tendent vers $Y_0$ quand $\eps$ tend
vers $0$.
$$\eqref{eq:vdpnonsing}
\qquad YY'=2XY\Bigl(1-\frac{\eps'}{2}X\Bigr)+\alpha+1
-\eps'X,\qquad\text{avec }\eps'=\eps^{1/3}$$
On écrit $Y$ sous la forme $Y=Y_0+\eps'Z$. On rappelle \eqref{eq:Y0}:
$Y_0Y'_0=2XY_0+2$.

Alors
\begin{multline*}
Y_0Y'_0+\eps'Z'(Y_0+\eps' Z)+\eps'Z\left(2X+\frac{2}{Y_0}\right)\\
= 2XY_0-\eps'X^2 Y_0+\eps'ZX(2-\eps'X) +\alpha-1+2-\eps'X
\end{multline*}
\begin{align*}
\Longleftrightarrow\ Z' &= -X^2 +\frac{a/\eps'-X-2Z/Y_0}{Y_0+\eps'Z}\\
  &\text{ avec}  \quad a=\alpha-1=\gO(\eps)=o(\eps')\,.\\ 
\end{align*}

Après linéarisation de cette équation différentielle, on obtient
\begin{align}
Z'=&-X^2+u_0\left(\frac{a}{\eps'}-X\right) -Z\Bigl( u_0^2(2+a-\eps'X)
\Bigr) +\eps' h(Z,X,\eps')\label{eq:Zvdp}\\
\qquad&\text{où }\qquad u_0(X)=\frac{1}{Y_0(X)}
= -X+\f{1}{2X^2} +o\left(\f{1}{X^2}\right) \notag\\
\qquad&\text{et }\qquad h(Z,X,\eps')=Z^2 u_0^3\,
\frac{2+a-\eps'X}{1+\eps' u_0 Z}\,.\notag\\
\end{align}

On introduit les fonctions 
\begin{eqnarray*}
F(x) &=& \int^x -u_0^2(t) (2+a-\eps' t) dt,\\
\text{et } R(x) &=& \Re\bigl(F(x)\bigr)\,.
\end{eqnarray*}

On désignera par $\gamma_X$ un chemin dans le plan complexe partant
d'un point (éventuellement infini) à déterminer et allant jusqu'à $X$.
 
En intégrant l'équation différentielle avec la méthode de variation de
la constante, on voit que $Z$ satisfait à l'égalité:
$$Z(X)=\int_{\gamma_X} e^{F(X)-F(t)} \Biggl( -t^2+ u_0(t)
\left(\frac{a}{\eps'}-t\right) +\eps'h(Z,t,\eps') \Biggr) dt\,.$$

\vskip 10pt
On souhaite montrer qu'on peut appliquer un théorème de point fixe,
pour montrer qu'il existe une et une seule solution à l'équation
\eqref{eq:vdpnonsing} qui tend vers $Y_0$ quand $\eps$ tend vers $0$.
On se place pour cela dans un domaine ouvert $\DD$ (qui
sera précisé par après) et on considère deux solutions notées $Y_1$ et
$Y_2$ arbitrairement proches de $Y_0$. 

Soit alors $h_i=h(Y_i,t,\eps')$, pour $i=1,2$.

On trouve
$$h_1-h_2=\Bigl[u_0(t)\bigl(2+a-\eps't\bigr)\Bigr]
\Bigl[u_0(t)\bigl(Y_1-Y_2\bigr)\Bigr] 
\frac{(Y_1+Y_2)u_0+\eps'Y_1Y_2u_0^2}{(1+\eps'Y_1u_0)(1+\eps'Y_2u_0)}$$

On est donc amené à introduire les espaces fonctionnels suivants:
\begin{align*}
{\EuFrak H} &=\Bigl\{\text{fonctions holomorphes $h$ sur $\DD$}\Bigl/
\sup_{t\in \DD} \frac{|h(t)|}{|u_0(t)|\cdot|2+a-\eps't|} 
\text{ soit fini}\Bigr\}\\
{\EuFrak Y} &=\Bigl\{\text{fonctions holomorphes $Y$ sur $\DD$}\Bigl/
\sup_{t\in \DD} |Y(t)|\cdot|u_0(t)| \text{ soit fini}\Bigr\}\\
\intertext{que l'on munit des normes correspondantes 
$\| h \|_\EuH \text{ et } \| Y \|_\EuY$ respectivement.}
\end{align*}
\begin{alignat*}{5}
\text{Soit }\qquad {\EuFrak 0}_H :\quad&\tilde{\EuY}&
\rightarrow &{\EuH} &\quad\text{ et }\qquad&
{\EuFrak 0}_Y :\quad&\tilde{\EuH}&\rightarrow &&{\EuFrak Y}\,. \\ 
& Y &\mapsto & h(Y,t,\eps') &&& h &\mapsto && Z(X)
\end{alignat*}
Ces deux opérateurs sont bien définis au moins sur des sous-ensembles
$\tilde{\EuH}$ et $\tilde{\EuY}$ des espaces ${\EuH}$ et $\EuY$ (ceci
sera reprécisé plus loin).

Nous allons montrer que l'opérateur ${\EuFrak 0}_Yo {\EuFrak 0}_H$ est
une contraction sur ces ensembles.

$$\|h_1 -h_2\|_\EuH \leq \|Y_1 -Y_2\|_\EuY \frac{ \|Y_1\|_\EuY+
\|Y_2\|_\EuY+|\eps'|\cdot\|Y_1\|_\EuY \|Y_2\|_\EuY}{( 1 -|\eps'|
\cdot\|Y_1\|_\EuY)  ( 1 -|\eps'|\cdot\|Y_2\|_\EuY)}$$
donc si $\|Y_1\|_\EuY$ et $\|Y_2\|_\EuY$ sont majorées par
$\delta$, et pour $\eps'<1/\delta$, on a 
\begin{equation}
\|h_1 -h_2\|_\EuH \leq \|Y_1 -Y_2\|_\EuY \frac{2\delta
  +|\eps'|\delta^2}{ (1-|\eps'|\delta)^2 }\,.\label{ineg1}
\end{equation}

Par ailleurs, si $Y_k(X,\eps')=\EuO_Y(h_k(X,\eps'))$,
\begin{align*}
Y_1-Y_2 &=\int_{\gamma_X} e^{F(X)-F(t)} \eps'(h_1-h_2) dt\\
\intertext{d'où}
\|Y_1 -Y_2\|_\EuY\! &\leq \sup_{X\in \DD} |u_0(X)| |\eps'|
\int_{\gamma_X} e^{R(X)-R(t)} |h_1-h_2| dt\\
 &\leq  |\eps'|.\|h_1 -h_2\|_\EuH\sup_{X\in \DD} e^{R(X)}
\negthickspace\!\! \int_{\gamma_X}  \!\!\negthickspace e^{-R(t)}
\Bigl|u_0^2(t)\bigl(2+a-\eps't\bigr)\Bigr|  
\biggl|\frac{u_0(X)}{u_0(t)}\biggr| dt
\end{align*}

Choisissons un $X$ tel que: pour tout $X'$ tel que 
$|X'|\geq |X|$, $u_0(X')$ est peu
différent de $-X'$. Plus précisément, on va choisir $\gamma_X$ et $X$
tels que:
\begin{enumerate}
\item $\gamma_X(t) \in \DD$ pour tout $t$.
\item $\gamma_X$ va de l'infini à $X$ avec $\gamma_X(t)=e^{i\theta} t$,
où $t$ est un réel allant de $+\infty$ à $|X|$, et $\theta=\arg X$.
\end{enumerate}

Alors, pour $\theta$ dans l'intervalle $\Bigl]\frac{3\pi}{8}+\delta,
\frac{\pi}{2}-\delta\Bigr[$, 
\begin{align*}
F(X)-F(t) =
&\int^{|t|e^{i\theta}}_{|X|e^{i\theta}} -u_0^2(\phi) (2+a-\eps' \phi)
d\phi \qquad\quad t=|t|e^{i\theta}\\
     =&\int^{|t|}_{|X|} (2+a)(-\phi^2) d\phi +\eps'\int^{|t|}_{|X|} 
           \phi^3 d\phi +o(|t|^2) \\
     =&\frac{2+a}{3} e^{i(3\theta-\pi)} |t|^3 +\frac{\eps'}{4}
     e^{4i\theta} |t|^4  +o(|t|^2)
\end{align*}
La partie réelle de cette dernière expression est positive, (pour
$|t|$ assez grand donné indépendamment de $|\eps'|$, pour tout
$\arg(\eps')$ assez voisin de $0$) et croissante avec $|t|$. Donc si
on choisit un $X$ assez grand, il existe un $\mu>0$ tel que
\begin{equation}
\mu \Bigl|\frac{d}{dt} R\bigl(\gamma_X(t)\bigr)\Bigr| \geq  \Bigl|
F'\bigl(\gamma_X(t)\bigr)\Bigr|\,, \text{ pour tout $t$,}\label{eq:Haddock}
\end{equation}
d'où on déduit
\begin{align}
\|Y_1-Y_2\|_\EuY &\leq |\eps'|. \|h_1-h_2\|_\EuH \sup_{X\in \DD}
e^{R(X)} \int_{\gamma_X} e^{-R(t)} \mu \biggl| \frac{d}{dt}\notag
R\bigl(\gamma_X(t)\bigr)\biggr| dt \\ 
\|Y_1-Y_2\|_\EuY &\leq |\eps'|\cdot \mu \cdot \|h_1-h_2\|_\EuH\label{ineg2}
\end{align}

\penalty -500
Nous pouvons maintenant préciser le domaine $\DD$:
\begin{list}{$\bullet$}{}
\item 
$\DD \subset (\text{Domaine de }Y_0)$ (pour que $u_0$ soit définie).
\item 
Tout point de $X\in\DD$ peut être joint de l'infini par un chemin
$\gamma_X$ vérifiant~\eqref{eq:Haddock}.
\item
$(2+a-\eps't)$ ne s'annule pas pour $t$ dans ce domaine.
\end{list}

\vskip 10pt
Le domaine $\DD$ n'est pas vide: avec les conditions données, on peut
prendre par exemple pour $\DD$ un secteur
$S\left(\infty,2|X_f|,\f{3\pi}{8}+\delta ,\f{\pi}{2}-\delta\right)$.

\smallskip
Alors, pour $\eps$ assez petit, les inégalités $\eqref{ineg1}$ et
$\eqref{ineg2}$ signifient que ${\EuO}_Y o {\EuO}_H$ peut
être un opérateur contractant.

Il ne reste qu'à vérifier qu'on a bien ${\EuO}_Y({\EuH}) 
\subset {\EuY}$ et que  ${\EuO}_H({\Euy}) 
\subset {\EuH}$ où ${\Euy}$ est une boule ouverte de
${\EuY}$. 

Si $Y\in \EuY$ et $\|Y\|_{\EuY}<\f{1}{2\eps'}$ (on a là la
description de l'ensemble ${\Euy}$), alors
\begin{align*}
{\EuO}_H(Y)=h(Y,X,\eps')&=u_0(X) (Y(X) u_0(X))^2 
\f{2+a-\eps'X}{1+\eps' u_0(X) Y(X)}\\
\| h(Y,X,\eps')\|_\EuH &\leq \sup_{X\in \DD} 
\left|\f{(y u_0)^2}{1+\eps' u_0 y}\right|
\end{align*}
qui est bien borné, et à $Y$ fixé est aussi borné pour tout $\eps'$
assez petit, puisque $|u_0 Y|$ est borné pour $X\in \DD$ quand $Y\in
\EuY$. 

Si $h\in \EuH$, 
\begin{align*}
{\EuO}_Y(h)&=Z(h,X,\eps')\\
                   &=\int_{\gamma_X} e^{F(X)-F(t)} \Biggl( -t^2+ u_0(t)
\left(\frac{a}{\eps'}-t\right) +\eps'h(t,\eps') \Biggr) dt\\
                   &=Z(0,X,\eps')+
\left(Z(h,X,\eps')-Z(0,X,\eps')\right)
\intertext{La différence dans la parenthèse ne posant pas de problèmes 
vu ce qui a été montré ci-dessus (la fonction nulle est clairement
dans l'espace $\EuH$), on regarde uniquement le premier terme.}
Z(0,X,\eps')&=\int_{\gamma_X} e^{F(X)-F(t)} \left( -t^2+ u_0(t)
\left(\frac{a}{\eps'}-t\right) \right) dt
\intertext{On retrouve quelque chose d'analogue à la démonstration
concernant la différence des fonctions $Y_1-Y_2$.}  
\|Z(0,X,\eps')\|_\EuY\! &\leq \f{1}{|\eps'|}\sup_{X\in \DD} |u_0(X)| 
\int_{\gamma_X} e^{R(X)-R(t)} \left|-t^2 \eps'+u_0(t)
\left({a}-\eps't\right)\right| dt\\
 &\leq  \f{1}{|\eps'|}.\sup_{X\in \DD} e^{R(X)}
\negthickspace\!\! \int_{\gamma_X}  \!\!\negthickspace e^{-R(t)}
\left|-t^2 u_0\eps'+u_0^2(t)\left(a-\eps't\right)\right|  
\left|\frac{u_0(X)}{u_0(t)}\right| dt\\
\|Z(0,X,\eps')\|_\EuY\! &\leq  \|-u_0\eps't(u_0(t)+t)+a
u_0^2(t)\|_\EuH \sup_{X\in \DD} e^{R(X)} \int_{\gamma_X} e^{-R(t)} \mu
\biggl| \frac{d}{dt} R\bigl(\gamma_X(t)\bigr)\biggr| dt 
\intertext{La fonction $-u_0\eps't(u_0(t)+t)+a u_0^2(t)$ est bien dans
$\EuH$: $u_0(t)=-t+o(1/t)$, donc la fonction est majorée par
$\left|u_0(t)(1+|at|+|\eps'|)\right|<\left|u_0(t)(2+a-\eps't)\right|$
pour tout $\eps$ assez petit.}
\end{align*}
$$\|Z(0,X,\eps')\|_\EuY\ \text{ existe bien, }\
Z(0,X,\eps')\in\EuY\,.$$ 

\bigskip
On peut donc appliquer le théorème du point fixe, ce qui nous donne
l'existence et l'unicité de la solution de
\eqref{eq:Zvdp}, $Z(X,\eps')$, dans l'espace $\EuY$: cette fonction
est holomorphe et tend vers $0$ quand $X\rightarrow\infty$, pour tout
$\eps'$.

On en déduit:

\begin{Th}\label{ThVdPPtfixe}
Il existe une solution $Y_\eps$ de \eqref{eq:vdpnonsing} sur un
domaine $\DD$ contenant un secteur ouvert centré en l'infini (de
taille indépendante de $|\eps'|$) d'ouverture un peu inférieure à
$\pi/8$.  $Y_\eps$ tend vers $0$ en l'infini et elle est telle que
$Y_\eps\underset{\eps\rightarrow 0}{\longrightarrow} Y_0$
uniformément.  

Ce domaine $\DD$ peut même être encore étendu.
\end{Th}
\bigskip

Montrons qu'on peut, en effet, prendre un domaine $\DD$ plus grand que
celui proposé. Comme l'équation~\eqref{eq:vdpnonsing} n'est pas
singulière, les chemins $\gamma_X$ peuvent être rallongés d'une
distance arbitraire (mais finie) à partir de $X_f$ (par exemple, on
peut aller jusqu'à $0$, ou $X_l$); l'intégrale $$\int_{\gamma_X}
e^{R(x)-R(t)}|u_0^2(t)\bigl(2+a-\eps't\bigr)|
\biggl|\frac{u_0(X)}{u_0(t)}\biggr| dt$$ restera effectivement bornée
si on évite les singularités de $Y_0$, puisque toutes les fonctions
qu'elle contient sont holomorphes et donc bornées sur tout compact
inclus dans leurs domaines de définition. Et puisqu'il n'y a pas dans
cette intégrale de singularité en $\eps$, on pourra choisir des bornes
indépendantes de $\eps$.

\medskip

On peut même étendre $\DD$ sur une partie de l'axe réel, et
ce jusque $|\eps|^{-1/3}$ par exemple. En effet, supposons que l'on
ait déjà prolongé $\gamma_X$ de $X$ à $X_l$, puis éventuellement à 
$X_0 \in \RR^+$, $X_0$ assez grand.

Alors 
\begin{align*}
R(x) =&\, \text{Cste} + \int^x_{X_0} -t^2 (2+a-\eps't)dt +o(x^2)\\
\frac{d}{dt} R(t) =& -t^2 (2+a-\eps't) +o(t) \in \RR^-\ \text{ si }\ t\leq
\frac{1}{\eps'} \\ 
\intertext{(pour $\eps'$ assez petit, $a$ étant de l'ordre de
  ${\eps'}^3$)}
\|Y_1-Y_2\|_\EuY \leq&\, |\eps'|. \|h_1-h_2\|_\EuH \sup_x e^{R(x)}
\int_0^{{1}/{\eps'}} e^{-R(t)} \Bigl(-\frac{d}{dt}R(t)\Bigr) dt\\
                  \leq&\, |\eps'|. \|h_1-h_2\|_\EuH \times \text{Cste}
\end{align*}
ce qui implique une convergence uniforme sur $[X_l,\eps^{-1/3}]$,
puisque l'inégalité $\eqref{ineg1}$ reste bien sûr vérifiée tant que
l'on ne s'approche pas des singularités où des zéros de la fonction
$u_0$.

\vskip 10pt

Cette solution $Y_\eps$ de \eqref{eq:vdpnonsing} correspond à une
solution $v$ de \eqref{eq:vdpdebase} qui existe donc, pour $\eps>0$, au
moins jusqu'en $2|X_f|e^{i7\pi/16} \times \eps^{1/3}$, sur le segment
$[-1+X_l\eps^{1/3},0]$ et qui coïncide au voisinage de l'infini avec
la solution $y$ étudiée dans \cite{FS}: ce sont des prolongements
holomorphes l'une de l'autre.

Cette solution $Y_\eps$ admet un développement asymptotique uniforme:
\begin{Cor}\label{th:serievdp}
Il existe une suite de fonctions $Y_n(X)$ holomorphes en $X$ et
tendant vers $0$ quand $X$ tend vers l'infini telle que, pour tout
$N\in\NN$, et $X\in\DD$,
$$\left|Y_\eps(X,\eps')-\sum_{n=0}^N Y_n(x){\eps'}^n\right|\leq 
{\eps'}^{N+1}|R_N(X)|\,,$$
où la fonction bornée $R_N$ vérifie en outre
$$\lim_{X\rightarrow\infty} R_N(X)=0\,.$$
\end{Cor}

Le cas $N=0$ correspond au théorème précédent. On va démontrer le cas
général par récurrence, à partir du cas $N=1$. On reprend l'équation
\ref{eq:Zvdp}:
$$Z'=-X^2+u_0\left(\frac{a}{\eps'}-X\right) -Z\Bigl( u_0^2(2+a-\eps'X)
\Bigr) +\eps' h(Z,X,\eps')\,.$$
On a montré que $Z(X,\eps')$ était une fonction qui existait, dans
l'espace $\EuY$, et que sa norme $\|Z\|_{\EuY}$ était bornée 
indépendamment de $\eps'$.

On définit $Y_1$ comme solution de l'équation différentielle ci-dessus
avec $\eps'=0$: 
$$Y_1'(X)=-X^2-X u_0(X) -2 Y_1(X) u_0^2(X)\,,$$ bornée dans $\DD$.
Dans le secteur $[3\pi/8+\delta,\pi/2-\delta]$ considéré, le
coefficient du terme homogène, $-u_0^2(X)$, a une partie réelle qui
tend vers $+\infty$.  Il existe donc une et une solution $Y_1$ de
l'équation différentielle qui soit bornée dans le secteur
considéré. De plus, le terme non homogène est d'ordre au plus $-X^2-X
u_0(X)=o(1)$; on peut en déduire, comme pour $Y_0$, que $Y_1(X)$ tend
vers $0$ en l'infini.

On pose maintenant $Z(X,\eps')=Y_1(X)+\eps'Z_1(X,\eps')$.

Alors $Z_1(X,\eps')$ est solution de
\begin{align*}
\eps'Z'_1&=u_0\left(\frac{a}{\eps'}\right)+\eps'h(Z,X,\eps')
-(Y_1+\eps'Z_1)\Bigl( u_0^2(2+a-\eps'X)\Bigr) +2Y_1u_0^2\\
Z'_1    &=-Z_1(2+a-\eps'X) u_0^2+u_0\left(\frac{a}{{\eps'}^2}\right)
+u_0^2 Y_1\left(X-\f{a}{\eps'}\right)+h(Z,X,\eps')\,,
\end{align*}
ce qui est une équation différentielle linéaire en $Z_1$ (le fait que
$Z$ puisse s'écrire en fonction de $Z_1$ n'est pas gênant: on sait que
cette fonction $Z$ existe, et $h(Z,X,\eps')$ peut être considéré comme
une donnée, comme fonction de $\EuH$).

On peut résoudre cette équation:
$$Z_1(X,\eps')= e^{F(X,\eps')}\int_{\gamma_X}
e^{-F(t,\eps')}\left(u_0(t)\left(\frac{a}{{\eps'}^2}\right) 
+u_0^2(t) Y_1(t)\left(t-\f{a}{\eps'}\right)+h(Z,t,\eps')
\right)dt \,,$$
en rappelant que $F(X,\eps')=\int^X  -u_0^2(t) (2+a-\eps' t) dt$.

Alors, puisque $|u_0(t)(2+a-\eps' t)|>0$ dans $\DD$, et que les
fonctions $\f{h\bigl(Z(t,\eps'),t,\eps'\bigr)}{u_0(t)}\ $ et
${Y_1(t)}{u_0(t)}$ sont bornées, il est clair que $Z_1(X,\eps')$ est
une fonction bornée, qui reste dans l'espace $\EuY$. Ce qui démontre
le cas $N=1$.

Pour la suite de la démonstration par récurrence, elle devient
évidente à partir de l'équation différentielle obtenue pour $Z_1$: on
va trouver à chaque fois $Z_{n-1}(X,\eps')=Y_{n-1}(X)+\eps'Z_n$, puis
$$Z'_n=-Z_1(2+a-\eps'X) u_0^2+ h_n(X,\eps')\,,$$ avec
$h_n(X,\eps')=\f{h_{n-1}(X,\eps')-h_{n-1}(X,0)}{\eps'} +2u_0^2(X)
Y_{n-1}(X)$. $h_n$ est une fonction de $\EuH$ avec l'hypothèse de
récurrence qu'on peut exprimer comme ceci: $h_n\in\EuH$ et $Z_{n-1}\in
\EuY$. Donc, en résolvant l'équation, on voit à nouveau que
$Z_n(X,\eps')$ est une fonction de $\EuY$.

Ce qui suffit à démontrer le théorème.

\section{Estimation de $\alpha^+-\alpha^-$}
\label{sec:calculvdp}
\subsection{Préliminaires}
\subsubsection{Premiers termes du développement de $v(u,\eps)$ et $Y(X,\eps)$}
\par\penalty 2000
Le calcul des premiers termes des séries pour $a$ et $v$ 
(cf \eqref{eq:serie}) nous donne:
 \begin{align*}
a(\eps)=a_1\eps+a_2\eps^2+\gO(\eps^3) \quad&\quad 
v(u,\eps)=v_0(u)+\eps v_1(u)+\eps^2 v_2(u) +\gO(\eps^3)\\ 
\intertext{où}
a_1=-\frac{1}{8} \quad&\quad a_2=-\frac{3}{32}\\
\intertext{et}
v_0(u) =\frac{-1}{u+1} \quad&\quad 
v_1(u)=-\frac{1}{8}\frac{u^2+4u+7}{(u+1)^4}\\
 \end{align*}
$$v_2(u) = - \frac {1}{32} \frac {3u^{5} + 21u^{4} + 66u^{3} + 126u^{2} +
159u + 121}{(u + 1)^{7}}$$ 
Par récurrence, on démontre que pour $n>0$ les $v_n(u)$ sont des
fractions rationnelles, dont le numérateur est un polynôme de degré
$3n-1$, et dont le dénominateur est $(u+1)^{3n+1}$.

Le simple changement de variable $u=-1+\eps^{1/3} X$ dans le
développement de $v(u,\eps)$ nous donne alors, après simplification:
\begin{multline*}
\eps^{1/3} v=-\frac{1}{X} -\frac{1}{2X^4} -\frac{5}{4X^7}
-\eps^{1/3}\left(\frac{1}{4X^3} +\frac{9}{8X^6}\right)  
-\eps ^{2/3}\left(\frac{1}{8X^2} +\frac{3}{4X^5}\right)\\ 
-\eps \frac{3}{8X^4} -\eps ^{4/3}\frac{3}{16X^3} 
-\eps^{5/3}\frac{3}{32X^2} +\gO\biggl(\frac{1}{|X|^8}+|\eps|\frac{1}{|X|^5}
+|\eps|^2\frac{1}{|X|^2}\biggr)\,. 
\end{multline*}

Ensuite, comme $Y(X,\eps)=\eps^{1/3}v$ se met sous la forme
$$Y=Y_0+\eps^{1/3}Y_1+\cdots+Y_n\eps^{n/3}+o\left(\eps^{n/3}\right),$$ 
après identification, on trouve pour les $Y_n$:
\begin{align}
 \begin{split}
Y_0&=-\frac{1}{X}-\frac{1}{2X^4}+\gO\left(\frac{1}{X^7}\right) \\
Y_1&=-\frac{1}{4X^3}+\gO\left(\frac{1}{X^6}\right) \\
Y_2&=-\frac{1}{8X^2}+\gO\left(\frac{1}{X^5}\right) \\
Y_3&=\gO\left(\frac{1}{X^4}\right) \\
&\text{et pour $n\geq 4$, on a } \\
Y_n&=\gO\left(\frac{1}{X^2}\right)\text{ ou } Y_n=o\left(\frac{1}{X^2}\right)\,.
 \end{split}\label{eq:SerY}
\end{align}

\bigskip
\subsubsection{Équivalent de $Y_0^+-Y_0^-$}
\par\penalty 2000
On introduit les fonctions suivantes:
$$\Ai_0=\Ai,\ \ \Ai_1:x\mapsto \Ai(jx),\ \ 
\Ai_2:x\mapsto\Ai(j^2x)\,.$$ 
On cherche à estimer la différence de deux des solutions de
\eqref{eq:Y0}. Pour cela, nous nous intéressons aux deux fonctions
suivantes: 
$$X_1(z)=-2\frac{\Ai'_1(z)}{\Ai_1(z)}\qquad\text{et}\qquad X_2(z)=-2
\frac{\Ai'_2(z)}{\Ai_2(z)}$$ 

Ces deux fonctions vérifient, à cause de la symétrie par rapport à
l'axe réel de l'équation \eqref{eq:vdpx}, la
relation $X_1(z)=\overline{X_2(\overline{z})}$. Elles correspondent à
deux solutions de l'équation \eqref{eq:Y0}:
$$Y_0^+(X)=X_2^2\left(z(X)\right)+z(X)\qquad\text{et}
\qquad Y_0^-(X)=X_1^2\left(z(X)\right)+z(X)\,.$$

En particulier, pour $X_1(z)$ sur l'axe réel, on a:
\begin{align*}
X:=&\, X_1(z)= X_1\Bigl(Y_0^-\bigl(X_1(z)\bigr)-X_1^2(z)\Bigr) 
=X_1\Bigl(Y_0^-(X)-X^2\Bigr)\\
 =&\, X_2(\overline{z})= X_2\Bigl(\overline{Y_0^-(X)-X^2}\Bigr) =
 X_2\Bigl(Y_0^+(X)-X^2\Bigr).\\ 
\intertext{On note dans ce cas}
z^+(X)=&\,Y_0^+(X)-X^2\qquad\text{et}\qquad z^-(X)=Y_0^-(X)-X^2
\end{align*}

On cherche un équivalent de $Y_0^+-Y_0^-$, et nous allons montrer le
résultat suivant:
\begin{Lemme}\label{lemme:Stokes-vdp}
Pour $X$ sur l'axe réel,
$$(Y_0^+-Y_0^-)(X) \underset{X\rightarrow+\infty}\sim  
\frac{4}{e} X^2 e^{-\frac{2}{3} X^3} i\,.$$
\end{Lemme}
Ce qui revient à réussir à calculer exactement un coefficient de
Stokes pour une équation différentielle non linéaire; dans le cas 
général, ce n'est pas possible.
\bigskip

Pour obtenir ce résultat, on utilise le fait que pour $X$ réel et
suffisamment grand, $Y_0^\pm -X^2$ est presque un réel négatif. Dans
ce cas $|\Ai_k(Y_0^\pm -X^2)|$ est grand pour $k=1$ et $k=2$,
exponentiellement petit par contre pour $k=0$, ce qui justifie
l'utilisation d'un développement limité pour estimer les deux
fonctions $\Ai_2(z)=-j^2 \Ai_1(z) -j\Ai_0(z)$ et $\Ai'_2(z)=-j^2
\Ai'_1(z)-j\Ai'_0(z)$. Dans le calcul qui suit, on travaillera en
variable $X$, avec aussi les fonctions $z^\pm(X)$ introduites plus
haut.

\begin{align*}
 X_1\left(Y_0^-(X)-X^2\right) &= X_2\left(Y_0^+(X)-X^2\right)\\
 -2 \frac{\Ai'_1(z^-)}{\Ai_1(z^-)} &= -2 \frac{\Ai'_2(z^+)}{\Ai_2(z^+)} \\
 -2 \frac{\Ai'_1(z^-)}{\Ai_1(z^-)} &= -2 \frac{\Ai'_1(z^+)+j^2
   \Ai'_0(z^+)}{\Ai_1(z^+)+j^2 \Ai _0(z^+)}\\  
 X_1(z^-) &= X_1(z^+) \biggl( 1+ j^2 \frac{\Ai'_0(z^+)}{\Ai'_1(z^+)}
 -j^2  \frac{\Ai_0(z^+)}{\Ai_1(z^+)} +o(\ldots) \biggr)\\
 X_1(z^-)-X_1(z^+) &=  X_1(z^+) \biggl(j^2 \frac{\Ai'_0(z^+)}{\Ai'_1(z^+)}
 -j^2  \frac{\Ai_0(z^+)}{\Ai_1(z^+)} +o(\ldots) \biggr)
\intertext{Il va être possible de calculer un équivalent de
cette différence, ce qui nous mènera à un équivalent pour
$Y_0^+-Y_0^-$. En effet, avec une formule de Taylor}
 X_1(z^-)-X_1(z^+) \sim & (z^--z^+)X'_1(z^-) =
 (Y_0^--Y_0^+)X'_1(z^-)
\intertext{sachant que d'après l'équation différentielle pour $X_1$,
on peut écrire sa dérivée} 
X'_1(z^-) &= \frac{X_1^2(z^-)+z^-}{2}\\ 
X'_1\Bigl(z^-(X)\Bigr) &=
\frac{X^2+(Y_0^-(X)-X^2)}{2} \sim -\frac{1}{2X}\,.
\end{align*}
\medskip

Nous allons utiliser les équivalents déjà vus précédemment
(cf~\eqref{eq:Ai}) pour estimer les deux fractions
$\f{\Ai'_0(z^+)}{\Ai'_1(z^+)}$ et $\f{\Ai_0(z^+)}{\Ai_1(z^+)}$,
puisque $z^\pm$ est proche du demi-axe réel négatif quand
$X\rightarrow\infty$:
\begin{align}
\xi =\frac{2}{3}\left( \mu\left(X^2-Y_0(X)\right)\right)^{3/2} &=
\frac{2}{3}\bigl(2^{-2/3}\bigr)^{3/2}
\biggl(X^2+\frac{1}{X}+\gO\Bigl(\frac{1}{X^4}\Bigr)\biggr)^{3/2}
\notag\\
&= \frac{1}{3} \Bigl(X^3+\frac{3}{2} +o(1)\Bigr)\notag\\
\xi_2=\frac{2}{3}\left( j\mu\left(X^2-Y_0(X)\right)\right)^{3/2} &= 
e^{i\pi}\xi =-\xi\notag
\end{align}
\begin{align} 
\frac{\Ai'_0(z^+)}{\Ai'_1(z^+)}-\frac{\Ai_0(z^+)}{\Ai_1(z^+)} &\sim 
\frac{-e^{-\xi}z^{1/4}}{1} \times \frac{1}{-j(jz)^{1/4}e^\xi} -
\frac{e^{-\xi}}{z^{1/4}} \times  \frac{(jz)^{1/4}}{e^\xi}  \notag\\
        &\sim  e^{-5i\pi/6} e^{-2\xi} +  e^{-5i\pi/6} e^{-2\xi} \notag\\
        &\sim  2e^{-5i\pi/6} e^{-2\xi} \notag
\intertext{En rappelant que}
(Y_0^+-Y_0^-)X'_1(z^-) &\sim  -X_1(z^+) j^2 \left(
\f{\Ai'_0(z^+)}{\Ai'_1(z^+)}- \f{\Ai_0(z^+)}{\Ai_1(z^+)} \right)\notag \\
\intertext{on trouve comme équivalent pour la différence le résultat
  annoncé au lemme~\ref{lemme:Stokes-vdp}:} 
(Y_0^+-Y_0^-)(X) &\sim  \frac{4}{e} X^2 e^{-\frac{2}{3} X^3} i\qquad
\text{ quand $X\rightarrow\infty$ dans $\RR$.}\label{eq:dY0}
\end{align}

\subsection{Établissement de l'équation vérifiée par  $\alpha^+-\alpha^-$}

On rappelle l'équation \eqref{eq:vdpdebase}, vérifiée pour 
$(v,\alpha)=(v^+, \alpha^+)$ et $(v^-,\alpha^-)$.
\begin{equation*}
\eqref{eq:vdpdebase}\qquad\eps\frac{dv}{du}=(1-u^2)+\frac{\alpha -u}{v}\,.
\end{equation*}

\medskip
Le but des paragraphes à venir va \^etre d'obtenir une estimation de $
\alpha^+ -\alpha^-$. Pour cela, on commence par chercher une équation
différentielle en $ v^+-v^-$, où le terme à calculer sera présent.  On
pose alors
\begin{alignat*}{5}
b& =&\,\alpha^+ -\alpha^-&\qquad & c& =& \,\alpha^+ +\alpha^-\\
d& =&\, v^+-v^-           &\qquad & e& =& \, v^+ +v^-
\end{alignat*}

Quand on effectue la substitution dans l'équation
\eqref{eq:vdpdebase}, on obtient
\begin{align*}
\eps d' &=(1-u^2)+\frac{\alpha^+-u}{v^+}-(1-u^2)-\frac{\alpha^--u}{v^-}\\
            &=\frac{ud+\alpha^+v^--\alpha^-v^+}{v^+v^-}                \\
            &=\frac{u}{v^+v^-}d+ \frac{\frac{b+c}{2}
               \cdot\frac{e-d}{2}+ \frac{b-c}{2}\cdot\frac{e+d}{2}}{v^+v^-}\\
            &=\frac{u}{v^+v^-}d+\frac{e}{2v^+v^-}b-\frac{c}{2v^+v^-}d  \\
            &=\frac{u-c/2}{v^+v^-}d+\frac{v^+ +v^-}{2v^+v^-}b          \\   
\end{align*}

En remplaçant $c=\alpha^+ +\alpha^-=1+a^+ +1+a^-$, on obtient en
définitive l'équation voulue.
\begin{Lemme}
On pose $d(u)=v^+(u)-v^-(u)$ et $b=\alpha^+-\alpha^-$. Alors
 l'équation différentielle suivante est vérifiée: 
 \begin{gather*}
\eps d'(u)=f(u)d(u)+g(u)b\\
\intertext{avec}
f(u)=\frac{u-1-a}{v^+v^-}\quad \hbox{et}\quad 
 g(u)=\frac{v^+ +v^-}{2v^+v^-}\\
\intertext{sachant que} 
a=\frac{a^++a^-}{2}=-\frac{1}{8}\,\eps+ \gO (\eps^2)\,.\\
 \end{gather*}
(toutes les fonctions ici dépendent aussi de $\eps$, cette dépendance
 n'est pas rappelée). 
\end{Lemme}
\medskip

On intègre ensuite cette équation différentielle avec la méthode de la
variation de la constante, entre $+\infty$ et $x$ (car $\lim_{+\infty}d=0$):
\begin{gather*}
d(x)=\frac{b}{\eps}e^{F(x)/\eps}\int_{+\infty}^x e^{-F(u)/\eps}
g(u) du,\\
\intertext{où}
F(t)=\int_{1}^t f(u) du\,.\\
\end{gather*}

Le choix de $x=x_l=-1+X_l\eps^{1/3}$ permet de trouver que
\begin{equation}
b=\eps d(x_l) e^{-F(x_l)/\eps} \frac{1}{\int_{+\infty}^{x_l}
  e^{-F(u)/\eps}g(u)du}\,.  \label{eq:b}
\end{equation}

\subsection{Calcul de l'équivalent pour $b$} 

On peut trouver une estimation, quand $\eps\rightarrow 0$, de chacun
des facteurs dans le produit de l'égalité \eqref{eq:b}; le but final
étant d'en obtenir une pour $b$.

Tout d'abord, on a vu que jusqu'au point $u=x_l$, on pouvait affirmer
que la fonction $v^+(u)=\eps^{-1/3}\:Y^+(\eps^{-1/3}(1+u))$ était
équivalente  à $\eps^{-1/3}\:Y_0^+(\eps^{-1/3}(1+u))$  quand
$\eps\rightarrow 0$. Ce qui nous mène à
\begin{equation}
 d(x_l)=v^+(x_l)-v^-(x_l) \underset{\eps\rightarrow 0}{\sim} 
 \eps^{-1/3}\: \Bigl(Y_0^+(X_l)-Y_0^-(X_l)\Bigr)\,. \label{eq:b2}
\end{equation}

\vskip 20pt
Ensuite, l'intégrale présente dans \eqref{eq:b} est susceptible
d'\^etre estimée avec la méthode du point col. On constate d'abord que
$$F'(t)=f(t)=0 \Longrightarrow\ t=1+a =\alpha\,.$$
Or on vérifie sans difficulté que $F(\alpha)=\gO(\eps^2)$.
D'où 
$$\frac{F(t)}{\eps}= \gO (\eps) +\frac{(t-\alpha)^2}{2\eps}f'(\alpha)
+ \frac{\gO(t-\alpha)^3}{\eps}\,.$$ 
Sachant que $v^\pm(\alpha)=-\frac{1}{2}  +\gO (\eps)$, on trouve  
\begin{align*}
f'(\alpha) &=\frac{1}{v^+(\alpha)v^-(\alpha)} 
     -(v^+v^-)'(\alpha)\frac{\alpha-1-a}{(v^+(\alpha)v^-(\alpha))^2}\\
           &=\left(\frac{1}{(-1/2)^2}+\gO (\eps)\right)+0\\
           &=4+\gO (\eps),
\intertext{donc}
\frac{F(t)}{\eps} &=\frac{2(t-\alpha)^2}{\eps} +
\frac{\gO\Bigl(\bigl|t-\alpha\bigr|^3
  +|\eps|.\bigl|t-\alpha\bigr|^2\Bigr)}{\eps} \\
\intertext{Par ailleurs}
g(\alpha)  &\sim -2\,.\\
\end{align*}
\begin{align*}
\intertext{La contribution principale dans l'intégrale vient du
  voisinage de $\alpha$,  on va donc écrire, après le changement de
  variable $x=(t-\alpha)\sqrt{2/\eps}$,}
-\int_{x_l}^{+\infty}\! g(t) e^{-F(t)/\eps} dt 
 &\sim -\sqrt{\f{\eps}2} \int_{(x_l-\alpha)\sqrt{2/\eps}}^{+\infty}\! 
g\left(\alpha+x\sqrt{\eps/2}\right) 
e^{-F\left(\alpha+x\sqrt{\eps/2}\right)/\eps} dx\\
 & \sim -\sqrt{\f{\eps}2}\ g(\alpha)
\int^{+\infty}_{(-2)\sqrt{2/\eps}} e^{-x^2 +\sqrt{\eps}.
  \gO\bigl(|x|^3+\sqrt{|\eps|}.|x|^2\bigr)}\; dx\\
 &\sim 2\sqrt{\frac{\eps}{2}} \int^{+\infty}_{-\infty} e^{-x^2} dx\,.
\end{align*}
Un équivalent de l'intégrale est donc:
\begin{equation}
 -\int_{x_l}^{+\infty}\! g(t) e^{-F(t)/\eps} dt \sim
 \sqrt{2{\pi\eps}}\,. \label{eq:b3}
\end{equation}

\vskip 20pt
Il reste à trouver une estimation pour $F(x_l)$. Pour calculer
$F(x_l)$, on peut couper l'intervalle $[x_l,1]$ en deux parties, en un
point $(-1+\lambda)$, où $\lambda$ est un réel fixé à choisir entre
$0$ et $2$.  À gauche, on pourra alors raisonnablement estimer $v(u)$
par l'expression $\eps^{-1/3}\bigl(Y_0(X) +\eps^{1/3}\: Y_1(X)
+\eps^{2/3}\: Y_2(X) +\gO(\eps)\bigr)$ d'après le corollaire
\ref{th:serievdp} , et à droite par $v_0+ \eps v_1+\gO(\eps^2)$.

\begin{equation*}
\frac{1}{\eps}F(x_l) =\underbrace{\frac{1}{\eps}
  \int^{-1+\lambda}_1 \frac{u-1-a}{v^+(u)v^-(u)} du}_{(A)} +
\underbrace{\frac{1}{\eps} 
 \int^{x_l}_{-1+\lambda} \frac{u-1-a}{v^+(u)v^-(u)} du}_{(B)}
\end{equation*}

Les calculs des premières fonctions $v_n$ et les approximations pour
les $Y_n$ (cf.~\eqref{eq:SerY}) données dans les préliminaires permettent 
de donner une estimation  des deux intégrales.
\bigskip
Pour l'expression $(A)$, on part de
$$\frac{1}{v^+(u)v^-(u)}= (u+1)^2 \biggl( 1- \frac{\eps}{4}
\frac{u^2+4u+7}{(u+1)^3} \biggr)+\gO(\eps^2)\,$$
{puisque $v^+$ et $v^-$ ont le même développement asymptotique, d'où}
\begin{align*}(A)&=\frac{1}{\eps} \int^{-1+\lambda}_1\!\!
    \frac{u-1-a}{v^+(u)v^-(u)} du\\ 
   &= \frac{1}{\eps}\!\! \int^{-1+\lambda}_1 (u-1-a)(u+1)^2
    du -\frac{1}{4} \int^{-1+\lambda}_1 (u-1-a)
    \frac{u^2+4u+7}{u+1} du +\gO(\eps)\,. \\
\intertext{On rappelle que $a=-\eps/8 +\gO(\eps^2)$. Le calcul des intégrales
  nous donne alors}
(A)&= \frac{1}{\eps} \Bigl(\frac{4}{3} -\frac{2\lambda^3}{3}
    +\frac{\lambda^4}{4}\Bigr)+ \Bigl( \frac{1}{3}-
    \frac{\lambda^3}{24} +2\ln \lambda-2\ln(2) \Bigr)
    +\gO(\eps)\,. 
\end{align*}

\vskip 20pt
Pour $(B)$, le calcul se révèle plus complexe. On sait que chacune des
deux fonctions $v^+$ et $v^-$ se prolonge à l'aide d'une fonction du
type $Y_\eps$ du théorème \ref{ThVdPPtfixe}: 
\begin{align*}
v^\pm(u,\eps)&=\eps^{-1/3}\:Y_\eps^\pm\left(\f{u}{\eps'}\right)
\intertext{avec $X=u/\eps'$, en utilisant le corollaire
\ref{th:serievdp}, }
v^\pm(u,\eps)&=\eps^{-1/3}\left(Y_0^\pm(X)+\eps^{1/3}\: Y_1^\pm(X)
             +\eps^{2/3}\: Y_2^\pm(X) +\eps Y_3^\pm(X)+ 
             \f{\eps^{4/3}}{X}\gO_\eps(1)\right) 
\end{align*}
où le $\gO_\eps(1)$ est borné quand $\eps\rightarrow0$ et au moins
borné quand $X\rightarrow\infty$. 

Comme $v^+$ et $v^-$ ont le même développement asymptotique pour tout
$u>\rho$ dans un certain secteur autour de l'axe réel, il est clair
que le développement en $1/X$ des $Y_n$ donné plus haut est valable
pour $Y_n^+$ comme pour $Y_n^-$. Ces deux fonctions sont donc
exponentiellement proches en $1/X$ quand $X$ devient grand. On peut
donc écrire,
\begin{align*}
(B)&= \frac{1}{\eps} \int^{x_l}_{-1+\lambda}
    \frac{u-1-a}{v^+(u)v^-(u)} du \\ 
   &=\int^{X_l}_{\lambda\eps^{-1/3}}\!
    \frac{-2+\eps^{1/3}X-a}{Y_0^+(X)Y_0^-(X)} \Bigl( 1 -\!
    \eps^{1/3}Q_1(X) + {\eps}^{2/3} Q_2(X) + \gR(X,\eps) \Bigr) dX\,, 
\intertext{en notant}
Q_1&=\f{Y_1^+}{Y_0^+}+\f{Y_1^-}{Y_0^-}\,,\ \ 
Q_2=\f{{Y_1^+}^2}{{Y_0^+}^2}+\f{Y_1^+Y_1^-}{Y_0^+Y_0^-}
    +\f{{Y_1^-}^2}{{Y_0^-}^2}-\f{Y_2^+}{Y_0^+}-
    \f{Y_2^-}{Y_0^-}\,; 
\intertext{le reste $\gR(X,\eps)$ se trouve être de la forme 
$\left(\f{\eps}{X}+{\eps^{4/3}}\right)\gO_\eps(1)$, puisqu'il
correspond à des termes $\eps \f{Y_k^\pm}{Y_0^\pm}$ ($1\leq k\leq 3$)
et $\f{\eps^{4/3}}X\f1{Y_0^\pm}$ respectivement équivalents quand $X$
devient grand à $\f{\eps}{X}$ et ${\eps^{4/3}}$.}
(B) &= \int^{X_l}_{\lambda\eps^{-1/3}}\!\!\!
     \f{-2+\eps^{1/3}X-a}{Y_0^+(X)Y_0^-(X)} dX -
     \eps^{1/3}\int^{X_l}_{\lambda\eps^{-1/3}}\!\!\!
     \f{-2+\eps^{1/3}X-a}{Y_0^+Y_0^-}\:Q_1\ dX\\
    &\phantom{=\,}+ {\eps}^{2/3}\int^{X_l}_{\lambda\eps^{-1/3}}\!\!\!
     \f{-2+\eps^{1/3}X-a}{Y_0^+Y_0^-}\:Q_2\ dX 
     +\int^{X_l}_{\lambda\eps^{-1/3}}\!\!
     \frac{(-2+\eps^{1/3}X-a)\gR(X,\eps)}{Y_0^+Y_0^-}\:dX
\end{align*}
On appellera les quatre intégrales précédentes, dans l'ordre,  $(I)$,
$(II)$, $(III)$ et $(IV)$.

\bigskip
Pour $(I)$, regardons d'abord $\int^{X_l}_{\lambda\eps^{-1/3}}
\frac{-2}{Y_0^+Y_0^-} dX$. D'après l'équation différentielle
\eqref{eq:Y0} vérifiée par $Y_0$, on a 
\begin{align*}
(Y_0^+-Y_0^-)'(X)               &= \frac{-2}{Y_0^+(X)Y_0^-(X)}
(Y_0^+-Y_0^-)(X)\\ 
\Longrightarrow (Y_0^+-Y_0^-)(X)&= (Y_0^+-Y_0^-)(0) e^{\int^X_0
  \frac{-2}{Y_0^+Y_0^-}dX}\qquad\text{(en intégrant l'équation).} \\ 
\intertext{Or, quand $X\longrightarrow\infty$, on a vu avec
  \eqref{eq:dY0} que}
(Y_0^+-Y_0^-)(X) &\sim  \frac{4}{e}i X^2 e^{-\frac{2}{3}X^3}\,,
\end{align*}
{donc}
\begin{align*}
\int^{X_l}_{\lambda\eps^{-1/3}} \frac{-2}{Y_0^+Y_0^-} dX 
 &= \ln\biggl(\frac{(Y_0^+-Y_0^-)(X_l)}{i}\biggr) 
-\ln\biggl(\frac{(Y_0^+-Y_0^-)(\lambda\eps^{-1/3})}{i}\biggr)\\
 &=\ln\bigl({-i(Y_0^+-Y_0^-)(X_l)}\bigr) +1 -2\ln 2
+\frac{2\lambda^3}{3\eps}  -2\ln(\lambda\eps^{-1/3})+o(1)\\
 &=\ln\bigl({-i(Y_0^+-Y_0^-)(X_l)}\bigr) 
+\frac{2\lambda^3}{3}\frac{1}{\eps} 
+\frac{2}{3} \ln\eps  -2\ln 2\lambda +1+o(1)\,.
\end{align*}

Pour tous les termes à venir (y compris pour les autres intégrales),
l'intégration au voisinage de $X_l$ et de tout point fini ne donne au
final que du $o(\eps)$, et il n'y a donc que ce qui se passe quand
$X\rightarrow\infty$ qui nous intéresse. Autrement dit, la différence
entre $Y_n^\pm$ et son développement en $1/X$ est négligée. Ainsi, pour
$\int^{X_l}_{\lambda\eps^{-1/3}}\frac{\eps^{1/3}X}{Y_0^+Y_0^-} dX$, on
utilise le développement asymptotique de $Y_0$:
\begin{align*} Y_0=-\frac{1}{X}-\frac{1}{2X^4}+\gO\Bigl(\frac{1}{X^7}\Bigr)&
\Longrightarrow \frac{1}{Y_0^2}=X^2-\frac{1}{X}
+\gO\Bigl(\frac{1}{X^4}\Bigr)\\ 
\int^{X_l}_{\lambda\eps^{-1/3}}\!\frac{\eps^{1/3}X}{Y_0^+Y_0^-} dX&=
 \eps^{1/3} \int_{\lambda\eps^{-1/3}}\!\! (X^3-1) dX +\gO(\eps^{1/3}) \\
 &=-\eps^{1/3}\frac{\lambda^4}{4\eps^{4/3}}
 +\eps^{1/3}\lambda\eps^{-1/3} +\gO(\eps^{1/3})\\
 &=-\frac{1}{\eps}\frac{\lambda^4}{4} +\lambda +\gO(\eps^{1/3})\,.
\end{align*}

Et pour la dernière partie de $(I)$,
\begin{align*}
\int^{X_l}_{\lambda\eps^{-1/3}}\!\frac{a}{Y_0^+Y_0^-} dX &=
-\frac{\eps}{16}\int^{X_l}_{\lambda\eps^{-1/3}}\!\frac{-2}{Y_0^+Y_0^-} dX
\Bigl(1+\gO(\eps)\Bigr) \\
 &=-\frac{\lambda^3}{24}+o(1)\,.
\end{align*}

D'où finalement
\begin{multline*}
(I)=\frac{1}{\eps}\Bigl(\frac{2\lambda^3}{3}-
\frac{\lambda^4}{4}\Bigr) +\frac{2}{3}\ln\eps -2\ln\lambda -2\ln 2 +1
+\lambda\\
 -\frac{\lambda^3}{24} +\ln\biggl(\frac{(Y_0^+-Y_0^-)(X_l)}{i}\biggr)
+o(1)\,.
\end{multline*}
 
\bigskip
On peut procéder de manière tout a fait identique pour $(II)$
$$(II)= \eps^{1/3}\int^{X_l}_{\lambda\eps^{-1/3}}
   \frac{-2+\eps^{1/3}X-a}{Y_0^+Y_0^-}\:Q_1(X)\ dX$$
en utilisant la série asymtotique de $Y_1$ dans $Q_1$ et ne négligeant
ce qui est en facteur de $a=\gO(\eps)$.

\begin{gather*}
 Y_1(X)=-\frac{1}{4X^3}+\gO(X^{-6})\quad\Longrightarrow\quad
-\frac{2Y_1}{Y_0^3}=-\frac{1}{2}+o(1)\\
(II)=\eps^{1/3}\int^{\lambda\eps^{-1/3}} \!\frac{-2+\eps^{1/3}X}{2}dX
+o(1)= \frac{\lambda^2}{4}-\lambda+o(1)
\end{gather*}

\bigskip
Et de m\^eme pour $(III)$.
\begin{gather*}
(III)= {\eps}^{2/3}\int^{X_l}_{\lambda\eps^{-1/3}} 
   \!\frac{(-2+\eps^{1/3}X)}{Y_0^2} \Bigr( \frac{Y_1^2}{Y_0^2} -
   2\frac{Y_2}{Y_0} \Bigl) dX \\
Y_2(X)=-\frac{1}{8X^2}+\gO\Bigl(\frac{1}{X^5}\Bigr)\Longrightarrow
-\frac{2Y_2}{Y_0^3}=-\frac{X}{4} +\gO\Bigl(\frac{1}{X^2}\Bigr)
\text{; et } 
\frac{Y_1^2}{Y_0^4}=\gO\Bigl(\frac{1}{X^2}\Bigr)\\
(III)= {\eps}^{2/3}\int^{\lambda\eps^{-1/3}} \!\! (-2+\eps^{1/3}X)
\left(-\frac{X}{4}\right) +o(1) = \frac{\lambda^3}{12}
-\frac{\lambda^2}{4} +o(1)
\end{gather*}

\bigskip


Pour la dernière intégrale,
\begin{align*}
(IV) &\leq \eps \int_{\lambda\eps^{-1/3}}
\frac{-2+\eps^{1/3}X}{Y_0^2}\left(\f{1}{X}+\eps^{1/3}\right) 
\gO_\eps(1)\ dX\\
     &\leq \eps  \int_{\lambda\eps^{-1/3}}\left( \gO_\eps(X)+\eps^{1/3}     
      \gO_\eps(X^2)\right) dX\\
     &\leq \eps \:  \gO(\eps^{-2/3})\\
     &= o(1)
\end{align*}

\bigskip
En fin de compte, quand on additionne $(A)+(I)+(II)+(III)+(IV)$, on obtient
\begin{align}
\frac{1}{\eps}F(x_l) &= \frac{1}{\eps}\frac{4}{3} +\frac{2}{3}\ln\eps
+\frac{4}{3} -4\ln 2 +\ln\Bigl(\frac{(Y_0^+-Y_0^-)(X_l)}{i}\Bigr)
+o(1) \notag\\
\intertext{Donc}
e^{-F(x_l)/\eps} &\sim e^{-{4}/{3\eps}} \eps^{-2/3}
  \frac{16}{e^{4/3}} \frac{i}{\bigl(Y_0^+-Y_0^-\bigr)(X_l)}\label{eq:b4}
\end{align}
\smallskip

{\small{\bf Remarque}
\par\penalty 10000
On peut signaler que le calcul de $(A)$ ci-dessus n'est pas
indispensable et qu'il suffit de poser $\lambda =0$ dans $(B)$ pour
obtenir le bon résultat. Le calcul complet de $(A)$ a cependant été
conservé parce qu'{\sl a priori}, il est indispensable, qu'il peut
l'être dans d'autres exemples et parce qu'il permet de détecter
d'éventuelles erreurs. Dans l'autre sens, remplacer $\lambda$ par $X_l
\eps^{1/3}$ dans $(A)$ ne suffit pas pour obtenir le résultat complet,
mais donne déjà les termes $e^{-{4}/{3\eps}} \eps^{-2/3}$.

On peut en revanche choisir effectivement n'importe quel $\lambda$
fixé dans $]0,2[$, ou même un $\lambda$ variant très peu avec $\eps$
comme par exemple $\eps^{1/20}$, on aboutira au bon résultat, avec 
quelques complications dans le second cas.}

\vskip 20pt
En combinant les résultats intermédiaires \eqref{eq:b2},
\eqref{eq:b3} et \eqref{eq:b4}, cela nous donne finalement:
\begin{align*}
\eqref{eq:b}\qquad b \sim& i\eps e^{-{4}/{3\eps}}
\eps^{-2/3} \frac{16}{e^{4/3}}
\frac{\eps^{-1/3}\bigl(Y_0^+-Y_0^-\bigr)(X_l)}{\bigl(Y_0^+-Y_0^-\bigr)(X_l)}
\frac{1}{\sqrt{2\pi\eps}}\\
b \sim&  i\:\frac{e^{-\frac{\scriptstyle 4}{\scriptstyle
                3\eps}}}{\sqrt\eps} \frac{8\sqrt 2}{\sqrt\pi e^{4/3}}
\end{align*}
Soit le
\begin{Th}
$$ \label{eq:B}  \alpha^+-\alpha^- \sim  i\:\frac{e^{-\frac{\scriptstyle 4}{\scriptstyle
                  3\eps}}}{\sqrt\eps} 
               \frac{8\sqrt 2}{\sqrt\pi e^{4/3}}$$
\end{Th} 
On retrouve le résultat de \cite{FS} de la proximité exponentielle des 
deux $\alpha$. Et on confirme que la borne trouvée: 
$\exp\left(R(-1+\rho)/\eps\right)$ était bien optimale. 

\section{Conséquence sur les coefficients de la série $\mathbf{\hat a(\eps)}$}
\label{sec:an}
D'après \cite{FS}, la fonction $a(\eps)$ est holomorphe dans le
domaine 
$${\mathcal M}=\biggl\{ \eps \Bigl/ \arg(\eps) \in \Bigl]
\frac{-5\pi}{2} +\delta,  \frac{\pi}{2}-\delta\Bigr[, 0<|\eps|\leq|\eps_1|
\biggr\}.$$ 
Pour $\eps_0$ assez petit ($|\eps_0|<|\eps_1|$), le chemin $\gamma$,
défini par un arc de cercle de rayon $| \eps_0|$ et deux segments
proches de $[0, \eps_0]$ (voir figure~\ref{fig:integ}) est inclu dans
${\mathcal M}$, et est rétractile en un point dans ce domaine..

\bigskip
\begin{figure}[ht!]
\center{\includegraphics[scale=0.8]{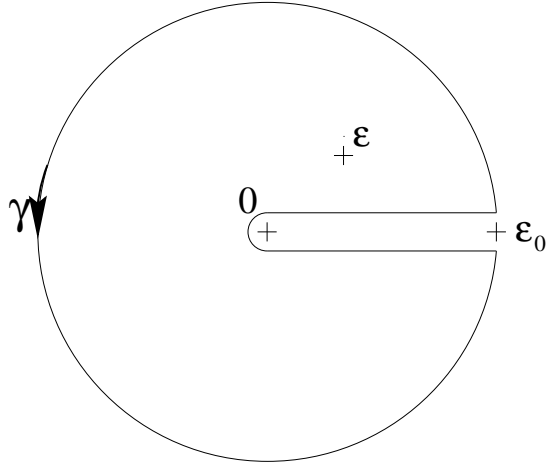}}
\caption{ Chemin d'intégration $\gamma$.}\label{fig:integ}
\end{figure}

Donc pour $\eps \in {\mathcal M}$, on a 
\begin{align*}
a(\eps) &=\frac{1}{2i\pi}\int_\gamma \frac{a(u)}{u-\eps}du\\
        &=\frac{1}{2i\pi}\int_{|u|=\rho} \frac{a(u)}{u-\eps}du +
\frac{1}{2i\pi}\int^{\eps_0}_0
\biggl(\frac{a(ue^{-2i\pi})-a(u)}{u-\eps}\biggr) du\\ 
\intertext{puis en développant en série de $\eps$ et en identifiant les
  coefficients obtenus} 
a_n &=\frac{1}{2i\pi}\int_{|u|=\rho} a(u)u^{-n-1} du
+\frac{1}{2i\pi}\int^{\eps_0}_0
\biggl(a(ue^{-2i\pi})-a(u)\biggr)u^{-n-1} du\\
\intertext{comme le premier terme n'est que de l'ordre de $\rho^{-n}$}
a_n &\sim \frac{1}{2i\pi}\int^{\eps_0}_0\Bigl(a^-(u)-a^+(u)\Bigr)u^{-n-1}
    du\\
\intertext{et que dans cette intégrale, la contribution principale
  vient du voisinage de $0$, le théorème \ref{eq:B} mène à}
a_n &\sim -\frac{1}{2i\pi} i\frac{8\sqrt 2}{e^{4/3}\sqrt\pi}
\int^{\eps_0}_0 \frac{e^{-4/3u}}{\sqrt u}u^{-n-1} du\,.\\
\intertext{On effectue alors le changement de variable
  $t=\frac{4}{3u}, du=-\frac{4}{3}\frac{dt}{t^2}$} 
a_n &\sim \frac{1}{2\pi}\frac{8\sqrt 2}{e^{4/3}\sqrt\pi}
\int_{+\infty}^{4/3\eps_0} e^{-t}t^{n-1} \sqrt t \; dt \times
\sqrt{\frac{3}{4}} \biggl(\frac{3}{4}\biggr)^n\,.\\
\intertext{Pour $n$ assez grand, $\int_{4/3\eps_0}^0
 \! e^{-t}t^{n-1/2} dt$ est négligeable face à
  $\int_{+\infty}^{4/3\eps_0} \!e^{-t}t^{n-1/2} dt$.}
\intertext{Donc} 
a_n &\sim \frac{-1}{2\pi}\frac{8\sqrt 2}{e^{4/3}\sqrt\pi}
    \sqrt{\frac{3}{4}} \biggl(\frac{3}{4}\biggr)^n  
    \Gamma\biggl(n+\frac{1}{2}\biggr) \\
\intertext{Or}
\Gamma\biggl(n+\frac{1}{2}\biggr) &\sim \sqrt{2\pi n}
\biggl(\frac{n-1/2}{e}\biggr)^{n-1/2} \sim \sqrt{2\pi}
\biggl(\frac{n}{e}\biggr)^n\\ 
\intertext{Donc}
a_n &\sim \frac{-1}{2\pi}\frac{8\sqrt 2}{e^{4/3}\sqrt\pi}
    \sqrt{\frac{3}{4}} \biggl(\frac{3}{4}\biggr)^n \sqrt{2\pi} 
\biggl(\frac{n}{e}\biggr)^n
\end{align*}
\begin{Th}
On a un équivalent exact pour les coefficients des $\alpha$ \og à
canards\fg:
$$a_n  \sim \frac{-4\sqrt 3}{\pi e^{4/3}}  \biggl(\frac{3n}{4e}\biggr)^n$$
\end{Th}
On rappelle que tous les $\alpha$ donnant des solutions \og canards\fg
sont exponentiellement proches en $\eps$. Ils ont donc tous le
m\^eme développement asymptotique. Par ailleurs, le résultat obtenu (à
$\eps$ donné) en sommant la série \og au plus petit terme\fg donne
aussi une approximation de $\alpha$ exponentiellement bonne, et est
donc une valeur à canard~\cite{FS}. 

\section{Résultats numériques}

On peut facilement implémenter le calcul de la suite des $a_n$ en
utilisant les formules de récurrence \eqref{eq:serie}. En particulier,
en calculant les 155 premiers termes (calcul effectué par Franck
Michel et amicalement transmis) et en regardant le résultat de la
multiplication $b_n=a_n\times(4e/3n)^n$, on obtient:

\bigskip
\begin{tabular}{|c|c||c|c||c|c|}
\hline
$b_{135}\ $&$-0,5417512651\ $&$b_{136}\ $&$-0,5418690317\ $&$b_{137}\
$&$-0,5419854885\ $\\
\hline
$b_{138}\ $&$-0,5421006603\ $&$b_{139}\ $&$-0,5422145711\ $&$b_{140}\
$&$-0,5423272443\ $\\
\hline
$b_{141}\ $&$-0,5424387024\ $&$b_{142}\ $&$-0,5425489682\ $&$b_{143}\
$&$-0,5426580621\ $\\
\hline
$b_{144}\ $&$-0,5427660064\ $&$b_{145}\ $&$-0,5428728208\ $&$b_{146}\
$&$-0,5429785257\ $\\
\hline
$b_{147}\ $&$-0,5430831405\ $&$b_{148}\ $&$-0,5431866841\ $&$b_{149}\
$&$-0,5432891757\ $\\
\hline
$b_{150}\ $&$-0,5433906324\ $&$b_{151}\ $&$-0,5434910728\ $&$b_{152}\
$&$-0,5435905137\ $\\
\hline
$b_{153}\ $&$-0,5436889722\ $&$b_{154}\ $&$-0,5437864645\ $&$b_{155}\
$&$-0,5438830066\ $\\
\hline 
\end{tabular}
\null
\vskip 10pt
Si on cherche alors les $b_n$ sous la forme $b_n=C+a/\sqrt{n}$ et
$b_n=C+a/\sqrt[3]{n}$, on obtient (en utilisant une formule de
linéarisation aux moindres carrés):
$$b_n=-0,5736898877+\frac{0,3710889332}{\sqrt n}   ,\qquad
b_n=-0,5891153498+\frac{0,2429519729}{\sqrt[3]{n}}.$$
Et la constante calculée ci-dessus vaut
$$-0,5813148764.$$
On constate donc que la constante calculée est compatible avec les
résultats numériques, mais que les «termes correctifs» du
développement des $a_n$ restent non négligeables.

\vskip 30pt

Il est aussi possible d'estimer directement la valeur de
$b=\alpha^+-\alpha^-=2\Im(\alpha^+)$ en cherchant les valeurs numériques
de $\alpha^+$ pour différentes valeurs de $\eps$.

On intègre numériquement, pour certains $\eps$ raisonnablement petits, 
l'équation différentielle~\eqref{eq:vdpdebase}, par exemple
d'une part le long du chemin comprenant les segments $[-1+10i,0]$
puis $[0,1]$ et d'autre part le long de $[9,1]$. Ces deux chemins
ont la propriété de descendre le relief (cf.~fig.~\ref{fig:relvdp}, et 
\cite{FS}), c'est à dire qu'une erreur dans l'estimation des valeurs 
initiales (en $-1+10i$ et en $9$) devient exponentionnellement petite
une fois arrivés en $1$; on prendra donc les valeurs $v(-1+10i)=i/10$
et $v(9)=-1/10$. De plus, un de ces chemins part de la montagne Nord,
l'autre de la montagne Est, et la définition de $\alpha^+$ est qu'il
s'agit de la valeur de $\alpha$ pour laquelle il existe une solution
tendant vers $0$ en l'$\infty$ sur ces deux montagnes (N.B.: on
retrouvera la manière de construire les $\alpha$ dans la partie
suivante), ce qui sera vrai si et seulement si les deux $v(1)$ trouvés
sont égaux.

\begin{figure}[ht!]
\center{\includegraphics[scale=0.64]{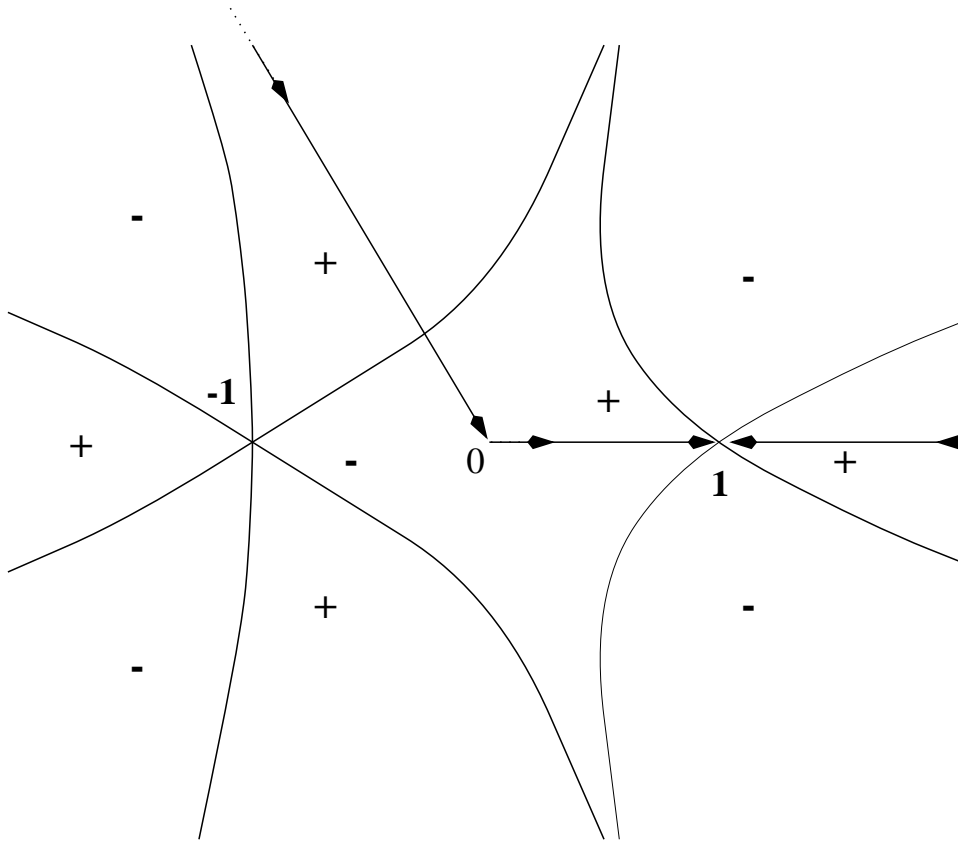}}
\caption{Relief pour l'équation de Van der Pol, et chemins pour 
l'intégration numérique}\label{fig:relvdp}
\end{figure}

Une recherche numérique de ce type donne les approximations suivantes:
\vglue 10pt
\begin{tabular}{|c|l|c|}
\hline
$\eps$ & $\qquad\alpha^+$ & $2\Im (\alpha^+)\times e^{4/3\eps}\sqrt\eps$\\
\hline\hline
$\ 0,20\ $  & $\ 0,9684+0,00153i$          & 1,07\\
\hline
$0,17$ & $\ 0,9733+0.00055i$           & 1,16\\
\hline
$0,14$ & $\ 0,9800+0.000120i$          & 1,23\\
\hline
$0,08$ & $\ 0,9893+1,40\cdot 10^{-7}i$ & 1,37\\
\hline
$0,06$ & $\ 0,9921+6,48\cdot 10^{-10}i$& 1,42\\
\hline
$0,05$ & $\ 0,9935+8,5\cdot 10^{-12}i $& 1,44\\
\hline
$0,04$ & $\ 0,9948+1,23\cdot 10^{-14}i\qquad$& 1,47\\
\hline
\end{tabular}
\vskip 10pt
Ces résultats sont compatibles avec la valeur de la constante trouvée
(égale environ à 1,68), avec des termes correctifs en $\sqrt\eps$ et
$\sqrt[3]\eps$.

\chapter{Généralisation}\section{Préliminaires}

\subsection{Rappel : estimation de Cauchy}

On rappelle pour les démonstrations à venir la propriété de
Cauchy, qu'on utilisera sous la forme suivante:
\begin{Lemmetech}\label{th:derhol}
Soit $f(x_1,\ldots,x_n)$ une fonction holomorphe en
$x=(x_1,\ldots,x_n)$, bornée pour $x$ dans l'ouvert $U=U_1\times
U_2\times\ldots\times U_n)$, et soit $F_1$ un fermé inclus dans
l'ouvert $U_1$. Alors
$$\frac{\partial f}{\partial x_1}(x)\ 
\text{est bornée pour } x\in F_1\times U_2\times\ldots\times U_n $$
\end{Lemmetech}

\smallskip
Un corollaire immédiat de ce lemme indique que si ce n'est pas
$f(x)$, mais $\frac{f(x)}{g(x_k)}$ qui est bornée, alors la fonction  
$\frac{1}{g(x_k)}\frac{\partial f}{\partial x_1}(x)$ est bornée elle
aussi.

\subsection{Forme normale pour une équation différentielle}

Plut\^ot que travailler sur une équation écrite sous la
forme~\eqref{ed} vue dans l'introduction:
$$\eps y'=y f(x,\eps)+h(x,\eps)+\eps y^2 P(x,\eps,y)\,$$ 
nous allons manipuler dans cette section une équation sous une
forme~\eqref{eq:gen} apparemment moins générale:
\begin{equation}
\eps y'=\bigl(x^p f(x)+\eps
g(x,\eps)\bigr) y +h(x,\eps) +\eps y^2 P(x,\eps,\eps y)\,.
\tag{\ref{eq:gen}}
\end{equation}
En fait, on se ramène facilement à une équation de ce type.

\smallskip
Soit l'équation
$$\eps u'= \Psi(x,u(x),\eps)\,,$$ 
où on demande juste que $\Psi$ soit holomorphe en $x$, méromorphe en
$u$ et qu'elle admette un développement asymptotique en $\eps$ quand
$\eps$ tend vers $0$. On suppose aussi qu'il existe une courbe lente
$u_0(x)$, telle que $\Psi(x,u_0(x),0)=0$.

On pose $u(x)=u_0(x)+\eps y(x)$, ou, ce qui donne souvent un résultat
plus agréable à manipuler, $u(x)=u_0(x)\bigl(1+\eps y(x)\bigr)$, et on
linéarise l'équation en $y$.
\begin{align*}
\eps u'&= \Psi(x,u,\eps)\,,\\
\intertext{donne, en regroupant dans $\psi_0$ et $\psi_1$ tout ce qui
est linéaire en $u$, et dans le dernier terme tout ce qui contient
un pôle en $u=0$,}
       &=\psi_0(x,\eps) + u\psi_1(x,\eps) + u^2 \psi_2(x,\eps,u) 
+\f{1}{u^k}\psi_3(x,\eps,u)
\end{align*}
On écrit alors $u=u_0 (1+\eps y)$:
\begin{multline*}
\eps u'_0\bigl(1+\eps y(x)\bigr) +\eps^2 u_0 y'=\psi_0 +u_0\psi_1+
\eps u_0 y\\
+ u_0^2\psi_2+u_0^2(2\eps y +\eps^2 y^2)\psi_2(x,\eps,u_0(1+\eps y)) 
+\f{1}{u_0^k(1+\eps y)^k}\psi_3(x,\eps,u_0(1+\eps y))
\end{multline*}
{Comme $u_0$ est courbe lente, les termes de degré $0$ en $\eps$
s'éliminent; on divise alors par $\eps$, puis on linéarise autour de
$y=0$: les termes en $y^k$, $k>0$ ont donc comme facteur un terme
$\f{\eps^k}{\eps}$, ce qui donne bien en particulier, pour $k\geq 2$, 
$\eps y^2 P(x,\eps, \eps y)$. On obtient donc bien une équation du
type annoncé ci-dessus, sachant que les fonctions $f$, $g$, $h$ et $P$
contiennent éventuellement des termes avec des facteurs en
$\f{1}{u_0^{k+1}}$.}

Ainsi, dans l'exemple de l'équation de Van der Pol~\eqref{eq:vdpdebase},
cette procédure de normalisation nous donne:
\begin{align*}
\eps v\f{dv}{du}&=(1-u^2)v+\alpha-u\\
\eps \f{dv}{du}&=(1-u^2)+\f{1+a(\eps)-u}{v}\,.\\
\intertext{Puis en posant $v=v_0(1+\eps y)$, $v_0$ étant
égal à $\f{-1}{u+1}$,}
\eps\left(v_0(1+\eps y)\right)'&=1-u^2+\f{1+a-u}{v_0(1+\eps y)}\,,\\
\eps v_0 (\eps y')&=1-u^2+\f{1+a-u}{v_0}\left(1-\eps y+\f{\eps^2
y^2}{1+\eps y}\right) -\eps v'_0 (1+\eps y)\,,\\
\left(-\f{\eps}{1+u}\right)\eps y' &=
(1-u^2)\left(\eps y-\f{\eps^2 y^2}{1+\eps y}\right)-a(1+u)
\left(1-\eps y+\f{\eps^2 y^2}{1+\eps y}\right)\\
&\phantom{= +}+\f{\eps}{(1+u)^2}(1+\eps y)\,,\\
\eps y' &=y\left(-(1-u)(1+u)^2+
\eps\left(\f{a}{\eps}(1+u)^2-\f{1}{1+u}\right)\right)\\
&\phantom{= +}+\f{a}{\eps}(1+u)^2-\f{1}{1+u}+\eps y^2\f{1}{1+\eps y}
\left((1-u)(1+u)^2+a(1+u)^2\right)\,.
\end{align*}
Pour cette équation, l'étude a été faite au voisinage des deux points
cols du relief construit à partir de $-(1-u)(1+u)^2$: autour de $u=1$,
dont on a montré qu'il n'était pas un point tournant dans certains cas
(\cite{FS}); et au voisinage de $-1$, au chapitre~\ref{chap:vdp}. Ce
sont ces deux études que l'on souhaite généraliser ici.

\section{Existence de solutions d'une équation différentielle près
d'un point tournant}  

Dans les paragraphes \ref{sec:loin} et \ref{sec:pres}, nous allons 
généraliser le théorème \ref{ThVdPPtfixe} qui a été
démontré dans le cas de l'équation de Van der Pol, au cas suivant:
\begin{equation}\tag{{\sc i}}\label{eq:gen}
\eps y'=\bigl(x^p f(x)+\eps
g(x,\eps)\bigr) y +h(x,\eps) +\eps y^2 P(x,\eps,\eps y)\,.
\end{equation}
Il est clair que $0$ est un point col pour le relief $\Re\left(\int^x
t^p f(t) dt\right)$ correspondant, et c'est donc en général un point
tournant pour l'équation différentielle. On veut étudier le
comportement des solutions de cette équation au voisinage de ce point.
Le but va être de démontrer que si on a les bonnes hypothèses sur les
fonctions intervenant dans l'équation, on a l'existence d'une vraie
solution bornée, d'abord pour des $x$ tels que $|x|>\rho$, où $\rho>0$
est une constante quelconque fixée indépendante de $\eps$ (résultat
classique, redémontré au paragraphe~\ref{sec:loin} ci-dessous, pour
$\eps$ assez petit). Ensuite, que cette solution peut bien souvent
être prolongée dans des domaines de la m\^eme forme, mais avec un
$\rho$ tendant cette fois vers $0$ avec $\eps$:
$|x|>\rho(\eps)=X_l\eps^r$ ($X_l$ constante assez grande, mais
indépendante de $\eps$, et $r$ un rationnel strictement positif
donné). Ce deuxième résultat sera démontré dans un second temps, dans
le paragraphe~\ref{sec:pres}, et on verra que si la solution prolongée
n'est pas nécessairement bornée (quand $\eps\rightarrow 0$) dans tout
son domaine d'existence, sa croissance reste contrôlée et n'est pas,
en tous cas, exponentielle.

Chacun de ces deux résultats ne sera en général vrai que dans certains
domaines déterminés non seulement par la norme de $x$, mais aussi par le
relief; de fait, ces domaines correspondront à peu près à des secteurs
centrés en $0$, avec des $x$ vérifiant la condition donnée ci-dessus:
$|x|>\rho(\eps)$.

\subsection{Existence d'une solution {\protect «}loin{\protect »} de $0$}
\label{sec:loin}

On se fixe donc une constante $\rho$ positive, arbitrairement petite.

On suppose qu'on a, pour l'équation \eqref{eq:gen}, les hypothèses
(\HB) suivantes:
\begin{list}{$\bullet$}{}
\item  $f$ est holomorphe dans $\Dom_0$, où  $\Dom_0$ est un
domaine ouvert, fini ou infini, contenant $0$, et $f$ ne s'annule pas
sur $\Dom_0$.
\item  $g$, $h$ et $P$ sont holomorphes en $x$ sur $\Dom_0$,
sauf éventuellement en $0$ où $g$ ou $h$ peuvent n'être que
méromorphes et avoir un p\^ole. 
\item  $g$, $h$ et $P$ sont holomorphes en $\eps$ dans des
secteurs ouverts $S_0$, centrés en $0$, et bornées pour tout $x\in
\Dom_0$ fixé quand $\eps\rightarrow0$: on suppose qu'il existe $q$ tel 
que les fonctions $P$, $x^q g(x,\eps)$ et $x^q h(x,\eps)$ ont un
développement asymptotique en $\eps$ quand $\eps$ tend vers $0$.
\item  $P$ est holomorphe en $\eps y$ dans un voisinage de $0$.
\item  On pose $F(t)=\int^t (u^p f(u))du$ et
$R(t)=\Re\left(\frac{F(t)|\eps|}{\eps}\right)$. On considérera à
$\theta=\arg\eps$ fixé un domaine {\bf fermé} $\Dom_\theta \subset
\Dom_0$, tel que $\Dom_\theta$ soit accessible à partir d'un sommet
$s_\theta$, qui est soit infini, soit dans $\Dom_0\setminus\Dom_\eps$
avec le relief $R(t)$, et tel que $\forall x \in \Dom, |x|\geq \rho$;
on reprendra, pour ce qui concerne l'accessibilité, les notations de
la définition~\ref{def:acc}. Les domaines $\Dom_\eps$ peuvent être non
bornés, si $\Dom_0$ l'est.
\end{list}

L'équation \eqref{eq:gen} a une solution formelle $\hat
y(x,\eps)=\sum^\infty_{n=0} y_n(x) \eps^n$ pour $|x|>\rho$, dont tous
les coefficients $y_n(x)$ sont holomorphes en $x$. Cela se vérifie en
remplaçant dans l'équation $y$ par cette série et en identifiant les
coefficients de $\eps^n$: on trouve alors des relations de récurrence
donnant $y_n$ en fonction des $y_p$ et de leurs dérivées (pour
$p<n$). La première relation, obtenue en posant $\eps=0$ donne
$$y_0(x)=-\frac{h(x,0)}{x^p f(x)}$$

\bigskip

Nous allons démontrer le théorème suivant:
\begin{Th}\label{th:loin}
On suppose que pour l'équation \eqref{eq:gen}, les hypothèses {\rm
(\HB)} sont vérifiées et que les fonctions
$$\f1{x^p f(x)}\ ,\ \frac{1}{x^p f(x)}\sup_{|y|\leq
\delta}\left(P(x,\eps,\eps y)\right)\ ,\ \frac{g(x,\eps)}{x^p f(x)}\
\text{ et }\ \ \frac{ h(x,\eps)}{x^p f(x)}$$  
sont bornées\footnote[2]{autrement dit, la partie linéaire homogène de
l'équation différentielle ne devient jamais négligeable par rapport
aux autres termes ou aux constantes} pour un certain $\delta$,
$x\in\Dom_0, |x|>r$ ($r$ constante arbitrairement petite) et $\eps\in
S_0,|\eps|<\eps_0$.  On fait varier $\eps$ dans un secteur $S_0$
d'ouverture arbitrairement petite.

Alors, pour tout réel positif $\rho> r$, il existe un domaine fermé
éventuellement non borné $\Dom\subset\Dom_0$, de la forme 
$$\Dom=\bigcap_{\eps\in S_0}\Bigl\{x\in\CC\ \bigl/\ |x|\geq\rho,\,
x\text{ accessible d'un sommet $s_\eps\!$ avec le relief }
\Re\frac{F(x)|\eps|}{\eps}\Bigr\}$$ et une solution $y(x,\eps)$ de
\eqref{eq:gen} existant pour $x\in\Dom$ qui soit bornée indépendamment
de $|\eps|$ sur ce domaine. Cette fonction $y(x,\eps)$ tend vers la
fonction $y_0(x)$ quand $\eps$ tend vers $0$ dans $S_0$, uniformément
pour $x\in\Dom$.
\end{Th}
Quelques remarques sur les hypothèses de ce théorème:
\begin{enumerate}
\item L'hypothèse $1/\Bigl(x^p f(x)\Bigr)$ bornée implique bien
entendu que $f$ ne s'annule pas sur $\Dom$.
\item Comme il existe un $\lambda_\rho$ tel que
$\lambda_\rho|u^p f(u)|\geq |g(u,\eps)|, \forall u\in\Dom$ et $\eps$
assez petit, cela signifie que le (vrai) relief correspondant à
 $\Re\left(\int^t \Bigl(u^p f(u)+\eps g(u,\eps)\Bigr) du\right)$ ne serait qu'une
petite perturbation du relief effectivement considéré.
\item Comme $h(x,\eps)/\Bigl(x^p f(x)\Bigr)$ est bornée pour tout
$\eps$ assez petit et élément d'un secteur $S_0$, la fonction $y_0(x)$
est elle aussi bornée sur $\Dom$. Ensuite, d'après l'inégalité de Cauchy,
$y'_0(x)$ fait aussi partie des fonctions bornées sur $\Dom$.
\item Les points $s_\eps$ et les reliefs ne dépendent que de
l'argument de $\eps$, et pas de sa norme. L'intersection se fait donc
uniquement sur le (petit) intervalle des arguments de $\eps$ possibles
pour $S_0$. Comme le relief et les sommets dépendent contin\^ument de
$\arg\eps$, en prenant $S_0$ suffisamment peu ouvert, on est assuré
que le domaine $\Dom$ sera non vide.
\end{enumerate}

\bigskip 
Le théorème revient à montrer qu'il existe une vraie solution
holomorphe correspondant à la solution formelle donnée ci-dessus dans
certains domaines $\Dom$ pour tout $\eps$ assez petit. On commencera
par considérer le cas où $\arg\eps$ est une constante, avec le relief
et le sommet $s$ correspondants. On omettra le plus souvent dans la
suite de rappeler la dépendance de la fonction $y$ en $\eps$.

\medskip

On cherche cette solution holomorphe sous la forme $y=y_0+\eps z$. 
L'équation \eqref{eq:gen} devient alors
\begin{multline*}
\eps z'(x)=\bigl(x^p f(x)+\eps g(x,\eps)\bigr) z(x) +y_0 (x) g(x,\eps)
-y'_0(x) +\frac{h(x,\eps)-h(x,0)}{\eps}\\
 +(y_0(x)+\eps z(x))^2 P(x,\eps,\eps y_0(x)+\eps^2 z(x))
\end{multline*}
En utilisant la formule de variation de la constante, on obtient 
$$z (x) =\frac{1}{\eps} \int_{\gamma_x}
e^\frac{F(x)-F(t)}{\eps}
\Bigl(y_0 (t) g(t,\eps)-y'_0(t) +\frac{h(t,\eps)-h(t,0)}{\eps}
+H(t,\eps,z(t)) \Bigr) dt$$ 
avec ${\gamma_x}$ un chemin partant du sommet $s$ et descendant le
relief jusqu'à $x$, et avec la fonction $H(t,\eps,z)=\eps z g(t,\eps)
+(y_0(t)+\eps z)^2 P(t,\eps,\eps y_0(t)+\eps^2 z)$.

\medskip

On souhaite appliquer un théorème de point fixe, pour $\eps$ assez
petit fixé. On considère pour cela les deux opérateurs suivants:

\begin{align*}
{\HHH}:\tilde{\ZZ}\subset{\ZZ}&\rightarrow{\HH}\\
           z& \mapsto {\HHH}(z),
\intertext{la fonction ${\HHH}(z)$ étant telle que 
${\HHH}(z)(x,\eps)=H(x,\eps,z(x,\eps))$ pour tout 
$(x,\eps)\in\Dom\times S_0$,}
\end{align*} 
et\label{Operateurs}
\begin{align*}
{\ZZZ}:\tilde{\HH}\subset{\HH}&\rightarrow{\ZZ}\\ 
           H& \mapsto z=\frac{1}{\eps}
\int_{\gamma_x} e^\frac{F(x)-F(t)}{\eps} \Bigl(y_0 (x)
g(t,\eps)-y'_0(t) +\frac{h(t,\eps)-h(t,0)}{\eps}
+H(t) \Bigr) dt
\end{align*}
où $\HH$ et $\ZZ$ sont les espaces fonctionnels suivants:
$$\ZZ=\{z(x,\eps)/\ z \hbox{ holomorphe et bornée dans } \Dom\times S_0\}$$ 
\begin{multline*}
\HH=\Bigl\{ H(x,\eps)\,/\ H\hbox{ est holomorphe dans } \Dom\times
S_0,\\
\hbox{ et }\sup_{x\in \Dom} \frac{H(x,\eps)}{x^p f(x)}
\hbox{ est borné indépendamment de }\ \eps \Bigr\}
\end{multline*}
Ces deux espaces sont munis des normes correspondantes $\|\cdot\|$ et
$\|\cdot\|_{\HH}$. Avec les hypothèses du théorème, les fonctions
$g(x,\eps)$, $P(x,\eps,\eps\delta)$ et $h(x,\eps)$ sont dans $\HH$,
ainsi que $\Dp{P}{(\eps y)}$, et avec une norme dans cet espace
que l'on peut majorer indépendamment de $\eps$.
\medskip

Alors $\HHH$ et $\ZZZ$ sont des opérateurs contractants. 
En effet, on peut écrire (en omettant de rappeler les dépendances de
$g$, $y_0$, $z_k$ et$H_k$  en $x$ et $\eps$), 
\begin{multline*} 
{\HHH}(z_1)-{\HHH}(z_2)= \biggl[ \eps g+ 
\bigl(2\eps y_0 +\eps^2 (z_1+z_2)\bigr)P(x,\eps,\eps{y_0}+\eps^2 z_2)\\
+(y_0^2+2\eps y_0 z_1+\eps^2 z^2_1) \frac{P(x,\eps,\eps{y_0}+\eps^2 z_1)- 
P(x,\eps,\eps{y_0}+ \eps^2 z_2)}{(z_1-z_2)}\biggr] (z_1-z_2)
\end{multline*}
et on constate que, si $\|z_k\|\leq \delta$ (on choisit ici $\delta$
pour que ce soit vrai) et pour $\eps$ assez petit, on se trouve devant
une combinaison linéaire de fonctions de $\HH$, avec $g$, $P$ et
$\Dp{P}{(\eps y)}$, multipliées par des fonctions bornées sur $\Dom$, telles
que $y_0$ et $z_1$, $z_2$. Cette combinaison linéaire est donc bornée,
et comme on peut mettre un $\eps$ en facteur, on a:
\begin{equation}\label{CH}
\|{\HHH}(z_1)-{\HHH}(z_2)\|_{\HH} \leq |\eps|\text{ C }
\|z_1-z_2\|
\end{equation}
Pour l'autre opérateur, 
\begin{align*}
{\ZZZ}(H_1)-{\ZZZ}(H_2)&=\frac{1}{\eps}\int_{\gamma_x} 
e^\frac{F(x)-F(t)}{\eps} \Bigl(H_1(t)-H_2(t)\Bigr) dt\\
\| {\ZZZ}(H_1)\!-\!{\ZZZ}(H_2) \| &\leq \f1\eps
\int_{\gamma_x}\!\!\!e^\frac{R(x)-R(t)}{\eps} \|H_1-H_2\|_{\HH}
\left|t^p f(t)\right| dt
\intertext{en utilisant la condition d'accessibilité du domaine $\Dom$,}
 & \leq \frac{\|H_1-H_2\|_{\HH}}{\eps}
\int_{\gamma_x} e^\frac{R(x)-R(t)}{\eps}\frac{1}{C_{\gamma_x}} 
\frac{-d}{dt} R(t) dt\\ 
 & \leq \|H_1-H_2\|_{\HH} \f{1}{C}
\Bigl[e^\frac{R(x)-R(t)}{\eps}\Bigr]^x_s
\end{align*}
Donc l'opérateur ${\ZZZ}$ est tel que
\begin{equation}\label{CZ}
\| {\ZZZ}(H_1)-{\ZZZ}(H_2) \| \leq
M. \|H_1-H_2\|_{\HH}
\end{equation}
où $M$ peut bien être choisi indépendant de $\eps$.

\medskip
Vérifions que les ensembles d'arrivée des deux opérateurs sont bien
ceux qui sont donnés.

Si $z(x,\eps)\in \ZZ$, 
\begin{equation}\label{eq:OH}
\left\|{\HHH}(z)(x,\eps)\right\|_{\HH}\leq |\eps|\cdot \|z\|
\cdot \|g\|_{\HH} + \left\|y_0+\eps z^2\right\|\cdot 
\left\|P(x,\eps,\eps\delta)\right\|_{\HH}\,,
\end{equation}
ceci, si pour tous les $z$ tels que $\|y_0+\eps z\|\leq \delta$
(d'ailleurs, l'opérateur n'est défini que pour de tels éléments de
$\ZZ$, et non sur tout l'espace). 

Si $H(x,\eps)\in \HH$,
\begin{align}
{\ZZZ}(H)&={\ZZZ}(0)+\left({\ZZZ}(H)-{\mathcal
Z}(0)\right)\notag\\ 
\left\|{\ZZZ}(H)\right\| & \leq \left\|{\ZZZ}(0)\right\|
+\left\|{\ZZZ}(H)-{\ZZZ}(0)\right\|\notag\\
 & \leq \f1\eps \left\| y_0g-y'_0+\Dp{h}\eps\right\|_{\HH} e^{\f{R(x)}\eps} 
\int_{\gamma_x} e^{-\f{R(t)}\eps}\left|t^p f(t)\right| dt
 +M \|H-0\|_{\HH}\notag\\
\left\|{\ZZZ}(H)\right\| & \leq  
M \left\| y_0g-y'_0+\Dp{h}\eps\right\|_{\HH}+M \|H\|_{\HH}\label{eq:OZ}
\end{align}
En effet, $y_0 g$, $y'_0$ et $\Dp{h}\eps$ sont dans l'espace $\HH$,
puisque $y_0$ et $y'_0$ sont bornées, et que $g\in\HH$ et $h\in \HH$.

\medskip
Donc, vu les inégalités \eqref{CH} et \eqref{CZ} le théorème du point
fixe est applicable dans ce cas pour l'opérateur ${\ZZZ} o
{\HHH}$ pour tout $\eps$ dans $S_0$, assez petit.


\bigskip

Nous venons de démontrer exactement qu'il existe une solution $y$ de
l'équation différentielle étudiée qui tend uniformément (pour $x\in
\Dom$) vers $y_0$ quand $\eps\rightarrow 0$ dans $S_0$. Cette solution
$y(x,\eps)$ est holomorphe en ses deux variables dans $\Dom\times
S_0$. En regardant un peu plus loin dans le développement de $y$ on
arrive au théorème suivant:
\begin{Th}\label{th:loin,c}
La solution $y$ du théorème~\ref{th:loin} admet la solution formelle
$\hat y =\sum y_n(x)\eps^n$ comme développement asymptotique, uniforme
pour $x\in\Dom$, quand $\eps$ tend vers $0$.
\end{Th}

Nous avons déjà montré que la vraie solution $y(x,\eps)$ existe et
qu'elle est bornée dans $\Dom$. Par hypothèse, on sait pouvoir
déterminer, pour tout $N$, les $N$ premiers coefficients de la série
formelle $\hat y$: $y_0(x),\ldots, y_{N-1}(x)$. Posons $R_N(x,\eps)$ tel
que 
$$y(x,\eps)=y_0(x)+\eps y_1(x)+\cdots+\eps^{N-1} y_{N-1}(x)
+\eps^N R_N(x,\eps)\,.$$ 
On sait de la fonction $R_N(x,\eps)$ qu'il s'agit, pour tout $\eps\in
S_0$ fixé, d'une fonction holomorphe et bornée en $x\in\Dom$, et que
$\eps^N R_N(x,\eps)$ est bornée, uniformément pour $x\in\Dom$, quand
$\eps$ tend vers $0$ dans $S_0$. On souhaite montrer en plus que
$R_N(x,\eps)$ est bornée sur $\Dom\times S_0$, en utilisant une 
équation différentielle pour $R$.

En remplaçant dans \eqref{eq:gen}, on obtient une équation
différentielle pour $R_N$, puisque tous les autres termes sont
connus. Dans le cas $N=1$, pour lequel on sait en plus que
$R_1(x,\eps)$ est bornée sur $\Dom\times S_0$,
\begin{align*}
\eps  \left(y_0+\eps R_1\right)'& = \left(x^P f+\eps g\right)
\left(y_0+\eps R_1\right) +h+\eps \left(y_0+\eps R_1\right)^2 
P(x,\eps,\eps y+\eps^2 R_1)\\
\eps R'_1&=\left(x^p f+\eps g\right) R_1
+y_0 g -y'_0+\f{h(x,\eps)-h(x,0)}\eps+ y^2 P(x,\eps,\eps y)\,,
\end{align*}
(on ne notera plus les dépendances en $x$ et/ou $\eps$ des fonctions
là où il n'y a pas d'ambigu\"ité)

On obtient bien une équation différentielle en $R_1(x,\eps)$, les
fonctions $y$ et $y_0$ en particulier pouvant \^etre considérées comme
données du problème. Cette équation différentielle est linéaire en
$R_1$, mais on va se contenter de vérifier que les fonctions
coefficients vérifient les propriétés du théorème \ref{th:loin}. 

Les fonctions $f$ et $g$ ne changent pas, leurs propriétés restent. La
nouvelle fonction pour $P$, $P_1$, est nulle. Le nouvel $h$ est
$$h_1(x,\eps)=y_0 g -y'_0+\f{h(x,\eps)-h(x,0)}\eps+ 
y^2 P(x,\eps,\eps y)\,,$$
une fonction holomorphe en $x$ sur $\Dom_0$, sauf éventuellement en
$x=0$, à cause de la présence de $h$ dans l'expression; elle est
aussi holomorphe en $\eps$ dans $S_0$, et reste bornée quand
$\eps\rightarrow 0$, d'après ce qu'on sait de $y$ et de la fraction en
$h$. Enfin, $h_1$ est bien dans l'espace $\HH$:
$$\left\|\f{h_1(x,\eps)}{x^p f}\right\|\leq \|y_0\| \cdot 
\left\|\f{g}{x^p f}\right\|+\left\|\f1{x^p f}\right\| \cdot \|y'_0\| +
\left\|\sup_\eps \Dp{h(x,\eps)}{\eps}\right\|_{\HH} + \|y\|^2 \cdot
\left\|\f{P(x,\eps,\eps y)}{x^p f(x)}\right\|\,,$$
puisque toutes les normes citées existent bien. 

On peut du coup appliquer l'intégralité du théorème \ref{th:loin}, la
fonction $R_1$ converge quand $\eps$ tend vers $0$, vers une fonction
$y_1(x)$, uniformément pour $x\in\Dom$: $R_1(x,\eps)=y_1(x)+\eps
R_2(x,\eps)$. La fonction $R_1$ est, pour tout $\eps$ assez petit,
bornée dans $\Dom$. 

Rien n'emp\^eche de refaire le m\^eme raisonnement pour $R_2$ à partir
de l'égalité $R_1=y_1+R_2$. On a donc, par récurrence directe,
$$\left\|y(x,\eps)-\sum_{i\leq N-1}y_i(x)\eps^i\right\|\leq 
\eps^N \|R_N(x,\eps)\,\|\ ,$$
où $R_N$ est, pour tout $N$, une fonction bornée en $\eps$
uniformément sur $\Dom$: $\hat y$ est un développement asymptotique de
$y$.

Pour \^etre complet dans la démonstration, il faut remarquer que dans
les hypothèses \HB, on demande que $h(x,\eps)$ admette un
développement asymptotique en $\eps$ quand $\eps$ tend vers $0$ dans
$S_0$, ce qui n'est pas le cas {\sl a priori} pour $h_1(x,\eps)$,
puisqu'on n'a pas d'abord cette propriété pour $y(x,\eps)$. Mais ce
n'est pas un problème: la démonstration du théorème \ref{th:loin}
demande seulement que l'on puisse écrire $h(x,\eps)=h(x,0)+\eps
h_\eps(x,\eps)$, avec $h_\eps$ bornée dans $\Dom\times S_0$, et cela,
on sait déjà pouvoir le faire pour $y(x,\eps)=y_0(x)+\eps
R_1(x,\eps)$. Ce qui suffit pour appliquer la récurrence.

\subsection{Au voisinage de $0$}\label{sec:pres}

Dans cette section, nous allons montrer qu'il existe des solutions près
du point tournant $0$, très exactement pour des $x$ dans des domaines
$$\left\{\varrho\,\eps^r< |x| \leq \rho,\ \arg(x)\in
[\theta_1,\theta_2]\right\}\ ,\ r=\f1{p+1}\,.$$ 
Ces solutions seront aussi appelées solutions intérieures.  Comme il
s'agit de prolonger les solutions trouvées au paragraphe précédent, on
conserve les hypothèses \HB, et on fixe à nouveau (au moins dans un
premier temps) un argument pour $\eps$, dans $S_0$.

Pour arriver au résultat de prolongement, nous regardons maintenant ce
qui se passe au voisinage de $x=0$, en utilisant une «loupe». Ce
changement de variable donne une équation non singulièrement perturbée
permettant effectivement d'analyser ce qui se passe au voisinage de
$0$. L'équation obtenue est appelée équation intérieure.

\bigskip
On pose $$Y=\eps y \hbox{ , et } X=\frac{x}{\eps'},$$ sachant que
$\eps'=\eps^{1/(p+1)}$, en prenant (par exemple) la valeur
principale de la racine. En partant de $\eps\in S_0$, on arrive ainsi
à $\eps'\in S'_0$.

L'équation \eqref{eq:gen} transformée s'écrit alors
\begin{equation} \label{Yg}
\frac{dY}{dX}=\Bigl(X^pf(\eps' X)+\eps'g(\eps'X,\eps)\Bigr) Y 
+\eps'h(\eps'X,\eps)+ \frac{1}{{\eps'}^p} Y^2 P(\eps'X,\eps,Y)\,.
\tag{{\sc ii}}
\end{equation}

On va se placer dans un domaine $\Dom_C$, qui contient une partie de
la montagne qui prolonge celle de $\Dom$ vers $0$ et de ses deux
vallées adjacentes, soit, à peu de choses près, l'intersection d'un
secteur ouvert centré en $0$, d'amplitude $\frac{3\pi}{p+1}$ et de
taille $\rho/\eps'$ et de la boule $|X|>\varrho=|X_l|$ (on se
contentera d'une constante $\varrho$ arbitrairement grande):
$$\Dom_C=\left\{X\in\CC\ \bigl/ \hbox{ }\varrho<|X|\leq\rho/\eps'\hbox{ et
}\arg X\! \in\, \Bigl]\theta,\theta+\frac{3\pi}{p+1}\Bigr[\, \right\}$$

\begin{figure}[ht!]
\caption{Domaine $\Dom_C$}
\centerline{\input domaine.pstex_t }
\label{DC}
\end{figure}

À $X$ donné, on définit le chemin $\Gamma_X$ comme étant le chemin
descendant le relief entre le point
$x_\rho/\eps'$ et $X$. $x_\rho$ est pris dans l'ensemble des points
\og au-dessus\fg de $\eps'X$ et de module $\rho$ (cf. aussi 
fig.~\ref{DC}).

On a défini le domaine $\Dom_C$, de telle manière que, pour tout
$X\in\Dom_C$, il existe un tel chemin $\Gamma_X$. 

Une hypothèse nécessaire pour continuer est que les trois limites
suivantes existent, au moins en tant que limites ponctuelles, pour
tout $X$ dans un domaine indépendant de $\eps$ incluant $\Dom_C$,
comme l'ensemble $\Dom_\infty=\left\{X\in\CC\ \bigl/ \hbox{ }|X|>|X_l|
\hbox{ et }\arg X\!  \in\, \Bigl]\theta,\theta+\frac{3\pi}{p+1}\Bigr[\,
\right\}$, et pour tout $Y$, $|Y|<\delta$, pour un certain $\delta$
fixé:
\begin{equation}\tag{$H_{\rightarrow 0}$}\label{H0}
\left. \begin{split}
\eps' g(\eps'X,{\eps'}^{p+1})&\underset{\eps'\rightarrow 0}{\longrightarrow} 
G_0(X)\\
\eps' h(\eps'X,{\eps'}^{p+1})&\underset{\eps'\rightarrow 0}{\longrightarrow} 
H_0(X)\\
P_\eps(X,\eps',Y)=\frac{P(\eps'X,{\eps'}^{p+1},Y)}{{\eps'}^p}
&\underset{\eps'\rightarrow 0}{\longrightarrow}  P_0(X,Y)
\end{split}\ \ \ \right\}
\end{equation}

Si on regarde ce que l'existence de telles limites signifie, on voit
que ces hypothèses sont équivalentes à: 
\begin{enumerate}
\item $g$ et $h$ ont, au pire, un pôle en $x=0$. Si $g$ a effectivement 
un pôle, alors
\begin{align*}
  g(x,\eps)&=g_0(\eps)+g_1(\eps) x +g_2(\eps) x^2+\ldots+g_n(\eps)
x^n+\ldots+\f{g_{-1}(\eps)}{x}\\ &
+\eps\!\left(\f{g_{-2}(\eps)}{x^2}+\ldots+
\f{g_{-p-2}(\eps)}{x^{p+2}}\right)\! 
 +\eps^2\!\left(\f{g_{-p-3}(\eps)}{x^{p+3}}+ \ldots+
\f{g_{-2p-3}(\eps)}{x^{2p+3}}\right)\!+\ldots\,, 
\intertext{en écrivant $g$ sous sa forme de série de Laurent en
    $x=0$ (en mettant $\eps$ en facteur là où il est nécessaire), dont
    tous les coefficients $g_k$, $k\in\Zed$ sont encore des fonctions
    holomorphes pour $\eps$ dans $S_0$, bornées quand $\eps\rightarrow
    0$; d'où, ensuite} 
G_0(X)&=\frac{\lim_{\eps\rightarrow 0}g_{-1}(\eps)}{X}+
\frac{\lim_{\eps\rightarrow 0}g_{-p-2}(\eps)}{X^{p+2}}+
\frac{\lim_{\eps\rightarrow 0}g_{-2p-3}(\eps)}{X^{2p+3}}
+\ldots
\end{align*}
Si au contraire $g$ est holomorphe en $x=0$,
alors $G_0(X)\equiv 0$. On obtient les mêmes résultats pour $h$ et
$H_0$. Il est clair que, dans tous les cas,
$\eps' g(\eps' X,\eps)$ converge uniformément vers $G_0(X)$, puisque
tous les termes des séries de Laurent restent bornés dans $\Dom_C$: 
$\eps'X$ et $1/X$ sont en effet bornés dans ce domaine.
\item $P(x,\eps,\eps y)$ peut s'écrire sous la forme $x^p P_1(x,
\eps y)+\eps P_2 (x,\eps,\eps y)$, et on trouve dans ce cas-là que
$$P_\eps(X,\eps',Y)=X^p P_1(\eps'X, Y)+\eps' P_2(\eps'X,\eps,Y)\,,$$
puis $P_0(X,Y)=X^p P_1(0,Y)$.
\end{enumerate}
\bigskip\bigskip

Pour commencer, nous allons montrer qu'il existe, pour $X$ assez grand, 
une solution $Y_0$, bornée dans un tel domaine $\Dom_C$, de l'équation
obtenue en posant $\eps'=0$ dans \eqref{Yg}
\begin{equation}
\frac{dY_0}{dX}=\Bigl(X^p f(0)+G_0(X)\Bigr) Y_0 +H_0(X)+Y_0^2
P_0(X,Y_0)\label{eq:Y_0g}
\end{equation} 
Il s'agit là d'une équation irrégulière singulière de rang $p+1$ (en
l'infini), dont on sait qu'elle a des solutions asymptotiques. On
redémontre ici sommairement ce résultat.

D'après leurs définitions, $G_0(X)$, $H_0(X)$ et
$P_0(X,\delta)/X^{p+1}$ tendent vers $0$ quand $X$ tend vers
l'infini. On peut donc montrer qu'il existe une solution formelle en
$1/X$ à cette équation\eqref{eq:Y_0g}: $\hat
Y_0(X)=-\frac{H_0(X)}{f(0) X^p}+ \cdots$. Nous allons montrer, avec un
point fixe, qu'il existe une et une seule solution à cette équation,
qui tende vers $0$ en l'infini, dans certains secteurs d'amplitude
inférieure à $3\frac{\pi}{p+1}$; la série formelle ci-dessus sera
série asymptotique à cette solution. On trouve d'abord que $Y_0$ est
solution de l'équation
$$Y_0(X)=\int_{\infty}^X \exp\left({\int_t^X \Bigl(u^p f(0)
+G_0(u)\Bigr) du}\right) \biggl(H_0(t)+Y_0(t)^2
P_0\left(t,Y_0(t)\right)\biggr) dt\,$$ 
ce qui définit un opérateur $\gY_0$ sur un sous-ensemble de l'espace
des fonctions bornées dans $\Dom_\infty$, qui tendent vers $0$ en
l'infini. Les images de telles fonctions par cet opérateur sont aussi
dans cet espace. Puis, si $Y_0^+$ et $Y_0^-$ sont deux fonctions de
cet espace on a
\begin{align*}
\begin{split}
\gY_0(Y_0^+)-\gY_0(Y_0^-)&
=\int_{\Gamma_X}  e^{\int_t^X\Bigl( u^p f(0) +G_0(u)\Bigr) du}
(Y_0^+-Y_0^-)\\
&\qquad \times\left((Y_0^++Y_0^-)P_0(t,Y_0^+)+Y_0^-
\frac{P_0(t,Y_0^+)-P_0(t,Y_0^-)}{(Y_0^+-Y_0^-)}\right) dt\\
\|\gY_0(Y_0^+)-\gY_0(Y_0^-)\| &\leq \|Y_0^+-Y_0^-\| \int_{\Gamma_X}
\exp\left({\Re\Bigl(\int_t^X \bigl( u^p f(0) +G_0(u)\bigr) du\Bigr)}\right)\\
&\qquad\times\biggl(2\delta \sup_{|\delta'|<\delta} |P_0(t,\delta')|+
\delta \sup_{|\delta'|<\delta} \Bigl|\frac{\partial P_0}{\partial
Y}(t,\delta')\Bigr|\biggr) dt
\end{split}
\end{align*}
ce qui montre bien, puisque $\frac{\partial P_0}{\partial
Y}(t,\delta)/t^p$ et $P_0(t,\delta)/t^p$ sont bornées, qu'on a un
opérateur contractant pour tous les $X$ accessibles, au moins pour
l'ensemble des fonctions $Y_0$ majorées par un certain $\delta$:
$$\delta\geq\|Y_0^\pm\|=\underset{|X|>\varrho}{\sup}|Y_0^\pm(X)|\,,$$
ce qui pour un $\delta$ fixé est vrai pour tous les $\varrho$ assez
grands, pusique les fonctions $Y_0$ tendent vers $0$. De m\^eme, si
$\|Y_0\|<\delta$, alors $\|\gY_0(Y_0)\|$ est inférieur à $\delta$ pour
un $\varrho$ assez grand, puisque $H_0$ tend aussi vers $0$ en l'infini.

Donc $Y_0$ existe bien dans un domaine de type $\Dom_\infty$.

\bigskip 

On remarque que si l'équation~\eqref{Yg} est effectivement non
singulièrement perturbée, on souhaite montrer un résultat d'existence
non trivial, puisqu'il s'agit de solutions dans un domaine non borné
(parce que dépendant de $\eps$), et avec des conditions initiales
dépendant de $\eps$ (puisqu'on veut prolonger des solutions du
paragraphe précédent).

Nous allons démontrer l'existence de solutions $Y$ de \eqref{Yg} pour
tous les $X$ accessibles à partir d'un point $x_\rho/\eps'$ {\em
fixé}. Nous les cherchons sous la forme $Y=Y_0+\eps' Z$. Alors
\begin{multline}\label{eq:enZ}
\frac{dZ}{dX}=\Bigl(X^p f(\eps'X) +\eps' g(\eps' X,\eps)\Bigr)
Z\\ +\frac{\eps'h(\eps'X,\eps)-H_0(X)}{\eps'}
+\left(\frac{\eps'g(\eps'X,\eps)-G_0(X)}{\eps'}+
X^p\frac{f(\eps'X)-f(0)}{\eps'}\right)Y_0\\ 
+Y_0^2\frac{{\eps'}^{-p} P(\eps'X,\eps,Y_0+\eps'Z) -P_0(X,Y_0)}{\eps'}
+\Bigl(2Y_0 Z+\eps'Z^2\Bigr)
\frac{P(\eps'X,\eps,Y_0+\eps'Z)}{{\eps'}^p}
\end{multline}
(on omet de rappeler systématiquement la dépendance de $Y_0$ en $X$ et
de $Z$ en $X$ et $\eps'$)

On pose
\begin{align*}
F_\eps(X)=&\int^X (f(\eps't)t^p+\eps'g(\eps't,\eps))dt\\
G(X,\eps')=&\left(\frac{\eps'g(\eps'X,\eps)-G_0(X)}{\eps'}+
X^p\frac{f(\eps'X)-f(0)}{\eps'}\right)Y_0
+\frac{\eps'h(\eps'X,\eps)-H_0(X)}{\eps'}\\
H(X,\eps',Z)=&Y_0^2\frac{P_\eps(X,\eps',Y_0+\eps'Z)-P_0(X,Y_0)}{\eps'}+
(2Y_0Z+\eps'Z^2)P_\eps(X,\eps',Y_0+\eps'Z)\\
            =&Y_0^2(X)X^p\f{P_1\bigl(\eps'X,Y_0(X)+\eps'Z\bigr)-
P_1\bigl(0,Y_0(X)\bigr)}{\eps'}\\ 
&\phantom{+}+P_2\bigl(\eps'X,\eps,Y_0(X)+\eps'Z\bigr) 
+\left(2Y_0(X)Z+\eps'Z^2\right)P_\eps\bigl(X,\eps',Y_0(X)+\eps'Z\bigr)
\end{align*}
Avec ces notations, on prend
\begin{equation}
Z_H(X,\eps')=\f1{\eps'}\bigl(\eps y(x_\rho,\eps) -Y_0(x_\rho/\eps')\bigr)+
\int_{\Gamma_X} e^{F_\eps(X)-F_\eps(t)} 
(G(t,\eps')+H(t,\eps',Z))dt\,,\label{eq:zgamma}
\end{equation}
où à $\eps'$ fixé, le premier terme correspond à une condition
initiale aux limites du domaine $\Dom_C$.

On va utiliser un théorème de point fixe avec des opérateurs
ressemblants à ceux du paragraphe précédent:
\begin{align*}
{\mathcal H}:\tilde{\Zz}\subset{\Zz}&\rightarrow{\Hh}\\
             Z& \mapsto {\mathcal H}(Z),
\end{align*}
{la fonction ${\mathcal H}(Z)$ étant telle que ${\mathcal
H}(Z)(X,\eps')=H\bigl(X,\eps',Z(X,\eps')\bigr)$ pour tout
$X\in\Dom_C$,} et
\begin{align*}
{\mathcal Z}:\tilde{\Hh}\subset{\Hh}&\rightarrow{\Zz}\\ 
             H& \mapsto Z_H
\end{align*}
la fonction $Z_H(X,\eps')$ étant définie par la relation
\eqref{eq:zgamma}. 

On peut noter que ${\mathcal H}(Z)$ contient ici encore quelques
termes linéaires en $Z$, qui resteront toutefois bornés par rapport à
ceux présents dans $F_\eps$.

Les deux opérateurs agissent à nouveau sur des sous-ensembles bornés
d'espaces $\Hh$ et $\Zz$. Ces espaces fonctionnels sont les suivants:
\begin{equation*}
  \Hh=\Biggl\{ H(X)\,/\,H \hbox{ est holomorphe pour tout }X\in\Dom_C,
  \hbox{ et } \sup_{X}\left|\f{H(X)}{X^p f(0)}\right|\hbox{ est
    borné}\Biggr\}
\end{equation*}
et
$$\Zz=\left\{ Z(X)\,/\,Z \hbox{ est holomorphe sur }\Dom_C ,
\hbox{ et tel que }
\sup_{X} {|Z(X)|}\hbox{ est borné}\right\} $$
(les $\displaystyle\sup_X$ sont bien sûr pris pour $X$ dans le
domaine $\Dom_C$).

On voudrait qu'avec les normes correspondantes $\|.\|_{\Hh}$ et  
$\|.\|$, l'inégalité suivante soit vraie: 
$$\|\HHH(Z_1)-\HHH(Z_2)\|_{\Hh}\leq C_{\HHH}\|Z_1-Z_2\|\,.$$ 
Or, si $\|Z_k\|\leq \delta$,
\begin{multline*}
|\HHH(Z_1)(X,\eps')-\HHH(Z_2)(X,\eps')|(X)\leq \|Z_1-Z_2\|\times
\bigl(|Y_0(X)|+|\eps'|\delta\bigr)\\
 \Biggl(\sup_{|\delta'|\leq \delta}\left|\Dp{P_\eps}{Y}
(X,\eps',Y_0+\eps'\delta')\right|
\left(|Y_0(X)|+|\eps'\delta'|\right)
+ 2\sup_{|\delta'|\leq \delta}
|P_\eps(X,\eps',Y_0+\eps'\delta')|\Biggr)
\end{multline*}
Une condition suffisante est que, pour tout $X$ dans $\Dom_C$ et
tout $\eps'$ assez petit la condition suivante est vérifiée:
\begin{multline}\label{HX}
\sup_{|\delta'|\leq \delta}
\left|\Dp{P_\eps}{Y}(X,\eps',Y_0+\eps'\delta')\right|
.\Bigl(|Y_0|+|\eps'\delta'|\Bigr) 
+ 2\sup_{|\delta'|\leq \delta}
\Bigl|P_\eps(X,\eps',Y_0+\eps'\delta')\Bigr|\\
\leq \text{M}_\varrho  \left|f(0)X^p\right|
\end{multline}
Elle l'est en particulier si
\begin{equation*}
\sup_{|\delta'|\leq \delta}\Bigl|P_\eps(X,\eps',Y_0+\eps'\delta')\Bigr|
\leq \text{M}'_\varrho  \Bigl|f(\eps'X)X^p\Bigr|
\end{equation*}
Il suffit de revenir à la variable $x$ pour s'apercevoir que cela
correspond à une des hypothèses de croissance comparée du
théorème~\ref{th:loin}:
$$\f{P(x,\eps,\eps y)}{f(x)x^p}=\f{P_1(x,\eps y)}{f(x)}
+\eps\f{P_2(x,\eps,\eps y)}{x^p f(x)}\ \text{ est bornée}\,.$$ 
Cette hypothèse reste vraie jusqu'en $x=0$ pour le premier terme
puisque $P_1(x,\eps y)/f(x)$ est holomorphe en ce point. En ce qui
concerne le second terme, il reste borné pour tout $x=\eps' C$, $C$
assez grand: $\eps/{x^p}$ est bornée pour les $x$ de cette forme, et
${P_2}/{f}$ est elle aussi holomorphe jusqu'en $x=0$.

On a alors
\begin{equation}\label{ineg:HHH}
\|\HHH(Z_1)-\HHH(Z_2)\|_{\Hh}\leq \Bigl(\|Y_0\|+|\eps'|\delta\Bigr) 
M_\varrho \|Z_1-Z_2\|\,.
\end{equation}

Par ailleurs,
\begin{align}
|\ZZZ(H_1)-\ZZZ(H_2)|(X,\eps') &\leq \int_{\Gamma_X} 
\left| e^{F_\eps(X)-F_\eps(t)} t^p f(0)\right|dt 
\times \|H_1-H_2\|_{\Hh}\,;\notag
\intertext{comme $\eps' g(\eps'X,\eps)$ est une fonction bornée sur 
$\Dom_C$ et que $|f(\eps't) t^p|= \gO\left(|t|^p\right)$, ce qui
implique $K^1_F |t|^p <|F'_\eps(t)|< K^2_F |t|^p$, on voit que}
\|\ZZZ(H_1)-\ZZZ(H_2)\|&\leq C_{\ZZZ}.\|H_1-H_2\|_{\Hh}\,.
\label{ineg:ZZZ} 
\end{align} 
Or ${\sup}_{|X|>\varrho}|Y_0(X)|+|\eps'\delta|$ peut être rendu aussi
petit qu'on veut quand $\varrho$ croît et que $\eps'$ est suffisamment
petit. Comme par ailleurs $C_{\ZZZ}$ est une constante qu'on peut rendre
indépendante de $\eps$ et $X$, et que $M_\varrho$ est
indépendante de $\eps$, les deux inégalités \eqref{ineg:HHH} et
\eqref{ineg:ZZZ} permettent d'envisager que l'opérateur $\ZZZ o \HHH$
soit contractant, au moins pour les $X$ dans un certain domaine
$\Dom_C\subset \{X\ /\ |X|>\varrho\}$.

\medskip
Il reste à vérifier que les images de ces opérateurs (restreints sur
des boules ouvertes de $\Hh$ ou $\Zz$) sont bien incluses dans les
espaces voulus.

Si $Z$ est une fonction de $\Zz$, dont la norme est majorée par un
certain $\delta>0$, la fonction image $H\bigl(t,\eps',Z(t,\eps')\bigr)$ 
est bien dans l'espace $\Hh$: 
\begin{list}{$\blacktriangleright$}{}
\item $P_\eps$ est dans $\Hh$ (puisque la condition \eqref{HX} est 
vérifiée).
\item $P_1(\eps'X,Y_0+\eps'Z)$ et $1/X^p P_2(\eps'X,\eps,Y_0+\eps'Z)$ sont 
holomorphes en toutes leurs variables, et donc bornées pour
$X\in\Dom_C$ et $\|Z\|<\delta$.
\item comme $P_1$ est holomorphe en $\eps'X$, y compris au voisinage de 
$0$, on peut écrire  $P_1(\eps'X,Y_0+\eps'Z)=\sum_{n\geq 0} (\eps'X)^n 
P_{1,n}(Y_0+\eps'Z)$, donc 
\begin{align*}
\f1{\eps'}&\Bigl(P_1(\eps'X,Y_0+\eps'Z)-P_1(0,Y_0)\Bigr)\\
&=\f1{\eps'} \left(P_{1,0}(Y_0+\eps'Z)-P_{1,0}(Y_0)
+\eps' X \sum_{n\geq 1}(\eps'X)^{n-1}P_{1,n}(Y_0+\eps'Z)\right)\\
&=Z\Dp{P_1}{Y}(0,Y_0)\Bigl(1+ o_{\eps'}(1)\Bigr)
+X \sum_{n\geq 1}(\eps'X)^{n-1}P_{1,n}(Y_0+\eps'Z)
\end{align*}
et puisque la somme dans cette expression est bornée par une constante
$C\|P_1\|$, et que la dérivée de $P_1$ l'est aussi, on peut majorer
cette différence, dans laquelle intervient toutefois à nouveau un facteur 
$X$.
\end{list}
Cela donne
\begin{multline*}
|H(X,\eps',Z)|\leq \|Y_0^2\|\cdot\delta \cdot 
\left\|\Dp{P_1}{Y}(0,Y_0)\right\|+ \|XY_0^2\|\cdot C\cdot\|P_1\| 
+\|Y_0^2\|\cdot\|P_2\|_{\Hh}\\
+(2\delta\cdot\|Y_0\|+\eps'\delta^2)\cdot
\|P_\eps(X,\eps',Y_0+\eps'Z)\|_{\Hh}\,.
\end{multline*}
On voit que non seulement $\HHH(Z)$ est dans l'espace $\Hh$, mais que
sa norme peut \^etre choisie aussi petite que l'on veut, si on prend
$\eps'$ assez petit et $\varrho$ assez grand, et ceci pour tout
$\delta$ fixé.

\medskip
Dans l'autre sens, regardons quelle est l'image de la fonction nulle
par $\ZZZ$. Cette fonction image est bornée si $G(X,\eps')$ est dans
$\Hh$ et si $\f1{\eps'}\bigl(\eps y(x_\rho,\eps)
-Y_0(x_\rho/\eps')\bigr)$ est bornée quand $\eps$ tend vers $0$.
\[G(X,\eps')=\left(\f{\eps'g(\eps'X,\eps)-G_0(X)}{\eps'}+
X^p\f{f(\eps'X)-f(0)}{\eps'}\right)Y_0
+\f{\eps'h(\eps'X,\eps)-H_0(X)}{\eps'}\,. \] 
Or, d'après ce qui a été dit sur les conséquences des hypothèses 
$\eqref{H0}$,
\begin{align*}
\eps'g(\eps'X,\eps)&= \eps'\sum_{k\geq 0} g_k(\eps)(\eps'X)^k
+\f{g_{-1}(\eps)}{X} \\
& \phantom{=} +\f{\eps}{{\eps'}^p}
\left({\eps'}^{p-1}\f{g_{-2}(\eps)}{X^2}+\ldots+ 
\eps'\f{g_{-p-1}(\eps)}{X^{p+1}}\right) 
 +\f{g_{-p-2}(\eps)}{X^{p+2}} +\ldots\\
\f{\eps'g(\eps'X,\eps)-G_0(X)}{\eps'}&=\sum_{k\geq 0}
g_k(\eps) (\eps'X)^k+\f{g_{-1}(\eps)-g_{-1}(0)}{\eps'X}\\
& \phantom{=}+
\left({\eps'}^{p-1} \f{g_{-2}}{X^2}+\ldots+ 
{\eps'}\f{g_{-p-1}}{X^{p+1}}\right) 
 +\f{g_{-p-2}(\eps)-{g_{-p-2}(0)}}{\eps'X^{p+2}} +\ldots
\end{align*}
Donc $\f{\eps'g(\eps'X,\eps)-G_0(X)}{\eps'}$, et de la m\^eme manière
$\f{\eps'h(\eps'X,\eps)-H_0(X)}{\eps'}$, sont des fonctions qui sont
bornées dans $\Dom_C$; elles sont dans l'espace $\Hh$, et leur norme
dans cet espace tend vers $0$ avec $\eps$.\label{PB}

En revanche, la fraction $\f{f(\eps'X)-f(0)}{\eps'}$ n'est pas
uniformément bornée pour $X\in\Dom_C$. On peut cependant démontrer le
lemme de majoration suivant:

\begin{Lemmetech}\label{lemme1}
Soit $F(x)$ une fonction holomorphe dans la boule ouverte
$\Dom_F={\mathcal B}(0,\rho)$. Alors, pour tout $X$ tel que
$(\eps'X)\in\Dom_F$,
$$\f{F(\eps'X)-F(0)}{\eps'}=X \tilde F(\eps'X)\,,$$ 
où $\tilde F$ est une fonction analytique; de plus, si $F$ est bornée,
$\tilde F$ est bornée elle aussi.
\end{Lemmetech}
On écrit d'abord $F(x)=\sum_{n\geq 0} f_n x^n$, puis on en tire
$$ \f{F(\eps'X)-F(0)}{\eps'} =
\sum_{n\geq 1} X f_n \times (\eps'X)^{n-1}\,.$$
On a donc $\tilde F(x)=\sum_{n\geq 1} f_n x^{n-1}$, qui est bien
analytique en $x$, et qui est bornée pour $x\in\Dom_F$ si $F$
l'est.~$\square$

\bigskip

Donc la fraction $\f{f(\eps'X)-f(0)}{\eps'}$  peut être d'ordre
$\gO(X)$; mais parce que $Y_0(X)$ tend vers $0$ quand
$X\rightarrow\infty$ au moins aussi vite que $1/X$, la fonction
$Y_0(X)\f{f(\eps'X)-f(0)}{\eps'}$ est, elle, bien bornée, et donc 
$X^p Y_0(X)\f{f(\eps'X)-f(0)}{\eps'}$ est une fonction de $\Hh$, elle
aussi. 

\smallskip
Il ne reste donc qu'une dernière condition à contrôler:
on veut montrer que la condition initiale
$Z(x_\rho/\eps',\eps')=\frac{1}{\eps'} \Bigl(\eps y(x_\rho,\eps)
-Y_0(x_\rho/\eps')\Bigr)$ est aussi petite que l'on veut pour tout
$\eps'$ assez petit.
Or $y(x_\rho,\eps)$ est borné en $\eps$, et $Y_0(x_\rho/\eps')$
est équivalent à 
$$Y_0(x_\rho/\eps')\underset{\eps'\rightarrow
0}\sim\f{H_0(x_\rho/\eps')}{f(0)\left(x_\rho/\eps'\right)^p} \sim
{\eps'}^p \f{H_0(x_\rho/\eps')}{f(0)\left(x_\rho\right)^p}\,,$$
sachant que $H_0(x_\rho/\eps')$ est au mieux un $\gO(\eps')$.  On est
donc assuré que $Z(x_\rho/\eps',\eps')$ est d'ordre au pire
${\eps'}^p$ quand $\eps'$ tend vers $0$, ce qui nous suffit largement.

\medskip
Donc $G(X,\eps')$ est une fonction de $\Hh$, on peut donc écrire
\begin{multline*}
\|\ZZZ(0)\|\leq 
\|Y_0\|\cdot\left\|\f{\eps'g(\eps'X,\eps)-G_0(X)}{\eps'}\right\|_{\Hh}+
\left\|X^p Y_0(X)\f{f(\eps'X)-f(0)}{\eps'}\right\|_{\Hh} \\
+\left\|\f{\eps'h(\eps'X,\eps)-H_0(X)}{\eps'}\right\|_{\Hh} + 
|Z(x_\rho/\eps',\eps')|
\end{multline*}
Tous ces termes sont bornés, pour tout $\eps'$ assez petit, et ils
décroissent avec $\varrho=|X_l|$.

\bigskip
On choisit donc $\delta\geq 2\|\ZZZ(0)\|$, puis on augmente la valeur
de $\varrho$ et on fait décro\^itre celle de $|\eps'|$ jusqu'à ce que 
$$C_{\ZZZ} M_\varrho (\|Y_0\|+\eps'\delta)<1\ \text{ et que }\ 
\sup_{\|Z\|<2\delta}\|H(X,\eps',Z)\|_{\Hh}<\delta/C_{\ZZZ}\,.$$ 
Dans ce cas, la composition des deux opérateurs $\HHH$ et $\ZZZ$ est
bien un opérateur contractant.

\medskip
On peut alors énoncer le théorème suivant.
\begin{Th}\label{th:pres}
On suppose que  pour l'équation~\eqref{eq:gen} les hypothèses
\eqref{H0} et celles du théorème~\ref{th:loin}, pour tout $r>0$, sont
vérifiées dans un $x$-domaine accessible (avec le relief
$\Re\frac{F(x)}{\eps}$), qu'on note à nouveau $\Dom\subset\Dom_0$.
Soit $X_l$ un réel assez grand, indépendant d'$\eps$ si
$|\eps|<\eps_0$. Soit alors pour un tel $\eps$ fixé des domaines
$$\Dom_C=\left\{X\in\CC \bigm/ \varrho<|X|\leq\rho/\eps'\hbox{ et
}\arg X\! \in\, \Bigl]\theta_1,\theta_2\Bigr[\,\right\}$$ 
{ $\theta_1$ et $\theta_2$ étant choisis, pour un $\eps$ donné, tels
que tout point $X$ de $\Dom_C$ soit accessible (avec le relief
$\Re\bigl(\eps' u^{p+1} f(0)\bigr)$) à partir d'au moins un point
de la forme $x/\eps'$, où $x\in\Dom$; et}
$$\Dom_{\eps'}=\Dom\,\cup\,\left\{x\in\CC \bigm/
(x/\eps')\in\Dom_C\right\}\,.$$ 

Alors il existe une solution holomorphe $y$ de \eqref{eq:gen} sur
$\Dom_{\eps'}$ telle que $\eps y$ est bornée quand $\eps$ tend vers
$0$, et $\eps y(X\eps',\eps)=Y(X,\eps')$ tend, uniformément pour $X$
dans $\Dom_C$ quand $\eps$ tend vers $0$, vers une fonction $Y_0(X)$
solution de l'équation \eqref{eq:Y0}.
\end{Th}
\bigskip

On sait qu'il existe une unique solution au problème de Cauchy
suivant: $Y(X,\eps')$ est solution de \eqref{Yg} et vérifie
$Y(x_\rho/\eps',\eps')=\eps y(x_\rho,\eps)$. On vient de montrer que
cette solution existe sur la partie du domaine $\Dom_C$ accessible à
partir du $x_\rho/\eps'$ choisi et qu'elle reste dans un certain
voisinage de $Y_0(X)$. Or, par unicité des solutions avec la condition
initiale au point $x_\rho$, la solution de \eqref{eq:gen} $y(x,\eps)$
pour $x\in\Dom, |x|\geq \rho$ se prolonge avec la fonction
$\frac{1}{\eps}Y(x/\eps',\eps)$ dans toute la partie de $\Dom_C$
accessible; et comme tout $\Dom_C$ est accessible à partir des
$x_\rho/\eps'$ (et m\^eme d'un nombre fini de tels points), et qu'il
ne peut y avoir qu'un seul prolongement dans ce domaine, $y(x,\eps)$
se prolonge en fait avec $\frac{1}{\eps}Y(x/\eps',\eps)$ pour
l'ensemble des $x\in \eps'\Dom_C$.

Ce qui clôt la démonstration du théorème.

\bigskip
Cette solution $y(x,\eps)$ est holomorphe en ses deux variables dans
$\Dom_{\eps'}\times S_0$. Pour démontrer la dépendance holomorphe en
$\eps'$, puis en $\eps$, de la fonction $Y(X,\eps)$ correspondante, on
procède comme pour le théorème \ref{th:loin}, en prenant un
$\Dom_\eps'$ accessible pour un intervalle d'arguments de $\eps'$ et en
considérant des espaces de fonctions sur $\Dom_C\times S_0$. On ne
détaillera pas la démonstration du théorème \ref{th:pres} modifiée
ici.

\subsection{Développement de la solution intérieure}

La solution $Y(X,\eps')$ trouvée tend vers $Y_0(X)$ quand
$\eps\rightarrow 0$, uniformément pour $X\in\Dom_C$. Admet-elle un
développement asymptotique de la forme $\hat Y=\sum_{n\geq 0}
Y_n(X){\eps'}^n$, uniformément dans tout le domaine $\Dom_C$  ?

La réponse n'est pas toujours positive. Si on tente de poser
$$Y(X,\eps')=\sum_{0\leq n\leq N} Y_n(X){\eps'}^n+{\eps'}^{N+1} \tilde
Y(X,\eps')$$ 
et de montrer que $\tilde Y$ est bornée de la même manière qu'on a
montré que $Z$ était bornée (théorème \ref{th:pres}), on se retrouve
face à un problème: au cours de la démonstration, on devait montrer,
page~\pageref{PB}, que $G$ était dans $\Hh$, ce qui n'est vrai que
parce que $Y_0$ tend vers $0$ en l'infini. Mais on peut seulement
montrer en général que $Z$ est bornée dans $\Dom_C$; une démonstration
par récurrence générale, valable sur tout le domaine $\Dom_C$ comme
dans le cas de l'équation de Van der Pol, est donc exclue.

\medskip
Par contre, on peut arriver à démontrer le corollaire suivant.
\begin{Cor}\label{cor:pabo}
  La fonction $Y(X,\eps')$ du théorème \ref{th:pres} vérifie la
  propriété suivante. Soit $\mu\in ]0,1[$, et $\varrho>0$ un réel
  assez grand avec le domaine $\Dom_C$ qui en dépend de la manière
  habituelle . Il existe une suite de fonctions $Y_0$, $Y_1$, \ldots,
  $Y_n$, \ldots telle que, pour tout $N$,
\[Y(X,\eps')-\sum_{n=0}^{N-1} Y_n(X){\eps'}^n={\eps'}^{N}
Z_N(X,\eps')\,,\] 
la fonction $Z_N(X,\eps')$ étant définie pour $X\in\Dom_C$, mais
bornée par ${\eps'}^{-(N-1)\mu}$ quand $\eps'$ tend vers $0$ et
seulement pour les $X$ tels que $\varrho<X<|\eps'|^{-\mu}$ (on notera
$\Dom_{C\mu}$ l'ensemble des $X$ bornés de cette manière et qui sont 
aussi dans $\Dom_C$).
\end{Cor}
Nous démontrerons cette propriété par récurrence sur $N$, en
réutilisant une partie de la démonstration du théorème précédent. On y
traite de fait le cas $N=0$: $Z_1(X,\eps')$ est une fonction bornée,
en particulier pour les $X<|\eps'|^{-\mu}$.

Mais auparavant, nous allons montrer certaines propriétés de
majoration. On introduit dans ce but, pour tout $m\in \Zed$ la classe de
fonctions suivante: $F(X,\eps',y,y_0) \in \gF_m$ si $F$ s'écrit comme une
somme finie telle que
$$F(X,\eps',y,y_0)=\sum_l E_l(X) F_l(\eps'X,\eps',y,y_0)$$
où
\begin{list}{$\blacktriangleright$}{}
\item $F_l(x,\eps',y,y_0)$ est une fonction bornée, holomorphe au 
voisinage de $x=0$, pour $\eps'$ dans des secteurs $S'_0$ et pour
$|y|\leq \delta$, $|y_0|\leq \delta$; elle admet un développement
asymptotique uniforme quand $\eps'\rightarrow0$.
\item $X^{-m} E_l(X)$ est holomorphe et bornée indépendamment de
$\eps'$ dans $\Dom_{C\mu}$.
\end{list}
Comme propriété de ces classes de fonctions, on peut noter que, si
$F\in\gF_m$ et $G\in \gF_n$, alors
\begin{list}{$\star$}{}
\item $F+G$ est dans $\gF_{\sup(m,n)}$.
\item $F\times G$ est dans $\gF_{m+n}$.
\item $\Dp{F}{Y}$ est une fonction de $\gF_m$, quitte à diminuer $\delta$.
\item $\eps'F\in \gF_{m-1}$.
\end{list}
Nous allons démontrer le lemme suivant:
\begin{Lemmetech}\label{lemme2}
Soit $F(X,\eps',y,y_0)$ une fonction de $\gF_m$. La fonction $\tilde
F$ définie comme
$$\tilde F: (X,\eps',y,y_0) \mapsto 
\f{F(X,\eps',y,y_0)-F(X,0,y,y_0)}{\eps'}$$
est une fonction de la classe $\gF_{m+1}$.
\end{Lemmetech}
La démonstration ressemble à celle du lemme~\ref{lemme1}. On écrit
$F_l$ comme la somme convergente (pour $\eps'X$ assez petit)
$$F_l(\eps'X,\eps',y,y_0)=\sum_{k\geq 0} (\eps'X)^k
F_{l,k}(\eps',y,y_0)\,.$$ 
On arrive alors à l'expression suivante (où la $\sum_k$ reste
convergente si $\eps'X$ est assez proche de $0$):
\begin{align*}
\tilde F (X,\eps',y,y_0) &= \sum_l E_l\times\! \left( \sum_{k\geq 1} X
F_{l,k}(\eps',y,y_0)
(\eps'X)^{k-1}\!+\f{F_{l,0}(\eps',y,y_0)-F_{l,0}(0,y,y_0)}{\eps'}\!
\right)\!;\\ 
\intertext{comme $F_{l,0}$ admet un développement asymptotique
uniforme en $\eps'=0$,}
\tilde F (X,\eps',y,y_0) &= \sum_l E_l(X)X 
\tilde F_{l}(\eps'X,\eps',y,y_0)+\sum_l E_l(X) 
\f{F_{l,0}(\eps',y,y_0)-F_{l,0}(0,y,y_0)}{\eps'}\,.
\end{align*}
$|X^{-m-1} E_l(X) X|$  et $|X^{-m-1} E_l(X)|$ étant bornées dans
$\Dom_{C\mu}$, la fonction $\tilde F$ est donc dans 
$\gF_{m+1}$.~$\square$

\medskip
Pour $N=1$, on reprend l'équation pour $Z$ \eqref{eq:enZ} vue plus
haut, qui est vérifiée par $Z_1$:
\begin{align*}
\f{dZ}{dX}&=f_\eps(X) Z+G(X,\eps')+H(X,\eps',Z)
\intertext{avec}
f_\eps(X)=& f(\eps'X)X^p+\eps'g(\eps'X,\eps)\\
G(X,\eps')=&\left(\frac{\eps'g(\eps'X,\eps)-G_0(X)}{\eps'}+
X^p\frac{f(\eps'X)-f(0)}{\eps'}\right)Y_0
+\frac{\eps'h(\eps'X,\eps)-H_0(X)}{\eps'}\\
H(X,\eps',Z)=&Y_0^2 \f{P_\eps(X,\eps',Y_0+\eps'Z)-P_0(X,Y_0)}{\eps'}
+(2Y_0Z+\eps'Z^2)P_\eps(X,\eps',Y_0+\eps'Z)\\
 =& Y_0^2 \f{P_\eps(X,\eps',Y_0+\eps'Z)-P_0(X,Y_0)}{\eps'} +(Y_0+Y)Z
P_\eps(X,\eps',Y_0+\eps'Z)
\end{align*}
Pour simplifier au maximum le terme $H$, nous allons réorganiser
l'équation, en utilisant le fait que $Y$ est maintenant connue, comme
fonction holomorphe bornée. On reprend la décomposition suivante de
$P_\eps$ donnée parmi les hypothèses:
\begin{align*}
P_\eps(X,\eps',Y)&=X^p P_1(\eps'X,Y)+\eps' P_2(\eps'X,\eps,Y)\,,
\intertext{puis}
X^p P_1(\eps'X,Y)-P_0(X,Y_0)&=X^p \left( P_1(\eps'X,Y)-P_1(0,Y)+
P_1(0,Y)-P_1(0,Y_0)\right)\\
        &=X^p \left( P_1(\eps'X,Y)-P_1(0,Y)\right)\\
&\qquad\qquad+X^p (Y-Y_0)\int_0^1 \Dp{P_1}{Y}\Bigl(0,Y_0+t(Y-Y_0)\Bigr) dt\,;
\end{align*}
ce qui donne
\begin{multline*}
\f{P_\eps(X,\eps',Y_0+\eps'Z)-P_0(X,Y_0)}{\eps'}=X^p
\f{P_1(\eps'X,Y)-P_1(0,Y)}{\eps'}+P_2(\eps'X,\eps,Y)\\
+X^p Z \int_0^1 \Dp{P_1}{Y}\Bigl(0,Y_0+t(Y-Y_0)\Bigr) dt\,.
\end{multline*}
On peut obtenir ainsi comme équation différentielle pour $Z_1$:
\begin{align*}
\f{dZ_1}{dX}&=f_1\Bigl(X,\eps',Y(X,\eps'),Y_0(X)\Bigr) Z_1+
G_1\Bigl(X,\eps',Y(X,\eps'),0\Bigr)
\intertext{avec}
f_1(X,\eps',y,y_0)=& f(\eps'X)X^p+\eps'g(\eps'X,\eps)+
\Bigl(Y_0(X)+y\Bigr) P_\eps(X,\eps',y)\\
        &\qquad\qquad+X^pY_0^2(X) \int_0^1 \Dp{P_1}{Y}
\Bigl(0,y_0+t\bigl(y-y_0\bigr)\Bigr) dt\\
G_1(X,\eps',y,y_0)=&\left(\!\frac{\eps'g(\eps'X,\eps)\!-G_0(X)}{\eps'}+\!
X^p\frac{f(\eps'X)\!-f(0)}{\eps'}\right)\!Y_0(X)\\
&\qquad+\frac{\eps'h(\eps'X,\eps)-H_0(X)}{\eps'}\\
&\qquad\qquad+Y_0^2(X) P_2(\eps'X,\eps,y)
+Y_0^2(X) X^p \f{P_1(\eps'X,y)-P_1(0,y)}{\eps'}
\end{align*}
{C'est à partir de cette équation devenue linéaire en $Z$ que nous
ferons la démonstration par récurrence. }

\smallskip\pagebreak[1]

D'après ce qui a été vu plus haut ces fonctions vérifient les propriétés
suivantes:
\begin{list}{$\blacktriangleright$}{}
\item $f_1(X,\eps',y,y_0)$ est de classe $\gF_p$.
\item $G_1(X,\eps',y,y_0)$ est aussi de la m\^eme classe $\gF_p$.
\item $Z_1$ est une fonction bornée dans $\Dom_\mu$, donc bornée dans
$\Dom_{C\mu}$. 
\end{list}

On note $Y_1(X)$ la solution de l'équation de $Z_1$ avec $\eps'=0$ qui
est bornée dans $\Dom_\infty$ ($Y_1\in \gF_0$); on montre l'existence
de cette solution comme on a montré celle de $Y_0$, sachant que
$X^{-p}f_1(X,0,Y_0(X),Y_0(X))$ et $X^{-p}G_1(X,0,Y_0(X),0)$ sont deux 
fonctions bornées dans $\Dom_{C\mu}$.

On pose ensuite
$$Z_1(X,\eps')=Y_1(X)+\eps' Z_2(X,\eps')\,.$$

\medskip

Nous aurons donc besoin de compléter le lemme~\ref{lemme2}, en
regardant ce qu'on peut en déduire pour les classes $\gG_m^k(Z)$
comprenant les fonctions $G(X,\eps',y,y_0,Z)$, où $G$ {est un polyn\^ome
en} $Z$ dont les coefficients $G_l$ de $Z^l$ {sont dans} $\gF_{k-lm}$.

L'intér\^et principal de ces classes $\gG$ de fonctions est que, si on
sait que $X^{-m}Z(X,\eps')$ est une fonction bornée dans
$\Dom_{C\mu}$, alors
$G\Bigl(X,\eps',Y(X,\eps'),Y_0(X),Z(X,\eps')\Bigr)$ cro\^it au plus
comme $X^{p+n-1}$ avec $X$ dans ce domaine.

\begin{Lemmetech}\label{bourrin1}
Soit $G(X,\eps',y,y_0,Z)$ une fonction de
$\gG_{m-2}^k(Z)$. Alors on a les propriétés suivantes:
\begin{enumerate}
\item $G$ multiplié par une fonction de $\gF_l$ est dans 
$\gG_{m-2}^{k+l}(Z)$.\label{p1}
\item Le produit de $G$ par $Z$ est dans 
$\gG_{m-2}^{k+m-2}(Z)$.\label{p2}
\item $\Dp{G}{Z}(X,\eps',y,y_0,Z)
\in\gG_{m-2}^{k-m+2}(Z)$.\label{p3} 
\item Une fonction $G(X,\eps',y,y_0,Z_{m-1})$ de $\gG_{m-2}^k(Z_{m-1})$ 
peut se transformer en posant $Z_{m-1}=Y_{m-1}(X)+\eps' Z_m$. La fonction
$\tilde G(X,\eps',y,y_0,Z_{m})$ obtenue est alors dans 
$\gG_{m-1}^k(Z_{m})$.\label{p4}
\end{enumerate}
\end{Lemmetech}
\begin{enumerate}
\item Ce premier point est trivial.
\item Pour démontrer cette propriété, on part de la définition:
\begin{align*}
G(X,\eps',y,y_0,Z)&=\sum_l G_l {Z}^l, \qquad\qquad
G_l(X,\eps',y,y_0)\in\gF_{k-l(m-2)}\\ G(X,\eps',y,y_0,Z)Z&=\sum_l
G_l {Z}^{l+1} =\sum_l G_{l-1} {Z}^{l},
\end{align*}
{et $G_{l-1}\in\gF_{k-(l-1)(m-2)}=\gF_{k+m-2-l(m-2)}$}
\item On procède comme ci-dessus:
\begin{align*}
\Dp{G}{Z}(X,\eps',y,y_0,Z)&=\sum_l G_l l {Z}^{l-1}
=\sum_l (l+1) G_{l+1} {Z}^{l}
\end{align*}
avec $(l+1) G_{l+1}\in \gF_{k-(l+1)(m-2)}=\gF_{k-m+2-l(m-2)}$. D'où la
propriété annoncée.
\item La démonstration ici se fait en écrivant
$Z_{m-1}=Y_{m-1}+\eps' Z_m$, où $Y_{m-1}\in \gF_{m-2}$.  On obtient
alors un polyn\^ome en $Z_m$ dont le terme en ${Z_{m}}^l$ est une
somme de termes de la forme
$$G_{j} \text{C}_j^{l} Y_{m-1}^{j-l}{\eps'}^l {Z_{m}}^l\,,$$
où $C^l_j$ désigne le nombre de combinaisons, et $G_j$ est par définition 
dans $\gF_{k-j(m-2)}$;
ce qui fait qu'on a en facteur de ${Z_{m}}^l$ une fonction dans
$\gF_{k-j(m-2)+(m-2)(j-l)-l}=\gF_{k-l(m-1)}$, pour tout $l$.  On
remarque en passant que $G(X,\eps',y,y_0,Y_{m-1}(X))$ peut se voir
comme un élément de $\gF_{k}$.
\end{enumerate}

Ces propriétés servent essentiellement à démontrer le lemme suivant:
\begin{Lemmetech}\label{bourrin2}
Soit $G\in \gG_{m-2}^k(Z)$ et $H$ la fonction définie par 
$$H(X,\eps')=G\Bigl(X,\eps',Y(X,\eps'),Y_0(X),Z_{m-1}(X,\eps')\Bigr).$$ 
Alors il existe $\breve G\in\gG_{m-1}^{k+1}(Z_m)$ telle que
\begin{equation*}
\f1{\eps'}\Bigl(H(X,\eps')-H(X,0)\Bigr)=
\breve G\Bigl(X,\eps',Y(X,\eps'),Y_0(X),Z_m(X,\eps')\Bigr)\,.
\end{equation*}
\end{Lemmetech}
On décompose
\begin{multline}\label{eq:Gm}
\breve G(X,\eps',Y,Y_0,Z_m)=
\f{G\left(X,\eps',Y,Y_0,Z_{m-1}\right) 
-G\left(X,0,Y,Y_0,Z_{m-1}\right)}{\eps'}\\
+\f{G\left(X,0,Y(X,\eps'),Y_0(X),Z_{m-1}\right) 
-G\left(X,0,Y(X,0),Y_0(X),Z_{m-1}\right)}{\eps'}\\
+\f{G\left(X,0,Y_0,Y_0,Z_{m-1}(X,\eps')\right) 
-G\left(X,0,Y_0,Y_0,Z_{m-1}(X,0)\right)}{\eps'}
\end{multline}
Dans cette somme,
\begin{enumerate}
\item Le premier terme correspond, si on utilise le lemme \ref{lemme2}, 
à un polyn\^ome en $Z_{m-1}$ dont le coefficient $H_l(X,\eps',Y,Y_0)$ pour
${Z_{m-1}}^l$ est dans $\gF_{k+1-l(m-2)}$. Si on remplace ensuite
${Z_{m-1}}^l$ par $(Y_{m-1}+\eps' Z_m)^l$, on retrouve le résultat de
la propriété \ref{p4} ci-dessus.  \par
Le premier terme est donc dans $\gG_{m-1}^{k+1}(Z_m)$.
\item On peut écrire la deuxième fraction de \eqref{eq:Gm}, à l'aide 
de la formule de Taylor avec reste intégral
$$\!\!\f{Y(X,\eps')-Y(X,0)}{\eps'} \int_0^1\!
\Dp{G}{Y}\Bigl(X,0,\bigl[Y_0+t(Y\!-\!Y_0)\bigr](X,\eps'),Y_0,
Z_{m-1}(X,\eps')\Bigr) dt\,.$$ 
L'intégrale st holomorphe en $Y$ et $Y_0$, et c'est aussi un
polyn\^ome en $Z_{m-1}$ qui reste dans $\gG_{m-2}^k(Z_{m-1})$. On le
multiplie par
\begin{align*}
\f{Y(X,\eps')-Y(X,0)}{\eps'}&=Z_1(X,\eps')\\
&=Y_1(X)+\eps' Y_2(X)+\ldots
+{\eps'}^{m-2}Y_{m-1}(X)+{\eps'}^{m-1}Z_{m-1}\,;
\end{align*} 
puisque ${\eps'}^{j-2}Y_{j-1}(X) \in \gF_0$ pour tout $j$, avec les
points \ref{p1} et \ref{p2} du lemme, on sait que le produit est dans 
$$\gG_{m-2}^{k}(Z_{m-1}) \cup \gG_{m-2}^{k-(m-1)+m-2}(Z_{m-1}) 
\subset \gG_{m-2}^{k}(Z_{m-1}).$$
La deuxième fraction est donc (largement) dans $\gG_{m-1}^{k+1}(Z_m)$.
\item La troisième fraction se réécrit aussi avec une formule de 
Taylor
\begin{multline*}
\f{Z_{m-1}(X,\eps')-{Z_{m-1}(X,0)}}{\eps'}\times\\
\int_0^1 \Dp{G(X,0,Y,Z)}{Z}\Bigl(X,0,Y_0,Y_0,
\left[Z_{m-1}+t(Z_{m-1}-Y_{m-1})\right](X,\eps')\Bigr) dt\,.
\end{multline*}
On intègre un polyn\^ome en $t$ dont les coefficients sont, d'après la
propriété \ref{p3} du lemme précédent, des éléments de
$\gG_{m-2}^{k-m+2}(Z_{m-1})$, ou (avec la propriété \ref{p4}) de
$\gG_{m-1}^{k-m+2}(Z_m)$. L'intégrale est donc dans
$\gG_{m-1}^{k-m+2}(Z_m)$. On la multiplie par la fraction, égale à
$Z_m$. On se retrouve donc dans $\gG_{m-1}^{k-m+2+(m-1)}=
\gG_{m-1}^{k+1}$.$\square$
\end{enumerate}

On fait comme hypothèses de récurrence que pour $m\leq n-1$, $Z_{m}$
est solution de l'équation différentielle
$$\f{dZ_{m}}{dX}=f_1\left(X,\eps',Y(X,\eps'),Y_0(X)\right)Z_{m}+
G_{m}\left(X,\eps',Y(X,\eps'),Y_0(X),Z_{m-1}(X,\eps')\right)$$ 
où $G_{m}(X,\eps',y,y_0,Z)$ est un élément de
$\gG_{m-2}^{p+m-1}(Z_{m-1})$. 
La croissance de la fonction $Z_{m}$ est majorée:
$X^{-m+1}Z_{m}(X,\eps')$ est bornée dans $\Dom_{C\mu}$. Il existe
alors pour tout $m\leq n-1$ une unique fonction $Y_{m}$ solution de
$$\f{dY_m}{dX}=f_1\Bigl(X,0,Y_0(X),Y_0(X)\Bigr)Y_{m}+
G_{m}\Bigl(X,0,Y_0(X),Y_0(X),Y_{m-1}\Bigr)$$
telle que $X^{-m+1}Y_{m}(X)$ est bornée ($Y_m \in \gF_{m-1}$).

\smallskip
Avec ces hypothèses, on pose 
$$Z_{n-1}(X,\eps')=Y_{n-1}(X)+\eps' Z_n(X,\eps')\,.$$
On remplace $Z_{n-1}$ par cette expression dans son équation
différentielle (on ne notera plus les dépendances de $Y(X,\eps')$ et 
$Y_0(X)$): 
\begin{multline}\label{eq:Zn}
\f{dZ_{n}}{dX}=
f_1(X,\eps',Y,Y_0)Z_{n}+
\f{f_1(X,\eps',Y,Y_0)-f_1(X,0,Y_0,Y_0)}{\eps'}Z_{n-1}\\
+\f{G_{n-1}(X,\eps',Y,Y_0,Z_{n-2})-G_{n-1}(X,0,Y_0,Y_0,Y_{n-2})}{\eps'}\,.
\end{multline}
On veut étudier la croissance de la fonction $G_n$ où l'on pose
\begin{multline}
G_n(X,\eps',Y,Y_0,Z_{n-1})=
\f{f_1(X,\eps',Y,Y_0)-f_1(X,0,Y_0,Y_0)}{\eps'}Z_{n-1}\\ 
+\f{G_{n-1}(X,\eps',Y,Y_0,Z_{n-2})-G_{n-1}(X,0,Y_0,Y_0,Y_{n-2})}{\eps'}
\label{eq:Gn}
\end{multline}

\medskip 

L'hypothèse de récurrence est donc que $G_{n-1}\in
\gG_{n-3}^{p+n-2}(Z_{n-2})$.

Si on applique le lemme \ref{bourrin2} à
$\f{f_1(X,\eps',Y,Y_0)-f_1(X,0,Y_0,Y_0)}{\eps'}$, en considérant $f_1$
comme élément de $\gG_{0}^{p}(Z_1)$, on voit que cette fraction est
dans la classe $\gG_{n-2}^{p+1}(Z_{n-1})$; en multipliant encore par
$Z_{n-1}$, on se retrouve, d'après le lemme
\ref{bourrin1} (propriété \ref{p2}), dans $\gG_{n-2}^{p+n-1}(Z_{n-1})$.

Quant à
$\f{G_{n-1}(X,\eps',Y,Y_0,Z_{n-2})-G_{n-1}(X,0,Y_0,Y_0,Y_{n-2})}{\eps'}$,
le lemme \ref{bourrin2} nous indique qu'il s'agit d'une fonction de
$\gG_{n-3+1}^{p+n-2+1}(Z_{n-1})$, soit $\gG_{n-2}^{p+n-1}(Z_{n-1})$.

Par conséquent, $G_n(X,\eps',Y,Y_0,Z_{n-1})\in
\gG_{n-2}^{p+n-1}(Z_{n-1})$.

On applique maintenant la formule de variation de la constante pour
l'équation différentielle \eqref{eq:Zn} concernant $Z_n$:
\begin{equation}\label{eq:Znvar}
Z_n(X,\eps')= Z_n\left(x_\rho{\eps'}^{\mu-1},\eps'\right) +
\int_{\Gamma^\mu_X} \!\! e^{F_1(X)-F_1(t)}
G_n\Bigl(t,\eps',Y(t,\eps'),Y_0(t),Z_{n-1}(t,\eps')\Bigr) dt\,,
\end{equation}
avec $F_1(t)=\int^t f_1\bigl(u,\eps',Y(u,\eps'),Y_0(u)\bigr)du$.
Le chemin $\Gamma^\mu_X$
est lui déduit du chemin $\Gamma_X$, à l'intérieur du domaine
$\Dom_{C\mu}$, descendant entre un point $x_\rho{\eps'}^{\mu-1}$ (au
lieu de $x_\rho{\eps'}^{-1}$ pour $\Gamma_X$) et $X$.

On en déduit que
\begin{align}
|Z_n(X,\eps')| &\leq
\sup_{X\in\Dom_{C\mu}} \left|\f{G_n\Bigl(X,\eps',Y(X,\eps'),Y_0(X),Z_{n-1}(X,\eps')\Bigr)}{f_1\Bigl(X,\eps',Y(X,\eps'),Y_O(X,\eps')\Bigr)}\right|\notag\\
&\leq \sup_{X\in\Dom_{C\mu}} |X|^{(n-1)}+
\left|Z_n(x_\rho{\eps'}^{\mu-1},\eps')\right|\label{ineq:Zn} \\
&\leq |\eps'|^{-(n-1)\mu}\,.\quad\square\notag
\end{align}

Cependant, cette démonstration doit encore être complétée: il reste à
montrer que la «condition initiale» que l'on choisit égale à
$$Z_n\left(x_\rho{\eps'}^{\mu-1},\eps'\right)=\f1{{\eps'}^n}\Bigl(\eps
y(x_\rho{\eps'}^\mu,\eps) -Y_0(x_\rho{\eps'}^{\mu-1})-
\ldots-{\eps'}^{n-1}Y_{n-1}(x_\rho{\eps'}^{\mu-1})\Bigr)$$ 
pour assurer la continuité entre solution intérieure et extérieure
n'est pas plus grande que les autres termes de l'inégalité
\eqref{ineq:Zn}. Cette propriété sera démontrée au paragraphe suivant 
(corollaire \ref{th:connexion}); pour le montrer il faudra chercher à
étendre le domaine de la solution extérieure.

On remarque que l'on construit avec cette démonstration du corollaire
une solution formelle $\hat Y(X,\eps')$ pour l'équation \eqref{Yg},
qui a une propriété ressemblant à celle d'une série asymptotique pour
la vraie solution $Y(X,\eps')$. Mais si on ne contr\^ole leur
croissance de manière intéressante que dans $\Dom_{C\mu}$, tous les
coefficients $Y_n(X)$ existent dans le m\^eme domaine que $Y_0(X)$,
c'est-à-dire pour $X\in\Dom_\infty$.

\bigskip
Comme l'équation~\eqref{Yg} n'est pas singulièrement perturbée, on
peut encore compléter le théorème~\ref{th:pres} et son premier
corollaire avec le corollaire suivant: 
\begin{Cor}[Principe de prolongement analytique]\label{th:+pres}
En plus de l'ensemble des hypothèses du théorème~\ref{th:pres},
supposons que la fonction $Y_0$ existe non seulement sur $\Dom_C$,
mais peut être prolongée dans un domaine, accessible ou non,
$\Dom_P\subset\Dom_0$, où $\Dom_P$ est un domaine borné, tel que
$\Dom_P\cap\Dom_C\neq
\emptyset$. On suppose aussi que toutes les fonctions $g$, $h$ et $P$
ont un prolongement analytique dans ce domaine. Il existe alors une
solution $Z(X,\eps')$ holomorphe sur $\Dom_P$ qui est un prolongement
de la solution $Y(X,\eps')$ du théorème~\ref{th:pres}. Cette solution
admet $\hat Y$ comme développement asymptotique uniforme sur
$\Dom_P$. Dans le domaine $\Dom_P\cup\Dom_C$, cette solution prolongée
$Y(X,\eps')$ tend uniformément vers $Y_0(X)$ quand $\eps'$ tend vers
$0$ et vérifie encore le corollaire \ref{cor:pabo}.

En particulier, si $0\in\Dom_P$, on peut prolonger les solutions
correspondantes de \eqref{eq:gen},
$y(x,\eps)=\frac{1}{\eps}Y(x/\eps',\eps)$, jusqu'à un voisinage de
$x=0$. Ce voisinage, cependant, n'est que de taille $\eps^{1/(p+1)}$
autour de $0$.

Par ailleurs si on est sûr que $y$ n'est pas exponentiellement
grande (en $\eps$) dans ce voisinage, elle est {\em a priori} de
taille $1/\eps$.
\end{Cor}

La démonstration de ce corollaire ne pose pas de problème. 

Soit $X_0$ un point de $\Dom_P\cap\Dom_C$. D'après le théorème, on a
une solution holomorphe $Y(X,\eps)$ sur $\Dom_C$, prolongeant
$y(x,\eps)$, et en particulier, $Y(X_0,\eps')$ est bien
défini. L'équation~\eqref{Yg}, régulièrement perturbée, a une solution
holomorphe sur $\Dom_P$, avec la condition initiale
$Z(X_0,\eps')=Y(X_0,\eps')$. Cette solution admet un développement
asymptotique $\sum Z_n{\eps'}^n$, qui est nécessairement égal à la
série formelle $\hat Y$, par construction: les $Y_n$ et $Z_n$ sont
définies par récurrence comme solutions d'équations différentielles
(les mêmes pour $Y_n$ et $Z_n$, quel que soit $n$) et ces fonctions
admettent la même condition initiale en $X_0$, donc $Y_N=Z_N$. Par
unicité, $Z$ est bien entendu un prolongement de $Y$.

\subsection{Connexion des développements asymptotiques}

On conserve dans ce paragraphe toutes les hypothèses des paragraphes
précédents.

En regardant alors le résultat du théorème \ref{th:loin,c} et du
corollaire \ref{cor:pabo}, on s'aperçoit qu'on n'a pas, pour
l'instant, de moyen de donner une bonne estimation de la solution
$y(x,\eps)$ trouvée pour certains $x$ au voisinage de $0$:
${|\eps'|}^{1-\mu}<|x|<|x_\rho|$ . Or on aimerait avoir une telle
estimation partout où la solution est définie.

Nous allons voir que la solution formelle extérieure reste un
développement asymptotique de la solution pour ces $x$, en adaptant
les démonstrations du paragraphe \ref{sec:loin}. Les hypothèses
\ref{H0} nous indiquent que si $|x|>{|\eps'|}^{\nu}$, $0<\nu<1$,
\begin{align*}
|g(x,\eps)| &< {\text B}_g {|\eps'|}^{-\nu}\\
|h(x,\eps)| &< {\text B}_h {|\eps'|}^{-\nu}
\intertext{puis, comme $f$ est bornée en $0$ et que $h(x,0)$ peut
être d'ordre $1/x$ en $x=0$,}
|y_0(x)|=\left|\f{h(x,0)}{x^p f(x)}\right| &<{\text B}_y{|\eps|}^{-\nu}
\end{align*}
On reprend les mêmes définitions pour les espaces et les normes
associées, et pour les opérateurs entre ces espaces, qu'à la
page~\pageref{Operateurs}, sauf que tout est cette fois défini pour
$x$ dans le domaine 
\[\Dom_\nu^\kappa=\Bigl\{x\in\Dom_0\ \Bigm/ \f{x}{\eps'}=X\in\Dom_C\, 
\text{ et }\ \kappa|x_\rho|>|x|>\f1\kappa{|\eps'|}^{\nu}\Bigr\}\]
Dans un tel domaine, les majorations pour les fonctions dérivées
seront les suivantes (quelles que soient les normes choisies sur des
espaces de fonctions holomorphes):
\begin{list}{{$\blacktriangleright$}}{}
\item si $|f(x)|< \text{A}$ pour tout $x\in\Dom^\kappa_\nu$, alors 
dans tout domaine $\Dom_\nu^\varkappa$ ($\varkappa<\kappa$), et pour
un certain $\text{B}>\text{A}$,
$\|f'(x)\|<\f{\text{B}}{|\eps'|^\nu}$. Il ne s'agit là que d'un cas
particulier de la formule de Cauchy.
\item si pour tout $\nu$, $\eps^m|h(x,\eps)|< \text{A}$, pour tout 
$x\in\Dom_\nu^\kappa$ et $\eps$ assez petit, alors on peut majorer la
différence ${h(x,\eps)-h(x,0)}$:
$\left\|\eps^m\f{h(x,\eps)-h(x,0)}\eps\right\| <
\f{\text{A}}{|\eps|^\nu}$.  En effet, loin de $x=0$, $h$ est
holomorphe en $x$ et admet un développement asymptotique en $\eps$
uniforme en $x$, ce qui fait que $\left(h(x,\eps)-h(x,0)\right)$
est borné par $\text{C} \eps$; au voisinage de $x=0$, où $h$ peut
avoir un p\^ole, on ne considère que la partie polaire de son
développement, qui nous emp\^eche d'avoir le m\^eme résultat:
$$h_{polaire}(x,\eps)=\f{h_{-1}}x +\cdots+f{h_{-n}}{x^n}+ 
\eps\left(\f{h_{-n-1}}{x^{n+1}}+\cdots+ \f{h_{N}}{x^{N}}\right) 
+\cdots$$ 
Ici $n$ est tel que $n\nu\leq m(p+1)$ pour tout $\nu<1$, donc
$n=m(p+1)$, et $N$ doit \^etre tel que $N\nu\leq (p+1)+m(p+1)$, donc
$N=n+p+1$; en continuant, on peut avoir dans $h_{polaire}$ des termes
en $\eps^2 x^{-2p-2-n}$, etc\ldots . En calculant alors la différence
$h(x,\eps)-h(x,0)$, on trouve alors le résultat, qui est vrai pour
toute fonction $h$ vérifiant la condition.
\end{list} 
\noindent{\bf Note}: Dans tout le reste de ce paragraphe, on omettra 
de noter les constantes $\kappa>1$ successives décroissantes qui
interviendront, et les constantes multiplicatives $\text{A}$,
$\text{B}$, etc.

On obtient comme majorations
\begin{align*}
\|g\|_{\HH}&=\sup_{x\in\Dom_\nu}\left|\f{g(x,\eps)}{x^p f(x)} 
\right|< {|\eps|}^{-\nu}\,\\
\|h\|_{\HH} &< {|\eps|}^{-\nu}\,\\
\|y_0\| &< {|\eps|}^{-\nu}\,
\intertext{et aussi}
\|y'_0\|_{\HH} &<\f{{|\eps|}^{-\nu}}{|\eps'|^\nu}{|\eps'|}^{-p\nu}
=|\eps|^{-2\nu} 
\intertext{alors que}
\|P\|_{\HH}&\ \text{ reste bornée quand $\eps$ tend vers $0$.}
\end{align*}
Cela va nous permettre de montrer le théorème suivant, qui prolonge le
théorème \ref{th:loin,c}:
\begin{Th}\label{th:-loin}
Si les hypothèses \ref{H0} et celles du théorème \ref{th:loin,c} sont
vérifiées, alors pour tout $\nu\in\, ]0,1[$, la solution $y(x,\eps)$
et la solution formelle $\hat y$ sont reliées dans le domaine
$\Dom_\nu^1$ (défini ci-dessus) par la relation:
$$\left|y(x,\eps)-\sum_{n=0}^N y_n(x)\eps^n\right|\leq C_N
|\eps|^{(N+1)}|\eps|^{-(N+2)\nu}$$
\end{Th}
Nous commencerons par démontrer le cas $N=0$, en reprenant
effectivement les mêmes opérateurs que ceux qu'on trouve en
page~\pageref{Operateurs}.

Alors, si $\|z_k\|< |\eps|^{-2\nu}$, ($\eps^2 z_k$ continue
heureusement à rester aussi petit que nécessaire), on obtient d'abord
l'inégalité suivante, déduite de \eqref{CH}:
\begin{align*}
\|\HHH(z_1)-\HHH(z_2)\|&\leq |\eps|\cdot\|z_1-z_2\|\times\\
                       &\phantom{\leq}\Bigl(\|g\|_{\HH}+
2(\|y_0\|+|\eps| \|z_k\|)\cdot\|P\|_{\HH} +
|\eps|(|\eps|\|z_k\|+\|y_0\|)^2 C\|P\|_{\HH} )\Bigr)\\
		       &\leq|\eps|
		       \left[{|\eps|}^{-\nu}+({|\eps|}^{-\nu}+|\eps|^{1-2\nu})
		       +|\eps|({|\eps|}^{-\nu}+|\eps|^{1-2\nu})^2\right]
\|z_1-z_2\|\\
\|\HHH(z_1)-\HHH(z_2)\|&\leq \left(|\eps|^{1-\nu}+|\eps|^{2-2\nu}
+|\eps|^{4-4\nu}\right)\|z_1-z_2\|\,,
\intertext{et on garde pour l'autre opérateur l'inégalité \eqref{CZ},
		       indépendante de $\eps$, vraie quelles que
		       soient les normes des $H_k$} 
\| {\ZZZ}(H_1)-{\ZZZ}(H_2) \| &\leq M. \|H_1-H_2\|_{\HH} 
\intertext{ce qui permet d'envisager que la combinaison des deux
		       opérateurs soit contractante, si les espaces
		       d'arrivée sont bien les bons. Or, si
		       $\|z\|<|\eps|^{-2\nu}$, d'après \eqref{eq:OH}}
\|\HHH(z)\|_{\HH} & \leq |\eps|\cdot \|z\|
\cdot \|g\|_{\HH} + \left\|y_0+\eps z^2\right\|\cdot 
\left\|P\right\|_{\HH}\\
		  & < |\eps|\cdot|\eps|^{-2\nu}\cdot{|\eps|}^{-\nu}+
{|\eps|}^{-\nu}\\
		  & < |\eps|^{1-3\nu}+{|\eps|}^{-\nu} \\
		  & < \sup\Bigl(|\eps|^{-\nu}, |\eps|^{1-3\nu}\Bigr)
\intertext{et dans l'autre sens, si 
$\|H\|_{\HH}<\sup\Bigl(|\eps|^{-\nu}, |\eps|^{1-3\nu}\Bigr)$, on
reprend l'inégalité \eqref{eq:OZ}} 
\|\ZZZ(H)\|&\leq  M \left\| y_0g-y'_0+\f{h(x,\eps)-h(x,0)}\eps
                  \right\|_{\HH}+ M \|H\|_{\HH}\\
           & < \Bigl(|\eps|^{-\nu}|\eps|^{-\nu}+
|\eps|^{-2\nu}|+|\eps|^{-\nu}|\eps|^{-\nu}\Bigr)+ 
\Bigl(|\eps|^{-\nu}+|\eps|^{1-3\nu}\Bigr)\\
           & < |\eps|^{-2\nu}
\end{align*}
Ce qui permet d'appliquer effectivement le théorème du point fixe et
de démontrer le théorème \ref{th:-loin} pour le cas $N=0$.

Pour le cas général, la méthode est exactement la même que pour le
théorème \ref{th:loin,c}; on obtient les mêmes équations pour les
$R_n$, la seule différence est que les fonctions ne sont plus bornées
indépendamment de $\eps$, mais sont bornées par une puissance de
$\eps$. Ainsi, en reprenant les notations de la démonstration,
\begin{align*}
\left\|\f{h_1(x,\eps)}{x^p f}\right\|&\leq \|y_0\| \cdot 
\|g\|_{\HH}+\left\|\f{y'_0}{x^p f}\right\| +
\left\|\sup_\eps \f{h(x,\eps)-h(x,0)}{\eps}\right\|_{\HH} + \|y\|^2 \cdot
\|P\|_{\HH}\\
	& <  {|\eps|}^{-\nu}{|\eps|}^{-\nu}+
|\eps|^{-2\nu} +{|\eps|}^{-\nu}
\sup\Bigl(|\eps|^{-\nu},|\eps|^{1-3\nu}\Bigr)+|\eps|^{-2\nu}\\ 
\left\|h_1\right\|_{\HH} & < |\eps|^{-2\nu}\,.
\end{align*}
Comme $R_1$ vérifie l'équation
$$\eps R'_1=(x^pf+\eps g) R_1+h_1\,$$
on voit que $\|R_1\|<|\eps|^{-2\nu}$ et $\|y_1\|<|\eps|^{-2\nu}$, puis
avec $\eps R_2=R_1-y_1$, 
$$\|R_2(x,\eps)\|<\|h_2\|_{\HH}=\|h_1\|_{\HH}+\|g\|_{\HH}\|y_1\|
<|\eps|^{-3\nu}\,.$$ 
Par récurrence on arrive alors sans difficulté au résultat complet du
théorème, puisqu'on aura toujours
$$h_n(x,\eps)=\f{h_{n-1}(x,\eps)-h_{n-1}(x,0)}{\eps}
-g(x,\eps)y_{n-1}(x)$$
et
$$\|R_n(x,\eps)\|\leq \|h_n\|_{\HH}\,,\text{ comme }\|y_n(x)\|
=\|h_n\|_{\HH}\,.$$

\bigskip

Ce théorème de prolongement de la solution extérieure va nous
permettre de finir la démonstration du paragraphe précédent, en
assurant la connexion à tous les niveaux $n$ entre la série formelle
intérieure $\sum Y_n(x/\eps'){\eps'}^n$ et la série formelle
extérieure $\sum y_n(x)\eps^n$. Nous appellerons cela le
\begin{Cor}\label{th:connexion}
Soit $N$ un entier positif et $\mu\in\, ]0,1[$. Si les $N$ premiers
termes de la série $(Y_0(X),Y_1(X),\ldots,Y_{N-1}(X))$ construite au
corollaire \ref{cor:pabo} existent et sont majorées jusqu'à l'ordre
$N-1$ comme le prévoit ce corollaire, alors on a encore
\[\left|Y\left(x_\rho{\eps'}^{\mu-1},\eps'\right)-\sum_{n=0}^{N}
Y_n(x_\rho{\eps'}^{\mu-1}){\eps'}^n \right|
\leq |{\eps'}|^{N+1} |{\eps'}|^{-N(1-\mu)}\,,\]
où $Y_N(X)=\lim_{\eps'\rightarrow 0} Z_N(X,\eps')$. 
\end{Cor}
Il est clair que ce théorème suffit à compléter la démonstration par
récurrence du corollaire \ref{cor:pabo}, puisque la condition initiale
pour $Z_{N+1}$ (qui correspond au premier membre divisé par
${\eps'}^{N+1}$) a alors au plus la m\^eme taille que celle prévue pour
la fonction elle-m\^eme.

\medskip
Soit $\nu<\mu$, avec pourtant $\nu$ peu différent de $\mu$. Il est
possible d'écrire pour tout $n\geq 0$, d'après le théorème
\ref{th:-loin} (on rappelle que ${\eps'}^{p+1}=\eps$).
\begin{align*}
Y\left(x_\rho{\eps'}^{\mu-1},\eps'\right)&=
\eps y\left(x_\rho{\eps'}^{\mu},\eps\right)\\
&=\eps y_0\left(x_\rho{\eps'}^{\mu}\right)+\eps^2
y_1\left(x_\rho{\eps'}^{\mu}\right) + \ldots 
+\eps^{n+1}y_n\left(x_\rho{\eps'}^{\mu}\right)
+o\left(\eps^{(n+1)-(n+2)\nu}\right)
\end{align*}
On décomposera $y$ jusqu'à un $n_0$ tel que $(n_0+1)-(n_0+2)\mu >
(N+1)-N(1-\mu)$. 

On conna\^it une majoration des $y_k$ sur $\Dom_\nu$:
$$\left|\eps^{k}y_k\left(x_\rho {\eps'}^\nu\right)\right| \leq
|\eps'|^{-(p+1)(k+1)\nu}\,.$$ 
Donc les $y_k(x)$, qui sont indépendants de $\eps$, ont en $x=0$ au plus
un p\^ole d'ordre $(k+1)(p+1)$.  Autrement dit, $y_k$ (pour tout $k\leq
n$) qui est méromorphe en $0$ comme toutes les fonctions de l'équation
\eqref{eq:gen}, s'écrit aussi, au voisinage de $0$
$$y_k(x)=y_{k,-(k+1)(p+1)} x^{-(k+1)(p+1)}+ y_{k,-(k+1)(p+1)+1}
x^{-(k+1)(p+1)+1}+\ldots +y_{k,0} +\ldots$$ 

L'hypothèse (de récurrence) principale du corollaire est
$$\left|Y\left(x_\rho{\eps}^{\mu-1},\eps'\right)-
\sum_{k=0}^{N-1}{\eps'}^k
Y_k\left(x_\rho{\eps'}^{\mu-1}\right)\right|< |\eps'|^N
|\eps'|^{(\mu-1)(N-1)}\,,$$ 
qui implique
$$\left|\eps\sum_{k=0}^{n_0} \eps^ky_k\left(x_\rho{\eps'}^{\mu}\right)
-\sum_{l=0}^{N-1}{\eps'}^l
Y_l\left(x_\rho{\eps'}^{\mu-1}\right)\right| < |\eps'|^N
|\eps'|^{(\mu-1)(N-1)}\,.$$ 

Le premier membre de l'inégalité ci-dessus correspond aussi à 
$${\eps'}^N Z_N\left(x_\rho{\eps'}^{\mu-1},\eps'\right)
+o\left(|\eps'|^N |\eps'|^{(\mu-1)(N-1)}\right)$$ par définition de
$Z_N$, qui est majorée dans $\Dom_{C\mu}$, d'après l'hypothèse de
récurrence.  On peut en déduire une définition de
$Y_N\left(x_\rho{\eps'}^{\mu-1}\right)$ à un $o\left(
|\eps'|^{(\mu-1)(N-1)}\right)$ près. En $X$, formellement, on a
d'abord
$$Y_N(X)=\lim_{\eps'\rightarrow 0}\f1{{\eps'}^N}\left(Y(X,\eps')-
\sum_{k\leq N-1} {\eps'}^k Y_k(X)\right)\,.$$
Puis,
\begin{align*}
{\eps'}^N Y_N(X) &= \text{Terme d'ordre } {\eps'}^N \text{ dans }
\sum_{k\geq 0} \eps^{k+1} y_k(X\eps')\\
       &= \text{Terme d'ordre } {\eps'}^N \text{ dans }
\sum_{k\geq 0} \eps^{k+1} \sum_{l\geq 0} y_{k,l-(k+1)(p+1)}
X^{l-(k+1)(p+1)} \f{{\eps'}^l}{{\eps}^{k+1}}
\intertext{ce qui mène à}
 Y_N(X) &\sim \sum_{k\geq 0} y_{k,N-(k+1)(p+1)}
X^{N-(k+1)(p+1)}\,:
\end{align*}
comme $Y_N$ est solution d'une équation différentielle avec une
irrégularité singulière en l'infini, on sait que ce développement
formel est effectivement un développement asymptotique pour $Y_N(X)$.

D'où
$$Y_N\left(x_\rho{\eps'}^{\mu-1}\right) 
= \sum_{k\leq n_0} y_{k,N-(k+1)(p+1)} 
\left(x_\rho{\eps'}^{\mu-1}\right)^{N-(k+1)(p+1)}\!+o\left(
|\eps'|^{(\mu-1)(n_0-(k+1)(p+1)}\right).$$
On peut prndre $n_0$ arbitrairement grand, donc dans la différence
$$\eps\sum_{k=0}^{n_0} \eps^ky_k\left(x_\rho{\eps'}^{\mu}\right)
-\sum_{l=0}^{N}{\eps'}^l Y_l\left(x_\rho{\eps'}^{\mu-1}\right)$$ 
il ne reste, à un $o\left(|\eps'|^{N+1} |\eps'|^{(\mu-1)N}\right)$
près, qu'une somme finie de sommes convergentes:
$$\sum_{k\leq n_0}\left( \sum_{n>N-(k+1)(p+1)} 
\eps^{k+1} y_{k,n} \left(x_\rho {\eps'}^\mu\right)^n\right)\,.$$
Or tous les termes de cette somme sont d'ordre, en $\eps'$, au moins
\begin{align*}
 & (k+1)(p+1)+\Bigl(N+1-(k+1)(p+1)\Bigr)\mu\\
=& (k+1)(p+1)(1-\mu) +(N+1) +(N+1)(\mu-1)\\
=& (N+1)-N(1-\mu)+(1-\mu)(kp+p+k)\,.
\end{align*}
Le dernier terme de cette expression étant strictement positif,
l'ordre est supérieur à $ (N+1)-N(1-\mu)$. Ce qui s'écrit
$$\left|Y\left(x_\rho{\eps}^{\mu-1},\eps'\right)-
\sum_{k=0}^{N}{\eps'}^k
Y_k\left(x_\rho{\eps'}^{\mu-1}\right)\right|< |\eps'|^{N+1}
|\eps'|^{(\mu-1)N}\,,$$ 
ce qu'on voulait démontrer.

\newpage\section{Existence de solutions bornées jusqu'en $x=0$}
\label{sec:alpha}

On reprend l'équation originale, en y introduisant explicitement un
multi-paramètre: 
$$\vec\alpha(\eps)=\Bigl(\alpha_0(\eps),\ldots,\alpha_i(\eps),\ldots,
\alpha_{p-1}(\eps)\Bigr)\,$$
qu'on supposera avoir une série asymptotique
$$\vec\alpha(\eps)\underset{\eps\rightarrow 0}{\sim}\sum_{n=0}^\infty
\vec a_n \eps^n\,.$$ 
On notera aussi
$$\alpha(x,\eps)=\sum_{i=0}^{p-1} x^i \alpha_i(\eps).$$ 

On considère donc un cas particulier de l'équation~\eqref{eq:gen},
qu'on écrira par exemple sous la forme
\begin{equation}
\eps y'=\bigl(x^p f(x)+\eps
g(x,\eps,\vec\alpha)\bigr) y +h(x,\eps,\vec\alpha) +\eps y^2 
P(x,\eps,\vec\alpha,\eps y)\label{eq:alpha}\tag{{\sc iii}}
\end{equation}
ou encore
\begin{equation}\tag{{\sc iii'}}\label{eq:alphaq}
\eps y'=\bigl(x^p f(x)+\eps
g(x,\eps)\bigr) y +h(x,\eps) +\eps y^2 P(x,\eps, \eps y) 
+Q(x,\eps,\vec\alpha ,y)\,,
\end{equation}
où on suppose que $Q$ comporte les termes suivants:
\begin{multline}\label{eq:decQ}
Q(x,\eps,\vec\alpha ,y)=\alpha(x,\eps) R_0(x) + 
\eps R_1(x,\eps,\vec\alpha)+\eps y R_2(x,\eps,\vec\alpha)\\ +
\eps y^2  R_3(x,\eps,\vec\alpha,\eps y)+x^p R_4(x,\eps,\vec\alpha)\,,
\end{multline}
avec $R_0(0)\neq 0$. Les deux premiers termes sont issus de $h$, qui
doit donc avoir une forme particulière; les termes $R_2$ et $R_3$
viennent de $g$ et $P$ respectivement, dans le cas le plus général.

Comme on va s'intéresser maintenant à des solutions $y(x,\eps)$ qui
soient bornées sur un voisinage entier de $0$, il paraît nécessaire
que cette fois toutes les fonctions qui interviennent dans l'équation
soient holomorphes en $x=0$.

\medskip
Nous allons dans cette partie tenter de montrer que, dans certaines
conditions, il existe de vraies solutions d'une équation
différentielle qui existent jusque dans un voisinage d'un point
tournant. Pour cela, nous allons montrer que les vraies solutions qui
existent sur des montagnes (cf.~\ref{sec:loin}) peuvent parfois \^etre
prolongées jusqu'au point tournant, en se servant des résultats du
paragraphe~\ref{sec:pres}. Ensuite, en jouant sur un
(multi-)paramètre, on fera en sorte que ces différents prolongements
aient la m\^eme valeur en $x=0$; cela signifiera que tous ces
prolongements correspondent à une seule solution de l'équation, qui
existera donc sur la réunion des domaines d'existence de ces solutions.

\subsection{Existence d'une solution formelle continue en $0$}

Dans le cas où l'équation~\eqref{eq:alphaq} dépend linéairement de
$\alpha(x,\eps)$, elle peut se réécrire
\begin{equation}
\eps y'=\bigl(x^p f(x)+\eps
g(x,\eps)\bigr) y +h(x,\eps) +\eps y^2 P(x,\eps,\eps y) +
\alpha(x,\eps) Q(x,\eps,y)\label{eq:alphas}
\end{equation}
Dans ce cas précis, on peut montrer qu'il existe une solution formelle
à l'équation qui est continue en $0$, à condition que $Q(0,0,0)$ soit
non nul.
Nous allons montrer que l'équation \eqref{eq:alphas} possède une
unique solution formelle $(\hat y,\ \hat a)$ qui
s'écrit comme série formelle
$$\hat y =\sum_{n\geq 0} \eps^n y_n(x)\ \text{ et }\
\hat a=\sum_{n\geq 0} \eps^n {a_n(x)},\qquad 
a_n(x)=\sum_{k=0}^{p-1} \alpha_{k,n} x^k\,, $$
où les $y_n$ eux-mêmes sont des fonctions de $x$ qui sont soit
holomorphes en $x$, soit au moins admettent, quand  $x$ tend vers $0$,
une série asymptotique bornée en $0$, notée $\hat y_n$. Tous les $a_n(x)$
sont, quant à eux, des polynômes de degré inférieur ou égal à $p-1$ en
$x$. Pour obtenir ce résultat, il suffit essentiellement de remplacer
dans \eqref{eq:alphas} $y$ et les $\alpha_k$ par leur série formelle, 
les fonctions holomorphes en $\eps$: $g$, $h$ et $P$ par leur série, 
puis d'identifier les coefficients de $\eps^n$.

On écrira la fonction $Q$ sous la forme 
$$ Q(x,\eps,y)= \sum_{n=0}^\infty Q_n(x,y) \eps^n\,,$$ et on
décomposera de la m\^eme manière la fonction $P(x,\eps,\eps y)$ par
rapport à $\eps$. Ce type de notation vaudra aussi pour $h(x,\eps)$ et
$g(x,\eps)$.
 
\bigskip
On commence pour $n=0$ par:
\begin{equation*}
0=x^p f(x) y_0(x) +h_0(x) + a_0(x) Q_0(x,y)
\end{equation*}
(on remarque que $Q_0(x,y)=Q(x,0,0)$).
D'où 
\begin{equation}
y_0(x)=-\frac{h_0(x) + a_0(x) Q(x,0,0)}{x^p f(x)}\label{eq:y0}
\end{equation}
On souhaite que $y_0$ soit continue en $0$. Il faut et il suffit pour
cela que la valuation de $h_0(x) + a_0(x) Q(x,0,0)$ soit supérieure ou
égale à $p-1$. Ce sera effectivement le cas si (et seulement si) on 
choisit pour $a_0(x)$ le polynôme obtenu par division des séries en
$x$ suivantes: $-h_0(x)$ par $Q(x,0,0)$. Cela nous donne bien une et
une seule valeur possible pour $a_0$ puis $y_0$.

Supposons à présent que l'on ait montré l'existence (et l'unicité) des
$n$ premiers termes des séries $\hat y$ et $\hat a$, donc regardé les
coefficients de $\eps^k$ pour $k\leq n-1$. Alors, on obtient l'égalité
suivante pour le coefficient de $\eps^n$:
\begin{multline*}
ny_{n-1}(x)=x^p f(x) y_n(x) +\! \sum_{k\leq n-1} g_k(x) y_{n-1-k}(x)
+h_n(x) +P_{n-1}(x,y)\\
+\sum_{k\leq n} a_k(x)Q_{n-k}(x,y)
\end{multline*}
On peut vérifier que $P_{n-1}(x,y)$, le coefficient de $\eps^{n-1}$
dans $y^2P(x,\eps,\eps y)$ ne dépend que de $x$ et de $y_0,
y_1,\ldots, y_{n-1}$. De même, $Q_k(x,y)$, d'après la forme de $Q$, ne
dépend que de $x$ et de $y_0, y_1,\ldots, y_{k-1}$.

Cela permet de trouver
$$y_n=\frac{\displaystyle ny_{n-1}-\sum_{k\leq n-1} g_k(x) y_{n-1-k}
+h_n(x) +P_{n-1}(x,y) +\sum_{k\leq n} a_k(x)Q_{n-k}(x,y)}{x^p f(x)}$$

Or l'expression
$$ny_{n-1}-\sum_{k\leq n-1} g_k(x) y_{n-1-k}
+h_n(x) +P_{n-1}(x,y) +\sum_{k\leq n-1} a_k(x)Q_{n-k}(x,y)$$
peut s'écrire (au moins formellement) comme une série en $x$ si on 
remplace tous les $y_k$ pour $k<n$ par leur série (formelle)
respective. On choisit alors le polynôme
$a_n(x)$ tel que les $p$ premiers termes du numérateur de la fraction
ci-dessus soient nuls, {i.e.} que $y_n$ reste bornée en $x=0$, à
nouveau en effectuant une division de séries. 

Il existe un unique polynôme $a_n$ de degré inférieur ou égal à
$p-1$ qui vérifie la propriété voulue, puis le $y_n$ calculé est bien
sûr unique lui aussi.

\bigskip

Nous pouvons ainsi construire des séries uniques $\hat a$ et $\hat y$,
solutions formelles de \eqref{eq:alphas}, telles que tous les
coefficients de $\hat y$ restent bornés quand $x$ tend vers $0$. En
revanche, rien ne prouve l'existence pour tout $\eps$ assez petit
d'une vraie solution $y(x,\eps), a(x,\eps)$ à l'équation 
différentielle qui soit continue et bornée quand $\eps\rightarrow 0$
dans un voisinage de $0$.

\bigskip
Pour le cas non linéaire en $\alpha$ de l'équation~\eqref{eq:alphaq},
la méthode de construction est analogue, mais plus
complexe. Cependant, dans le cas qui nous intéresse où $Q$ peut
s'écrire sous la forme \eqref{eq:decQ}, l'existence et l'unicité de la
solution formelle $(\hat y, \hat\alpha)$ reste assurée (pour une
démonstration complète, voir \cite{CRSS}).


\subsection{Existence de vraies solutions holomorphes jusqu'en 0}

\begin{Th}\label{th:y(0)borne}
On suppose que pour l'équation~\eqref{eq:alpha}, en plus des
hypothèses du théorème~\ref{th:pres}, on a que les deux fonctions $g$
et $h$ sont holomorphes en $0$. Alors il existe de vraies solutions
($y,\alpha$), où $\alpha$ s'écrit nécessairement
$\alpha(x,\eps)=a_0(x)+a_1(x)\eps +\gO(\eps^2)$ (les deux termes $a_0$
et $a_1$ étant imposés), telles que les $y(x,\eps)$ définies, bornées
sur une montagne vérifient en plus les propriétés suivantes:
\begin{enumerate}
\item $y(0,\eps)$ est définie et bornée quand $\eps$ tend vers $0$.
\item De plus $y(0,\eps)-y(0,0) =\gO(\eps)$ (ou $o(\eps)$).
\end{enumerate}
\end{Th}

D'après le théorème~\ref{th:pres}, il existe des solutions $y(x,\eps)$
de l'équation~\eqref{eq:alpha} dans des domaines tels que
$x>|X_l|\eps'$, ce qui ne permet pas encore de les définir en
$0$. Pour arriver jusqu'en $0$, vu le corollaire~\ref{th:+pres}, il
faut encore que $Y_0(X)$ puisse \^etre prolongée holomorphiquement
jusqu'en $X=0$.

Or précisément, dans le cas qui nous intéresse, où $g$ et $h$ sont
holomorphes en $0$, on a vu que $G_0(X)\equiv 0$ et $H_0(X)\equiv
0$. L'équation pour $Y_0$ est alors particulièrement simple:
$$\eqref{eq:Y_0g}\ \ \frac{dY_0}{dX}=\Bigl(X^p f(0)\Bigr) Y_0 +Y_0^2
P_0(X,Y_0)\,,$$ 
et la solution $Y_0$ qui s'impose est une solution évidente: la
fonction nulle, la seule solution de cette équation qui tend vers $0$
en l'infini dans un très large secteur; celle-ci existe bien s\^ur en
$X=0$, et on peut utiliser le corollaire~\ref{th:+pres} qui assure
alors l'existence de $y(x,\eps)$ jusqu'en $x=0$.

Mais ce résultat ne nous suffit pas, la solution $y(x,\eps)$ obtenue
n'étant pas bornée quand $\eps\rightarrow 0$. Pour aller plus loin, il
faut regarder exactement les premiers termes de la solution
correspondante $Y(X,\eps')$, que dans les limites du corollaire
\ref{cor:pabo} on mettra sous la forme
$$Y(X,\eps')\sim \sum_{n=0}^\infty Y_n(X) {\eps'}^n\,.$$

\medskip
Montrons que les $p+1$ premiers $Y_n$ sont tous identiquement
nuls.

On sait qu'on a une solution $y_0(x)$ pour $\eps=0$ qui est
continue en $x=0$. L'équation pour $y_0$ nous donne 
$$y_0(x)=-\frac{h\bigl(x,0,\vec a_0\bigr)}{f(x)x^p}$$ donc
$h\bigl(x,0,\vec a_0\bigr)$ est divisible par $x^p$. Si on pose
$\alpha=a_0+\eps\tilde\alpha$, la fonction $h$ peut se mettre sous la
forme $h(x,\eps,\vec\alpha)=x^p h_1(x)+\eps
h_2(x,\eps,\vec{\tilde\alpha})$. Dans toute la suite, on supposera
qu'on a déjà écrit $\vec\alpha$ sous cette forme et transformé $g$,
$h$ et $P$ en conséquence pour obtenir effectivement
l'équation~\eqref{eq:alpha}:
\begin{multline}\tag{\ref{eq:alpha}}
\eps y'=\bigl(x^p f(x)+\eps g(x,\eps,\eps\vec{\tilde\alpha})\bigr) y 
+x^p h_1(x)+\eps h_2(x,\eps,\vec{\tilde\alpha})\\ 
+\eps y^2 \Bigl(x^p P_1(x,\eps y)
+\eps P_2(x,\eps,\vec{\tilde\alpha}, \eps y)\Bigr)
\end{multline}

Le résultat pour l'équation en $Y$, \eqref{Yg}, est le suivant: 
\begin{multline}
\frac{dY}{dX}=\Bigl(X^p f(\eps' X)+\eps'
g(\eps'X,\eps,\vec{\tilde\alpha})\Bigr) Y\\ 
+\eps X^p h_1(\eps'X)+\eps'\eps h_2(\eps'X,\eps,\vec{\tilde\alpha})+ 
Y^2 P_\eps(X,\eps',\vec{\tilde\alpha},Y)\,.
\end{multline}
Donc, par récurrence finie, en faisant l'hypothèse de récurrence
$Y_{n-1}(X)=0$,  pour $Y_n$, $n$ allant de $0$ à $p$, les équations
sont équivalentes à 
$$\frac{dY_n}{dX}=\Bigl(X^p f(0)\Bigr) Y_n\,.$$ 
La fonction nulle est la solution de chacune de ces équations qui
vérifie la condition initiale voulue; c'est-à-dire, quand $\eps'$ tend
vers $0$,
$$Y_n\left(\f{x_\rho{\eps'}^\mu}{\eps'}\right)=\f1{{\eps'}^n}\left(\eps
y(x_\rho{\eps'}^\mu,\eps)-\sum_{k=0}^{n-1} {\eps'}^k
Y_k\left(\f{x_\rho{\eps'}^\mu}{\eps'}\right)\right)\ :$$
quand $\eps$ tend vers $0$, $\f\eps{{\eps'}^p}
y(x_\rho{\eps'}^\mu,\eps)$ tend vers $0$ et les autres termes à droite
sont nuls par hypothèse de récurrence, d'où la condition 
$$\lim_{X\rightarrow\infty} Y_n(X)=0,\ \ \forall n \leq p\,.$$
\begin{align*}
\intertext{Pour $n=p+1$, on obtient en revanche}
\frac{dY_{p+1}}{dX}&=\Bigl(X^p f(0)\Bigr) Y_{p+1} + X^p h_1(0)\\
\intertext{dont la solution qui nous intéresse vérifie }
\lim_{X\rightarrow\infty} Y_{p+1}(X)&=\lim_{\eps'\rightarrow 0}
y_0(x{\eps'}^\mu)+o\left({\eps'}^{1-\mu}\right)=y_0(0)
\intertext{c'est-à-dire}
Y_{p+1}(X)&=-\frac{h_1(0)}{f(0)}\,.\\
\end{align*}
On en conclut que $Y(X,\eps')$ est égal à
$$Y(X,\eps')=\eps \tilde Y(X,\eps') \sim \eps (Y_{p+1}(X)+\eps'
Y_{p+2}(X)+\ldots)\,, $$
d'où l'on déduit que
$$y(\eps'X,\eps)=\frac{Y(X,\eps')}{\eps}=\tilde Y(X,\eps')$$
reste bornée en $x=0$ (et dans un petit voisinage de taille $\eps'$
autour de ce point) quand $\eps\rightarrow 0$; ce qui démontre le
premier point du théorème.

\bigskip
Il reste à montrer que l'équivalent de $y(0,\eps)-y(0,0)$ est d'ordre
$\eps$ au moins. On continue dans ce but le calcul des premiers termes
de la série formelle de $\tilde Y$. 

Revenons d'abord à l'équation~\eqref{eq:alpha}. Comme il existe une
solution formelle, on sait calculer $y_1(x)$ avec l'équation suivante:
$$y'_0(x)=x^p f(x) y_1(x) +h_2(x,0,a_1(x))+ y_0^2(x) x^p P_1(x,0)
+y_0(x) g(x,0,0)\,,$$
donc pour que $y_1$ soit continue en $0$, il faut et il suffit que
\begin{equation}\label{prop-y1}
y'_0(x)-h_2(x,0,a_1)- y_0(x)g(x,0,0)\ \text{soit de valuation
$p$ au moins.}
\end{equation}
\medskip
Pour $\tilde Y$, on trouve comme équation différentielle non
singulièrement perturbée
\begin{multline*}
\frac{d\tilde Y}{dX}=\Bigl( X^p f(\eps' X)+
\eps' g(\eps'X,\eps,\eps\vec{\tilde\alpha}) \Bigr)\tilde Y +X^p h_1(\eps'X)
+\eps'h_2(\eps'X,\eps,\vec{\tilde\alpha})\\
+\eps \tilde Y^2 P_\eps(X,\eps',\vec{\tilde\alpha}, \tilde Y)
\end{multline*}
On ne cherche que les termes $\tilde Y_0$, $\tilde Y_1$,\ldots $\tilde
Y_p$, donc dans l'équation ci-dessus, on tronquera toutes les séries
pour ne garder que les termes de degré en $\eps'$ strictement
inférieur à $p+1$. En particulier, le terme $\eps \tilde Y^2
P_\eps(X,\eps',\vec{\tilde\alpha}, \tilde Y)$ peut \^etre oublié, et on
obtient pour ces fonctions $\tilde Y_j$ des équations différentielles
linéaires. 

Si on cherche à identifier les deux séries formelles, celle pour $y$
et celle pour $\tilde Y$, on s'aperçoit que les $\tilde Y_j$ (pour
$j\leq p$) ne devraient dépendre que de $y_0$, et \^etre de la forme
$\tilde Y_j = \lambda_j X^j$. Nous allons démontrer que les équations
pour les $\tilde Y_j$ admettent effectivement de telles solutions. À
partir de l'équation pour $\tilde Y$,
\begin{multline*}
\sum_{j=1}^p  \lambda_j j X^{j-1} {\eps'}^j \equiv \Bigl(X^p f(\eps'X)+
\eps'g(\eps'X,0,0)\Bigr)\sum_{j=1}^p  \lambda_j  X^{j} {\eps'}^j\\
+X^p h_1(\eps'X)+\eps'h_2(\eps'X,0,a_1)\qquad 
[{\text{mod }}{\eps'}^{p+1}]\,,
\end{multline*}
{en ne regardant dans l'égalité ci-dessus que les termes
dont le degré en $\eps'$ est inférieur ou égal à $p$.}
\begin{multline*}
\Longleftrightarrow
\eps'\Bigl(\sum_{j=1}^p  \lambda_j j X^{j-1} {\eps'}^j-
h_2(\eps'X,0,a_1) -g(\eps'X,0,0)\sum_{j=1}^p  \lambda_j  X^{j}
{\eps'}^j\Bigr)\\
\equiv X^p\Bigr(f(\eps'X)\sum_{j=1}^p  \lambda_j  X^{j}
{\eps'}^j +h_1(\eps'X)\Bigr)\qquad [{\text{mod }}{\eps'}^{p+1}] \,.
\end{multline*}
{À gauche, on ne dépasse pas le degré $p-1$ en $X$, alors
qu'à droite, on commence au degré $p$; la seule possibilité d'avoir
cette égalité est donc que les deux membres soient nuls.}
\begin{align*}
\intertext{À droite, cela donne (toujours en séries tronquées)}
\sum_{j=1}^p  \lambda_j X^{j}{\eps'}^j 
&\equiv -\frac{h_1(\eps'X)}{f(\eps'X)}=y_0(\eps' X)
\qquad [{\text{mod }}{\eps'}^{p+1}]\\
\intertext{puis à gauche, on retrouve}
&\eps'\Bigl(y'_0(\eps'X)-h_2(\eps'X,0,a_1)-g(\eps'X,0,0)y_0(\eps'X)\Bigr)\\
\intertext{qui d'après la propriété~\eqref{prop-y1} est de valuation $p$
en $\eps'X$, donc de valuation $p+1$ en $\eps'$.}
\end{align*}
Ces $\tilde Y_j$ doivent vérifier en plus une condition initiale, qui
est, quand $\eps'$ tend vers $0$,
\begin{align*}
\tilde Y_j\left(x_\rho{\eps'}^{\mu-1}\right)&\sim
\f1{{\eps'}^{p+1+j}}\left(\eps y_0\left(x_\rho{\eps'}^{\mu}\right)
+\gO(\eps^2) -\sum_{k=0}^{j-1} {\eps'}^k
\tilde Y_k\left(x_\rho{\eps'}^{\mu-1}\right)\right)\\
&\sim \f1{{\eps'}^{j}}
y_0\left(x_\rho{\eps'}^{\mu}\right)-\sum_{k=0}^{j-1} \lambda_k 
\left(x_\rho {\eps'}^\mu\right)^k
\intertext{soit, si $y_0(x)=\sum \lambda_k x^k$,}
\tilde Y_j\left(x_\rho{\eps'}^{\mu-1}\right) 
&\sim \lambda_j  \left(x_\rho{\eps'}^{\mu-1}\right)^j
\end{align*}
Les solutions $\tilde Y_j$ trouvées formellement sont donc bien aussi
effectivement les solutions des équations différentielles
correspondantes avec les conditions initiales voulues.

Pour $y(x,\eps)$, on en déduit alors
\begin{align}
y(\eps'X,\eps) &= \lambda_0 + \eps'\lambda_1 X +\cdots+{\eps'}^p
\lambda_p X^p +\gO(\eps)\notag\\
\label{prop:y}y(0,\eps)-y(0,0) &\sim \eps \tilde Y_{p+1}(0)\ \ \text{ voire
plus petit, si $\tilde Y_{p+1}(0)=0$.}
\end{align}

\subsection{Une estimation préliminaire}

On étudiera dans la suite de cette partie l'équation~\eqref{eq:alphaq}
dans laquelle on a déjà écrit $\vec\alpha=\vec
a_0+\eps\vec{\tilde\alpha}$.
\begin{equation*}
\eps y'=\bigl(x^p f(x)+\eps
g(x,\eps)\bigr) y +h(x,\eps) +\eps y^2 P(x,\eps, \eps y) 
+Q(x,\eps,\eps\vec{\tilde\alpha} , y)
\end{equation*}
la fonction $Q$ étant holomorphe en toutes ses variables; $Q$ s'écrit
comme dans l'égalité \eqref{eq:decQ}.  

On note encore $Q_{k}^0=R_0(0)\neq 0$.

On notera $\gamma_l$ un chemin descendant le relief à partir d'un
sommet de la $l-$ième montagne jusqu'en $x=0$ (on rappelle qu'il na\^it
$p+1$ montagnes autour du point $0$).

Nous aurons besoin plus loin de l'estimation suivante:
\begin{Lemmetech}\label{th:col}
Quand $\eps$ tend vers $0$,
$$\int_{\gamma_l}\!\! \exp\left({\frac{1}{\eps}\!\!\int_x^0\!\!\! t^p f(t) dt
+r(x)}\right) \frac{\partial Q}{\partial \tilde\alpha_k} (x,\eps,
\eps\vec{\tilde\alpha}, y)dx \sim \eps
C_k {\eps'}^{k+1} \exp\left({l(k+1)\frac{2i\pi}{p+1}}\right)$$
où $C_k$ ne dépend pas de $l$. 
\end{Lemmetech}
Il ne s'agit dans ce lemme que d'un variante de la méthode du point
col. La valeur de l'intégrale est essentiellement concentrée là où 
$\int_x^0\!\! t^p f(t) dt$ est nulle, c'est-à-dire en $x=0$. Comme
$\exp(r(x))$ et $f(x)$ sont holomorphes et non nulles
au voisinage de $x=0$, l'intégrale est peu différente de
$$\int_{\gamma_l}\!\! e^{r(x)}\exp\left(-\frac{1}{\eps} x^{p+1}
\frac{f(0)}{p+1}\right) \frac{\partial Q}{\partial \tilde\alpha_k} (x,\eps,
\eps\vec{\tilde\alpha}, y)dx$$
et on effectue le changement de variable 
$u=x\times\f{\sqrt[p+1]{f(0)/(p+1)}}{\eps'}$. 

On arrive à l'estimation
$$e^{r(0)}\eps'\sqrt[p+1]{(p+1)/f(0)}\int_{\gamma_l}\!\!e^{-u^{p+1}}
\frac{\partial Q}{\partial\tilde\alpha_k}
\left(u\eps'\sqrt[p+1]{(p+1)/f(0)},\eps,\eps\vec{\tilde\alpha},y\right)du$$
Or $$\f1\eps\frac{\partial Q}{\partial\tilde\alpha_k}
(x,\eps,\eps\vec{\tilde\alpha}, y)=x^k R_0
+\eps\left(\Dp{R_1}{\tilde\alpha_k} +y\Dp{R_2}{\tilde\alpha_k} +y^2
\Dp{R_3}{\tilde\alpha_k}\right)+x^p \Dp{R_4}{\tilde\alpha_k}\,, $$ 
et la dérivée partielle est donc équivalente à
$\Dp{Q}{\tilde\alpha_k}\sim\eps x^kQ_k^0$ quand $\eps$ puis $x$
tendent vers $0$. D'où la nouvelle approximation
$$e^{r(0)}\eps'\sqrt[p+1]{(p+1)/f(0)}\ \eps\, Q_k^0 
\int_{\gamma_l}\!\!e^{-u^{p+1}} \left(u\eps'\right)^k
\left(\sqrt[p+1]{\frac{p+1}{f(0)}}\right)^k du$$
aussi équivalente à
$$\eps\, C'_k{\eps'}^{k+1}\int_{\gamma_l}e^{-u^{p+1}} u^k du$$
Notons $I_k=\int_{\gamma_1}e^{-u^{p+1}} u^k du$, qui est une intégrale
bien définie. Pour passer du chemin
$\gamma_1$ à $\gamma_l$, au moins dans un voisinage de $0$, il suffit
d'effectuer une rotation d'angle $\frac{2(l-1)\pi}{p+1}$. Donc, en
posant $$u=v\exp\left(i \frac{2(l-1)\pi}{p+1}\right)\,$$
on obtient que
\begin{align*}
\int_{\gamma_l}e^{-u^{p+1}} u^k du&=\int_{\gamma_l}e^{-v^{p+1}} v^k
\exp\left(i \frac{2k(l-1)\pi}{p+1}\right) \exp\left(i
\frac{2(l-1)\pi}{p+1}\right)  dv\\
 &=\exp\left(i \frac{2(k+1)(l-1)\pi}{p+1}\right) I_k
\end{align*}
D'où le résultat final:
$$\!\int_{\gamma_l}\!\! \exp\left({\frac{1}{\eps}\!\!\int_x^0\!\!\! t^p
f(t) dt +r(x)}\right) \frac{\partial Q}{\partial \tilde\alpha_k}
(x,\eps, \eps\vec{\tilde\alpha},y)dx \sim 
C_k {\eps'}^{k+1}\eps\ \exp\!\left({(k+1)\frac{2li\pi}{p+1}}\right)\,,$$
la valeur de $C_k$ étant la constante non nulle
$$C_k=C'_k*I_k\exp\left({-(k+1)\frac{2i\pi}{p+1}}\right)\,.$$

\subsection{Existence de solutions $y(x,\eps)$ exceptionnelles qui
restent bornées dans un voisinage de taille fixe du point 0}

\begin{Th}\label{th:alphas}
Soit l'équation
\begin{equation}
\eps y'=\bigl(x^p f(x)+\eps g(x,\eps)\bigr) y +h(x,\eps) 
+\eps y^2 P(x,\eps,\eps y)
+ Q(x,\eps, \eps\vec{\tilde\alpha}, y)\tag{\ref{eq:alphaq}}
\end{equation}
où toutes les fonctions qui interviennent sont holomorphes dans un
voisinage $\Dom_0$ de $x=0$ (avec $f(0)\neq 0$), holomorphes en $\eps$
dans des secteurs centrés en $0$, $P$ étant aussi holomorphe en $\eps
y=0$. $Q$ se décompose comme en \eqref{eq:decQ}
\begin{multline*}
Q(x,\eps,\vec\alpha ,\eps y)=\alpha(x,\eps) R_0(x) + 
\eps R_1(x,\eps,\vec\alpha)+\eps y R_2(x,\eps,\vec\alpha) \\+
\eps y^2  R_3(x,\eps,\vec\alpha,\eps y)+x^p R_4(x,\eps,\vec\alpha)\,,
\end{multline*}
avec
$$\alpha(x,\eps)=
\sum_{k=0}^{p-1}\alpha_k(\eps)x^k=a_0(x)+\eps\tilde\alpha(x,\eps)\,.$$ 

Dans ce cas, il existe une unique solution formelle $(\hat y=\sum
y_n(x)\eps^n, \hat\alpha=\sum \vec a_n\eps^n)$ telle que $y_n(0)$ est
bien définie pour tout $n$.

On construit $p+1$ domaines $\Dom_l\subset\Dom_0$ contenant chacun une
des montagnes (pour le relief $\Re(\int^x \frac{f(t)t^p}{\eps} dt)$)
naissant au voisinage du point col $0$, chacun des domaines $\Dom_l$
étant  accessible à partir de son sommet, qu'on pourra éventuellement 
choisir infini (cf. fig~\ref{fig:alphas}).

Alors, si pour un certain $\delta$ fixé et tout $\eps$ assez petit
\begin{equation}\tag{H1}\label{C1}
\frac{|g(x,\eps)|+|h(x,\eps)|+|P(x,\eps,\eps \delta)|}{x^p
f(x)} \hbox{ est bornée sur } 
\bigcup_{l=1}^{p+1} \Dom_l,
\end{equation}
et
\begin{equation}\tag{H2}\label{C2}
\frac{Q(x,\eps,\delta,\delta)}{x^p f(x)}
\hbox{ est bornée sur }  
\bigcup_{l=1}^{p+1} \Dom_l,
\end{equation}
il existe, pour tout $\eps>0$ assez petit, des $\alpha_k(\eps)$ et une
fonction $y(x,\eps)$ holomorphe bornée
sur $\bigcup \Dom_l$ et sur un voisinage complet de $x=0$, 
$\mathcal{B}(0,\rho)$, tels que l'équation \eqref{eq:alphaq} soit
vérifiée. 

Cette solution $(y,\vec\alpha)$ correspond à la solution formelle,
c'est-à-dire qu'elle admet les séries formelles solutions comme
développements asymptotiques:
$$y(x,\eps)\sim\sum_{n=0}^\infty y_n(x)\eps^n,\ \ \ \ 
\vec\alpha(\eps)\sim\sum_{n=0}^\infty \vec a_n\eps^n\,.$$
\end{Th}

\begin{figure}[ht!]
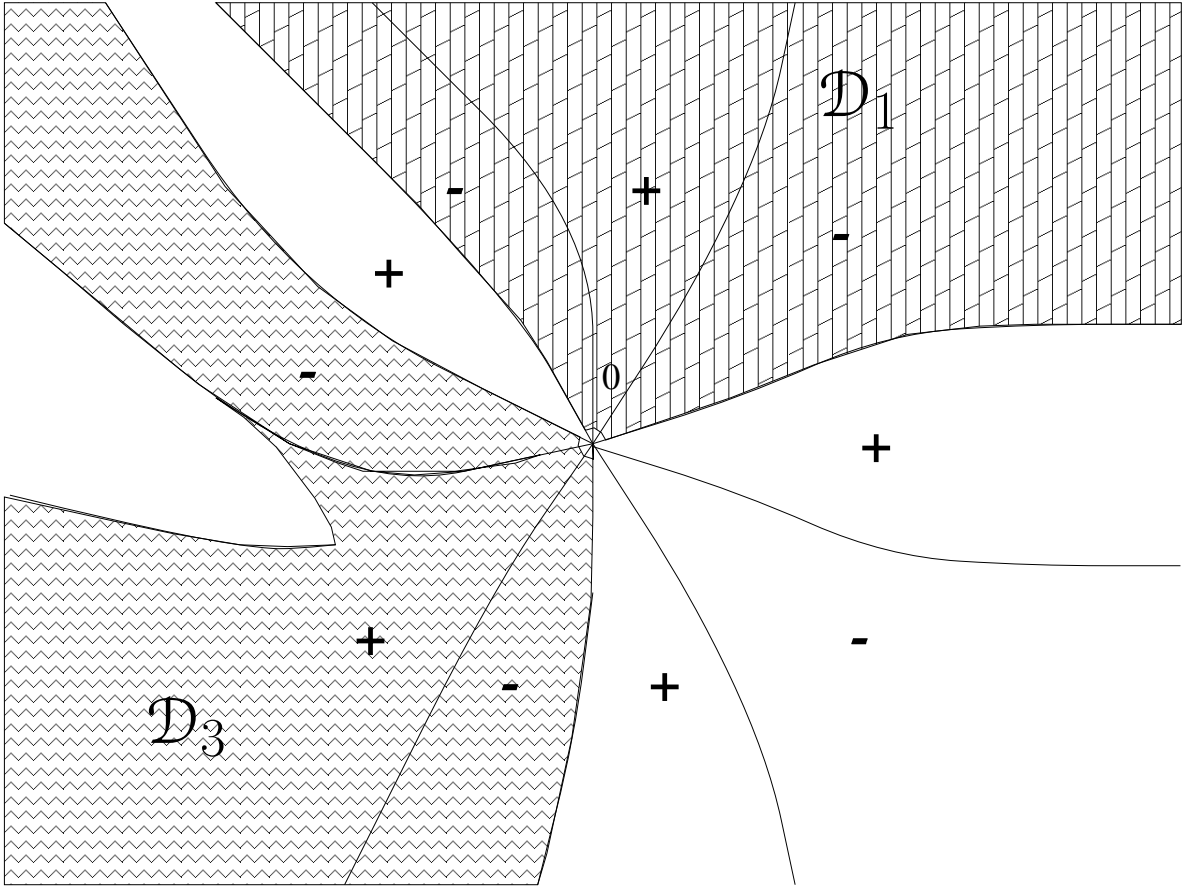

\centerline{\input Domalpha.pstex_t}
\caption{Exemple de voisinage de $0$, avec représentation de $\Dom_1$
et $\Dom_3$ ($p=4$)}\label{fig:alphas}
\end{figure}


Pour tout $\tilde\alpha$ dans un voisinage de taille $\eps$ de
$a_1(x)$ (trouvé formellement, comme $a_0$ précédemment), on sait que
sur chaque domaine $\Dom_l$, il existe une solution
$y^{[l]}(x,\eps)$. Si les hypothèses du théorème~\ref{th:alphas} sont
vérifiées, on a montré (au théorème~\ref{th:+pres}) qu'on peut
prolonger ces solutions jusqu'en $x=0$, et qu'elles restent bornées en
ce point quand $\eps\rightarrow 0$.  Pour démontrer le théorème, nous
allons montrer que pour tout $\eps$ assez petit, il existe un
$\tilde\alpha(x,\eps)$ tel que les fonctions $y^{[l]}(x,\eps)$
définies chacune sur leur montagne $\Dom_l$ soient toutes égales en
$x=0$: $y^{[1]}(0,\eps)=y^{[2]}(0,\eps)=\cdots=y^{[p+1]}(0,\eps)$.

\medskip
Regardons comment $y$ varie en fonction de
$\vec{\tilde\alpha}$. L'équation~\eqref{eq:alphaq} implique, si on la dérive par
rapport à $\alpha_i$:  
\begin{multline*}
\eps \left(\frac{\partial y}{\partial\tilde\alpha_k}\right)'=
\Dp{Q}{\tilde\alpha_k}(x,\eps,\eps\vec{\tilde\alpha},y)
+\frac{\partial y}{\partial \tilde\alpha_k}\times\\
\left[x^p f(x)+\eps
\left(g(x,\eps)+2yP(x,\eps,\eps y)+\eps y^2 \frac{\partial P}{\partial
y}(x,\eps,\eps y)+\f1\eps \frac{\partial Q}{\partial
y}(x,\eps,\eps\vec{\tilde\alpha},y)\right)\right]\,.
\end{multline*}
On voit qu'on a une équation différentielle linéaire en
$\frac{\partial y}{\partial\tilde\alpha_k}$; les solutions peuvent donc
s'écrire:
$$\frac{\partial y^{[l]}}{\partial\tilde\alpha_k}=\frac{1}{\eps}
\int_{\gamma_l} e^{r(x)}\exp\left(\frac{1}{\eps}\int^0_x t^p f(t)dt\right) 
\frac{\partial Q}{\partial\tilde\alpha_k}(x,\eps,\eps\vec{\tilde\alpha},y) 
dx\,,$$
avec
$$r(x)=\int^0_x \left( g(t,\eps)+2yP(t,\eps,\eps y)+\eps y^2
\frac{\partial P}{\partial y}(t,\eps,\eps y)+\f1\eps
\frac{\partial Q}{\partial y}(t,\eps,\eps\vec{\tilde\alpha},y)\right)dt\,,$$
une fonction bornée quand $\eps$ tend vers $0$, puisque $\Dp{Q}{y}$
est un $\gO(\eps^2)$.
On peut appliquer le lemme~\ref{th:col}, qui nous donne alors:
\begin{equation}
\frac{\partial
y^{[l]}}{\partial\tilde\alpha_k}\underset{\eps\rightarrow 0}\sim
C_k {\eps'}^{k+1}\,\eps\, \exp\left({l(k+1)\f{2i\pi}{p+1}}\right)\,. 
\end{equation}

En reprenant le résultat complet du théorème~\ref{th:y(0)borne}
(cf.~\eqref{prop:y}), on peut écrire une formule de Taylor pour chacun
des $y^{[l]}$ au voisinage de la fonction $y_0$ en $x=0$:
\begin{align*}
y^{[l]}(0,\eps,\vec{\tilde\alpha})-y_0(0)&=\eps \tilde Y_{p+1}(0)
+\eps\sum_{k=0}^{p-1} \tilde\alpha_k C_k {\eps'}^{k+1}
\exp\left({l(k+1)\frac{2i\pi}{p+1}}\right)+o(\ldots)
\intertext{qui donne, après division par $\eps$,}
\f1\eps \left(y^{[l]}(0,\eps,\vec{\tilde\alpha})-y_0(0)\right)
&=\tilde Y_{p+1}(0) +\sum_{k=0}^{p-1} \check\alpha_k C_k
\exp\left({l(k+1)\frac{2i\pi}{p+1}}\right)+o(\ldots)\,,
\end{align*}
où $\check\alpha_k=\tilde\alpha_k {\eps'}^{k+1}$. On rappelle que ni
$\tilde Y_{p+1}$ ni $y_0$ ne dépendent de $l$, $\eps$ ou $\alpha$.

\medskip
Nous allons regarder le vecteur 
$$\gY=\begin{pmatrix}
\frac1\eps\left(y^{[1]}(0,\eps,\vec{\tilde\alpha})-y_0(0)\right)-\eta\\
\frac1\eps\left(y^{[2]}(0,\eps,\vec{\tilde\alpha})-y_0(0)\right)-\eta\\
\vdots\\
\frac1\eps\left(y^{[p+1]}(0,\eps,\vec{\tilde\alpha})-y_0(0)\right)-\eta
\end{pmatrix}$$ 
de taille $p+1$, dépendant des $p+2$ variables $\eps$, $w_1=\eta$,
$w_2=\check\alpha_0$, $\ldots$, $w_{p+1}=\check\alpha_{p-1}$.  D'après
ce qui précède, la matrice jacobienne de $\cal Y$ par rapport aux
$p+1$ dernières variables est 
\begin{gather*}
\Dp{\gY}{w}_{\eps=0,w=(\tilde Y_{p+1}(0),0,\ldots,0)}= 
 \begin{bmatrix} 
  -1 &
C_0\,e^{\left(1\cdot 1\cdot \f{2i\pi}{p+1}\right)} & \ldots &
C_{p-1}\,e^{\left(1\cdot p\cdot \f{2i\pi}{p+1}\right)}\\
  -1 &
C_0\,e^{\left(2\cdot 1\cdot \f{2i\pi}{p+1}\right)} & \ldots &
C_{p-1}{\,}e^{\left(2\cdot p\cdot \f{2i\pi}{p+1}\right)}\\
\vdots & \vdots & \vdots & \vdots\\ 
  -1 &
C_0\,e^{\left((p+1)\cdot 1\cdot \f{2i\pi}{p+1}\right)} & \ldots &
C_{p-1}{\,}e^{\left((p+1)\cdot p\cdot \f{2i\pi}{p+1}\right)}
 \end{bmatrix}
\end{gather*}
Le déterminant correspondant est non nul. En effet, en factorisant le
jacobien obtenu par colonnes, on obtient 
\[
{\mathit Det}\left(\Dp{\gY}{w}\right)=
-\prod_{k=0}^{p-1} C_{k}
\begin{vmatrix}
1 & \left(\exp\f{2i\pi}{p+1}\right) &  \ldots & 
\left(\exp\f{2i\pi}{p+1}\right)^{p}\\
1 & \left(\exp\f{4i\pi}{p+1}\right) &  \ldots & 
\left(\exp\f{4i\pi}{p+1}\right)^{p}\\
\vdots & \vdots & &\vdots\\
1 & \left(\exp\f{2(p+1)i\pi}{p+1}\right) &  \ldots & 
\left(\exp\f{2i(p+1)\pi}{p+1}\right)^{p}\\
\end{vmatrix}
\]
Le déterminant qui reste est un déterminant de Van der Monde, qui est
bien non nul.

On en déduit que $\gY$ est localement inversible autour de $w=0$. 

On regarde la fonction $\gY(\eps,w)$; on sait que, avec $w_0=(\tilde
Y_{p+1}(0),0,\ldots,0)$, 
\begin{align*}
\gY(0,w_0)&=0 
\intertext{et}
\Dp{\gY}{w}_{(\eps=0,w=w_0)}\ \text{ est inversible}
\end{align*}
D'après le théorème des fonctions implicites, il existe donc une
unique fonction $w(\eps)$ telle que $\gY\bigl(\eps,w(\eps)\bigr)=0$,
pour tout $\eps$ dans un voisinage suffisamment petit de $0$.  Les
fonctions $\check\alpha(\eps)$ que l'on en déduit sont bornées en
$\eps$; donc les fonctions
$\eps\vec{\tilde\alpha}=\check\alpha{\eps'}^{p-k+1}$ tendent vers $0$
avec $\eps$. Cette propriété nous suffit pour pouvoir effectivement
appliquer le théorème des fonctions implicites: les $\tilde\alpha_k$
sont toujours accompagnés d'un $\eps$ en facteur dans toute ce
paragraphe. Les fonctions $\eps\vec{\tilde\alpha}(\eps)$ sont
holomorphes en $\eps$ dans des secteurs ouverts en $0$, puisque le
passage de $\eps$ à $\eps'=\eps^{\f{1}{p+1}}$ n'est pas bijectif dans
un visinage du point $0$.


Avec $\vec\alpha$, on trouve en m\^eme temps la fonction $y(x,\eps)$
correspondante à partir de $y(0,\eps)=y_0(0)+\eps \eta(\eps)$; dans
les paragraphes précédents, on a montré qu'il s'agissait d'une
fonction bornée, en $x=0$ et dans tous les domaines $\Dom_l$ autour de
$0$. On vient ainsi de finir de démontrer l'existence de la vraie
solution $(y,\alpha)$ du théorème, bornée quand $\eps$ tend vers $0$,
dans l'union des domaines $\Dom_l$.

\medskip
Il reste à justifier que la solution formelle $(\hat y,\hat\alpha)$
est le développement asymptotique de $(y,\vec\alpha)$. La
démonstration se fait de manière analogue à celle du
théorème~\ref{th:loin,c}. On sait que la solution $(y,\vec
a_0+\eps\vec{\tilde\alpha})$ existe, avec $\eps\vec{\tilde\alpha}$ bornée
en $\eps$, et on écrit
\begin{align*}
y(x,\eps)&=y_0(x)+\eps z(x,\eps)\\
a_0(x)+\eps\tilde\alpha(x,\eps)&=a_0(x)+\eps a_1(x)+\eps^2
\breve\alpha(x,\eps)\,.
\end{align*}
Si on introduit ces notations dans l'équation
différentielle~\eqref{eq:alphaq}, cela donne une équation
différentielle linéaire en $z$ avec un paramètre
$\breve\alpha(x,\eps)$ qui intervient lui aussi de manière linéaire:
\begin{multline*}
\eps z' =\left(x^p f(x)+\eps g(x,\eps)\right)z+
\eps\breve\alpha(x,\eps) R_0(x)+x^p 
\f{R_4(x,\eps,\eps\vec{\tilde\alpha})-R_4(x,0,0)}\eps +\\
\left(\f{h(x,\eps)-h(x,0)}\eps-y'_0+g(x,\eps) y_0 +y^2 P(x,\eps,\eps
y)+a_1(x,\eps) R_0+ R_1+ y R_2+ y^2 R_3\right)\!
\end{multline*}
(on rappelle que $R_0$ ne dépend que de $x$; $R_1$, $R_2$ et $R_4$ de
$x$, $\eps$ et $\eps\vec{\tilde\alpha}$; et que $R_3$ dépend en plus
de $\eps y$).

On vérifie très facilement que les hypothèses du théorème sont encore
vérifiées pour cette équation:
\begin{enumerate}
\item $f$ et $g$ ne changent pas.
\item la nouvelle fonction $\breve P$ est nulle.
\item la nouvelle fonction $\breve Q(x,\eps,\eps\vec{\breve\alpha},\eps y)$, 
quant à elle, se réduit à la somme des deux termes
$\eps\breve\alpha(x,\eps)R_0(x)+x^p
\f{R_4(x,\eps,\eps\vec{\tilde\alpha})-R_4(x,0,0)}\eps$. Le premier terme 
ne pose aucune difficulté. Et le second peut s'écrire aussi 
$x^p \breve R_{4}(x,\eps,\eps\vec{\breve\alpha})$.
\item $\breve h(x,\eps)$ correspond à la deuxième ligne dans l'équation 
ci-dessus. Cette fonction qui peut s'exprimer en fonction de $x$ et
$\eps$ uniquement (puisque $y$ et $\eps\tilde\alpha$ sont bien connus)
vérifie la condition de majoration voulue.
\end{enumerate}
Avec ce qui vient d'\^etre démontré, on sait donc qu'il existe
$z(x,\eps)$ et $\eps\breve\alpha(x,\eps)$ bornées en $\eps$ telles que
$z$ existe et est holomorphe en $x=0$, et $\eps\breve\alpha(\eps)$ est
holomorphe en $\eps$.

Une simple récurrence suffit ensuite pour clore la démonstration du
théorème~\ref{th:alphas}.

\newpage\chapter{L'équation du Brusselator}\label{chap:Brusselator}

Le système différentiel du Brusselator apparaît pour l'étude de
l'évolution cinétique de certaines réactions chimiques
auto-catalytiques. Il peut se mettre sous la forme suivante:
\[ \left\{ \begin{array}{rl}
   \eps \dot \chi &=a(1+\chi)^2 y +a\chi^2+(a-1)\chi\\
      \dot y   &=-\chi(1+\chi)-y(1+\chi)^2
           \end{array}
\right. \]
Qualitativement, il y a des solutions «canard» de ce système pour
des valeurs du paramètre $a$ voisines de $a_0=1$. Dans ce cas, la courbe
lente a pour équation $y=\frac{-\chi^2}{(1+\chi)^2}$, et les trajectoires
des solutions peuvent être amenées à suivre cette courbe lente y
compris sur sa partie instable (cf. figure \ref{fig:Brusselc}).
On s'intéressera en fait aux grands canards, c'est-à-dire aux trajectoires
qui suivent la partie instable de la courbe jusqu'en $x$ infini.

\begin{figure}[ht!]
\center{\includegraphics[scale=0.85]{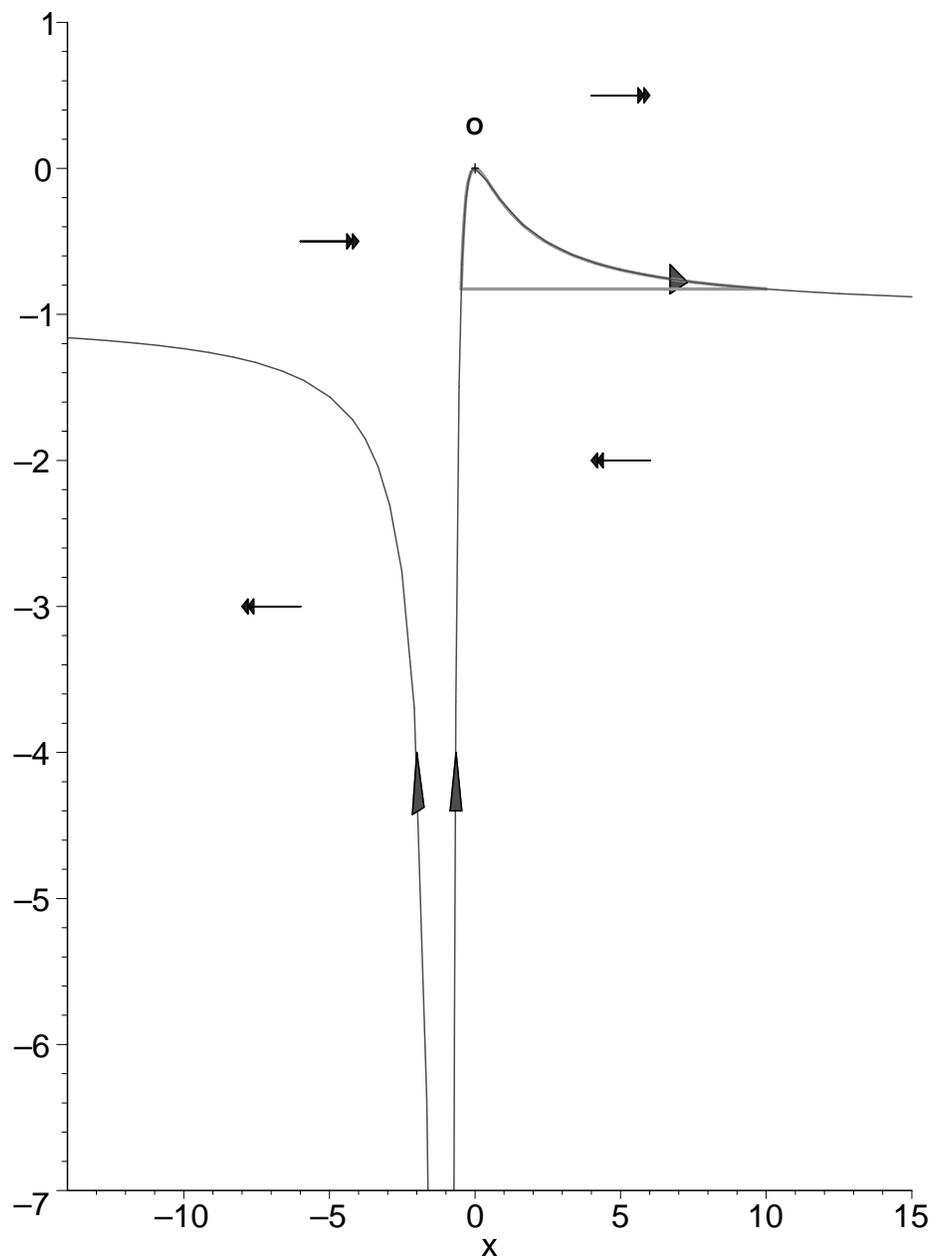}}
\caption{Courbe canard du Brusselator}\label{fig:Brusselc}
\end{figure}

\section{Existence de solutions canards}

Nous allons commencer par transformer le système pour obtenir une
équation sous une forme \og normale\fg, correspondant à celle étudiée
dans la partie précédente (en particulier, nous ramènerons le point
col pour le relief situé à l'infini en un point fini). Puis nous
verrons quel est le domaine d'existence des solutions surstables à
l'équation obtenue.

Pour étudier ces solutions, nous allons prendre la nouvelle variable
$z=\dot \chi$, ce qui ramène la première équation du système à 
\begin{align*}
a (1+\chi)^2 y&=-a\chi^2-(a-1)\chi+\eps z\\
\intertext{et la seconde à:}
z\frac{dy}{d\chi}=\frac{dy}{dt}&=-(\chi+1)\chi+\chi^2+\frac{a-1}{a}\chi-\frac{\eps}{a}z
=\frac{-\chi-\eps z}{a}\,.\\
\intertext{Or, la première équation ci-dessus nous donne aussi}
\frac{dy}{d\chi}&=\frac{-2\chi}{(1+\chi)^3}-\frac{a-1}{a}\frac{1-\chi}{(1+\chi)^3}+
\frac{\eps}{a}\frac{2z+(1+\chi)dz/d\chi}{(1+\chi)^3}\,,\\
\intertext{soit}
a(1+\chi)^2\frac{dy}{d\chi}&=-a+\frac{1-\chi}{1+\chi}+
\eps\Bigl(\frac{2z}{1+\chi}+\frac{dz}{d\chi}\Bigr)\,.
\end{align*}
Le système du Brusselator peut donc se ramener à l'équation suivante, 
que nous appellerons équation du Brusselator.
$$\eps z \frac{dz}{d\chi}=\frac{2\chi}{1+\chi}
\left(z-\frac{(1+\chi)^3}{2}\right)
+(a-1)z+z\left(z-\frac{(1+\chi)^3}{2}\right)\frac{2\eps}{1+\chi}$$

On voit déjà sur cette équation que la fonction 
$\Phi_0(\chi)=\frac{(1+\chi)^3}{2}$ est particulière pour cette 
équation: pour le paramètre $a=1$ elle est solution quand $\eps=0$, et 
c'est alors la seule solution (pour  $\eps=0$ quand on fait varier $a$) 
qui soit continue en $\chi=0$.

\medskip
On effectue le changement de variable suivant, qui ramène l'infini en $-1$
et laisse $0$ en $0$.
$$\chi=\frac{-x}{1+x}\ \ , \text{ puis } 
\frac{dx}{d\chi}=-(1+x)^2$$ 
Si on compare alors avec l'équation de Van der Pol, $0$ correspondra
ici à $1$ (le point col qui n'est, exceptionnellement, pas un point
tournant), $-1$ à $-1$ (l'autre point col du relief).

D'où l'équation
\begin{equation}
\label{Brus0}\eps z(x)\frac{dz(x)}{dx}=\frac{2x}{(1+x)^2} 
\Bigl(z-\frac{1}{2(1+x)^3}\Bigr)-\frac{a-1}{(1+x)^2}z 
-z\Bigl(z-\frac{1}{2(1+x)^3}\Bigr)\frac{2\eps}{1+x}
\end{equation}

Montrons que cette équation est bien du genre étudié lors de la
généralisation (voir partie précédente). 
On pose pour cela $z=\Phi_0(x)(1+\eps y)$, où
$\Phi_0(x)$ correspond à la courbe lente, ici $\frac{1}{2(1+x)^3}$, et
on remplace dans l'équation \eqref{Brus0} après en avoir divisé les
deux membres par $z$ 
$$\eps\Phi'_0(1+\eps y)+\eps^2 y'\Phi_0 =
\frac{2x\eps y}{(1+x)^2(1+\eps y)}
-\frac{a-1}{(1+x)^2}-\Phi_0\eps y\frac{2\eps}{1+x}$$
Donc
$$\eps y'=y\left[4x(1+x)+\frac{\eps}{1+x}\right]+\frac{3}{1+x}
-2(x+1)\frac{a-1}{\eps}-\eps y^2 4x(1+x)\frac{1}{1+\eps y}$$
On voit que autour de $0$, on a:
\begin{list}{$\bullet$}{}
\item $p=1$ et $f(x)=4(1+x)$
\item $g(x,\eps)=\frac{1}{1+x}$
\item $h(x,\eps)=\frac{3}{1+x}$
\item $P(x,\eps,\eps y)=-\frac{4x(1+x)}{1+\eps y}$
\item $Q(x,\eps,y)=-2(x+1)$ avec $\alpha=\frac{a-1}{\eps}$.
\item On prend alors pour $\Dom_0$ l'ensemble $\CC$ privé du
demi-axe réel $]-\infty,-1]$.
\item Toutes les hypothèses étant respectées, vu le relief
représenté figure \ref{fig:Brusselr}, le théorème~\ref{th:alphas}
s'applique: pour une valeur de $a$ bien choisie, soit $a^+(\eps)$, il
existe une solution holomorphe $z^+(x)$ qui sera continue en $0$ et
qui a pour domaine d'existence les montagnes Nord et Est et les
vallées adjacentes; voir fig.~\ref{fig:Brusselr}. Symétriquement, en
prenant les conjugués, on trouve une solution $z^-=\overline{z^+}$
correspondant au paramètre $a=a^-(\eps)=\overline{a^+}$.
\end{list}
Bien évidemment, ces deux solutions $(a,z)$ correspondent toutes les
deux aux séries formelles uniques
$$\hat a(\eps)=\sum a_n\eps^n \ \text{ et }\ \hat z(\eps,x)=\sum
z_n(x)\eps^n$$ construites de telle sorte que tous les termes $z_n$
sont continus en $x=0$.
\bigskip

\begin{figure}[ht!]
\caption{Relief correspondant à $\int_0^x 2t(1+t)dt$: domaine accessible 
à partir des sommets Est ou Nord, pour $a=a^+$; les courbes passant par 
$0$ sont les courbes de niveau $R(x)=0$, celles passant par $-1$ ont pour
équation $R(x)=1/3$.}
\center{\includegraphics[scale=0.7]{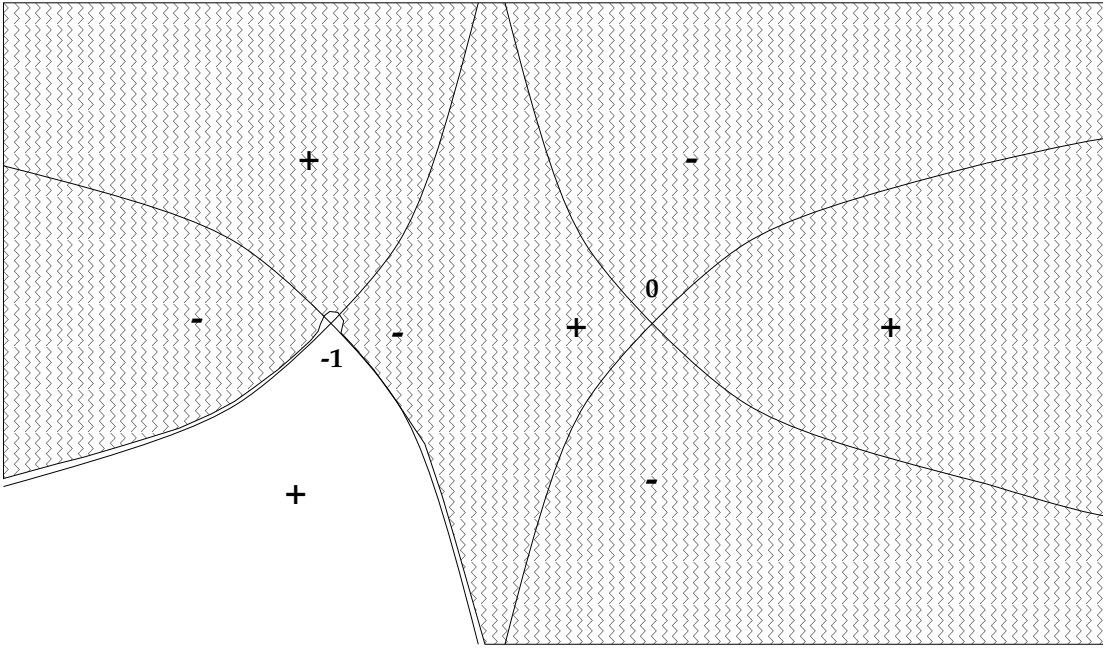}}
\label{fig:Brusselr}
\end{figure}

Si on regarde ensuite ce qui se passe autour de $-1$, en changeant de
variable: $\tilde x=x-1$, on trouve de manière analogue
\begin{list}{$\bullet$}{}
\item $p=1$, $f(\tilde x)=4(1-\tilde x)$ et $\eps'=\sqrt{\eps}$
\item $g(\tilde x,\eps)=
\frac{1}{\tilde x}$ d'où aussi $G_0(X)=\frac{1}{X}$
\item $h(\tilde x,\eps)=
\frac{3}{\tilde x}-2\tilde x\frac{a-1}{\eps}$ d'où aussi
$H_0(X)=\frac{3}{X}$ 
\item $P(\tilde x,\eps,\eps y)=
\frac{4\tilde x(1+\tilde x)}{1+\eps y}$ d'où aussi
$P_0(X,Y)=\frac{4X}{1+Y}$ 
\item On prend alors pour $\Dom$ la montagne Est (à partir de $-1$)
et éventuellement une partie des vallées adjacentes.
\item Les hypothèses des théorèmes~\ref{th:pres}, \ref{cor:pabo}, 
\ref{th:-loin} sont respectées: 
on peut prolonger le domaine jusqu'en $-1+X_l\sqrt\eps\in \RR$; et la
solution canard admet une sorte de développement asymptotique en
$\sqrt\eps$ pour des $x$ tels que $|1+x|<|\sqrt\eps|^{1+\nu}$
($\nu>0$), obtenu à partir de la solution formelle de l'équation
différentielle non singulièrement perturbée; tandis que le
développement asymptotique $\hat z(\eps,x)=\sum z_n(x)\eps^n$ reste
valable, dans le sens du théorème \ref{th:-loin}, pour les $x$ (sur la
montagne Est) tels que $|1+x|>|\sqrt\eps|^{1+\mu}$, $0<\nu<\mu<1$.
\end{list}

\section{Estimation de $a^+-a^-$}

L'existence de solutions là où on le souhaite est donc
acquise. Maintenant, pour des raisons de simplicités de calcul, nous
allons repartir de l'équation \eqref{Brus0}; la progression suivra alors 
exactement le même shéma que pour l'équation de Van der Pol. C'est à
dire que nous allons établir une équation reliant la différence
$a^+-a^-$ à la différence entre les fonctions $z^+-z^-$, étudiée au
voisinage du point singulier $-1$; pour pouvoir ensuite estimer chacun
des termes présents dans l'équation obtenue, nous aurons besoin en
particulier des équivalents (pour $\eps$ tendant vers $0$) de toutes
les fonctions solutions des équations intérieures comme extérieures,
et des coefficients de Stokes des solutions intérieures. Certains des
calculs seront faits ici en préliminaires. 

Pour étudier le voisinage de $-1$,
le changement de variable qui s'impose est en tous cas le suivant
$$x=-1+\eps^{1/2}X\ \ , \ \ z=\eps^{-3/2}Y\,.$$
D'où l'équation intérieure, non singulièrement perturbée
$$Y\frac{dY}{dX}=\frac{2(-1+\sqrt\eps X)}{X^2}\Bigl(Y-\frac{1}{2X^3}\Bigr)
-\frac{a-1}{X^2}Y-\frac{2Y}{X}\Bigl(Y-\frac{1}{2X^3}\Bigr)\,.$$

Cette dernière équation a pour solution formelle la série $\hat Y=
\sum Y_n(X) \eps^{n/2}$ (et on sait que cette série correspond à une
vraie solution dans certains domaines), avec $Y_0$ solution de 
l'équation obtenue en posant $\eps=0$, soit
\begin{equation}\label{eq:BY0}
Y_0 \frac{dY_0}{dX}=\frac{-2}{X}\Bigl(Y_0-\frac{1}{2X^3}\Bigr)(Y_0+
\frac{1}{X})\,.
\end{equation}

\subsection{Quelques équivalents}

Si on cherche des équivalents pour les premières fonctions des séries 
$\hat y$ et $\hat Y$, on commence par écrire dans l'équation 
\eqref{Brus0} la fonction $z$ et le paramètre $a$ sous la forme
$z(\chi)=\Phi_0(\chi)+\eps \Phi_1(\chi)+ o(\eps)$ et
$a=1+a_1\eps+o(\eps)$; on part de
$$\Phi_0(\chi)=\frac{(1+\chi)^3}{2}\,,\ 
\Phi_0(x)=\frac{1}{2(1+x)^3}\,,$$
puis, comme 
$$\Phi_0 \Phi'_0=\frac{2\chi}{1+\chi}\Phi_1+ a_1\Phi_0$$ 
et que $\Phi_1$ doit être continue en $0$, on trouve 
$$a_1=3/2\ \ , \ \Phi_1(\chi)=\frac{3(1+\chi)^4}{8}(2+\chi)$$
autrement dit
$$\Phi_1(x)=\frac{3}{8}\frac{2+x}{(1+x)^5}\,.$$
On continue pour obtenir le terme suivant: dans l'équation du
Brusselator, les termes de degré $2$ en $\eps$ mèneront à l'équation
$$\Phi_0\Phi'_1+\Phi'_0\Phi_1=\frac{2\chi}{1+\chi}\Phi_2+a_1\Phi_1+
a_2\Phi_0+\frac{2\Phi_0\Phi_1}{1+\chi} $$
d'où
$$a_2=15/8\ \ , \ \Phi_2(\chi)=\frac{3(1+\chi)^4}{32}(6\chi^3+29\chi^2
+48\chi+30)$$
autrement dit
$$\Phi_2(x)=\frac{3}{32}\frac{5x^3+23x^2+42x+30}{(1+x)^7}\,.$$
Ensuite, pour la série $\hat Y$, on fait le changement de variable et
on identifie, avec $x=-1+X\sqrt\eps$, les deux séries
$$\eps^{3/2}\Bigl(\Phi_0(x)+\eps\Phi_1(x)+\eps^2 \Phi_2(x)\Bigr)=
Y_0(X)+\sqrt\eps Y_1(X)+\eps Y_2(X)+o(\eps)\,, $$ ce qui mène à
$$Y_0(X)\underset{X\rightarrow\infty}\sim
\frac{1}{2X^3}+\frac{3}{8X^5}+\frac{9}{16X^7}\ ,\ 
Y_1(X)\sim\frac{3}{8X^4}+\frac{33}{32X^6}\ ,\ Y_2(X)\sim\frac{3}{4X^5}\,.$$

\subsection{Coefficient de Stokes pour l'équation \eqref{eq:BY0}}

Si l'équation~\eqref{eq:BY0} a deux solutions conjuguées $Y_0^+$ et $Y_0^-$, 
bornées toutes deux pour les $X\in\RR$ assez grand, elles ont 
nécessairement la même série asymptotique (dont les termes sont donnés 
ci-dessus) et la différence entre les deux est donc
exponentionnellement petite. Elle correspond à des termes de Stokes;
nous allons montrer que:
\begin{Lemme}
La différence entre deux solutions bornées en l'infini de
l'équation~\eqref{eq:BY0} est équivalente, quand $X$ tend vers
l'infini, à
$$(Y_0^+-Y_0^-)(X)\sim\, 32i\sqrt{{2\pi}}X^4 
e^{-2X^2}\label{Stokes:br}\,.$$
\end{Lemme}
Pour parvenir à ce résultat, nous allons essayer de donner, sous une
certaine forme (car on ne peut espérer en donner explicitement les
solutions exactes), les solutions de l'équation.
\begin{align*}
\eqref{eq:BY0}\ Y_0' =&
-\frac{2}{X}Y_0+\frac{1}{X^4}-\frac{2}{X^2}+\frac{1}{X^5Y_0} \\ 
\intertext{On commence par changer de variable en posant $t=1/X$}
\frac{dY_0}{dt} =& \frac{2Y_0}{t}-t^2+2-\frac{t^3}{Y_0}\,,\\
\intertext{puis on introduit la nouvelle fonction $u(t)$ telle que
$Y_0(t)=t^2 u(t)$}
2tu+t^2 u'=& 2tu-t^2+2-t/u\\
  u' =& -1+ \frac{2}{t^2} -\frac{1}{tu}\,,\\ 
\intertext{et $v(t)$ vérifiant $u(t)=-t-2t^{-1}+v(t)$}
 \frac{dv}{dt} =& \frac{1}{t^2+2-t v(t)}\,.\\
\intertext{À présent plutôt que de regarder $v$ en fonction de $t$, on
prend la fonction réciproque: on cherche $t$ en fonction de $v$}
\frac{dt}{dv} =& t(v)^2 +2 -v t(v)\\
\intertext{On reconnaît une équation de Ricatti. La manière habituelle
de résoudre ce type d'équation est de poser $t=-\frac{z'(v)}{z(v)}$.}
-\frac{z''}{z}+\frac{{z'}^2}{z^2} =& 2+\frac{z'}{z} v +
\frac{{z'}^2}{z^2}\\
z''+vz'+2z =&0\\
\end{align*}
Cette dernière équation est une équation différentielle
hypergéométrique, qui a les deux solutions indépendantes suivantes
(parmi d'autres possibles):
\begin{align*}
z_1(v)=v \exp\left(-\frac{v^2}{2}\right)\ \text{ puis }\ &
z_2(v)=v  \exp\left(-\frac{v^2}{2}\right) \int^v_{i\infty}
\frac{1}{w^2} e^{w^2/2} dw\,,\\
\intertext{de dérivées}
z'_1(v)=(1-v^2)e^{-v^2/2} \ \text{ et }\ &
z'_2(v)=(1-v^2)e^{-v^2/2} \int^v_{i\infty} \frac{e^{w^2/2}}{w^2} dw
+\frac{1}{v}\,.\\
\end{align*}
On peut finalement écrire toute solution $t(v)$ sous la forme:
$$t(v)=\frac{C_1 (v^2-1)e^{-v^2/2}
+C_2(v^2-1)e^{-v^2/2} \int^v_{i\infty} \frac{e^{w^2/2}}{w^2} dw
-C_2/v}{C_1 v e^{-v^2/2}+C_2 v  
e^{-v^2/2} \int^v_{i\infty} \frac{1}{w^2} e^{w^2/2} dw}$$
Mais toutes ces solutions $t$ ne nous conviennent pas. En effet, on a vu
précédemment que les solutions $Y$ considérées devaient vérifier la
propriété suivante:
\begin{align*} 
Y_0(X)\underset{X\rightarrow\infty}\sim&\, \frac{1}{2X^3}\,,\\
\intertext{soit, après changements de variables,} 
-t^3-2t+t^2v \underset{t\rightarrow 0}\sim&\,  t^3/2\,,\\  
\intertext{ce qui suppose au moins que} 
t(v)\underset{v\rightarrow\infty}\sim&\, 2/v\,.\\
\end{align*} 
Or si on prend $C_2=0$, on trouve $t=\frac{v^2-1}{v}$, qui ne
vérifie pas cet équivalent.
 
Par contre, pour toute autre solution (dans le secteur où elle est
définie à l'infini), on trouve le bon équivalent. On montre en
utilisant des intégrations par parties que, quand
$v\rightarrow\infty$,  
$$\int_{i\infty}^v\frac{e^{w^2/2}}{w^2} dw=\frac{1}{v^3}e^{v^2/2}
\left(1+\frac{3}{v^2}+o(1/v^2)\right).$$
Ce qui fait que, si $C_2\neq 0$,
\begin{align}
z'(v) &\underset{v\rightarrow\infty}\sim -\frac{2}{v^3}\label{eq:z'(v)}\\
z(v)  &\underset{v\rightarrow\infty}\sim\frac{1}{v^2}\label{eq:z(v)}\\
\intertext{et donc}
t(v) &\underset{v\rightarrow\infty}\sim \frac{2}{v}\notag
\end{align}
ce qui donne une solution $Y_0$ qui se comporte bien comme prévu à
l'infini. 

Comme $C_2\neq 0$ dans les cas qui nous intéressent, on écrira
désormais $t(v)$ avec une seule constante $C=C_1/C_2$:
$$t(v)=\frac{(v^2-1)e^{-v^2/2}\left(\int_{i\infty}^v
\frac{e^{w^2/2}}{w^2} dw +C\right)-\frac{1}{v}}{ve^{-v^2/2}
\left(\int_{i\infty}^v \frac{e^{w^2/2}}{w^2} dw +C\right)}$$

Si alors $C=0$, la fonction $t^+(v)$ est bien définie pour tout $v$
tel que $\arg(v)\in [-\pi/4+\delta,5\pi/4-\delta], \delta>0$, parce
que pour $\arg(v^2)\in ]-\pi/2,\pi/2[$, l'exponentielle $e^{-v^2/2}$
est bien bornée et que pour $\arg(v)\in [\pi/4,3\pi/4]$, ce sont les
intégrales qui sont effectivement exponentionnellement petites. Si on
repasse à la variable $X$, cela donne l'existence d'une solution $Y_0^+$ 
au moins pour tout $X$ assez grand sur l'axe réel et la montagne Nord
(cf.~figure~\ref{fig:Brusselr}) d'après l'équivalent $X=1/t(v)\sim 2v$. 

Pour $C=i\sqrt{2\pi}$, on peut réécrire $t(v)$ sous la forme suivante:
$$t^-(v)=\frac{(v^2-1)e^{-v^2/2}
\int^v_{-i\infty}\frac{e^{w^2/2}}{w^2}dw -\frac{1}{v}}{v e^{-v^2/2}
\int^v_{-i\infty} \frac{1}{w^2} e^{w^2/2} dw}$$ 
En effet, en intégrant par parties 
\begin{align*}
\int^{+i\infty}_{-i\infty} \frac{\exp(w^2/2)}{w^2} dw &=
\left[-\frac{\exp(w^2/2)}{w}\right]^{i\infty}_{-i\infty}
+\int^{+i \infty}_{-i\infty} \exp(w^2/2)dw\\
\intertext{et avec le changement de variable $t=iw/\sqrt{2}$,}
&= \sqrt{2}/i \int^{-\infty}_{+\infty} \exp(-t^2)dt=i\sqrt{2{\pi}}
\end{align*}
La solution $t^-$ donne une solution $Y_0^-$ pour des $X$ dans un
domaine symétrique (par la conjugaison complexe) au domaine image de
$t^+$. 

\smallskip
On considère à présent des variables complexes $v^+$ et $v^-$ telles
que $$t^+(v^+)=t^-(v^-) \in \RR.$$
Pour des raisons de symétrie, ces deux nombres $v^\pm$ sont {\it a
priori} conjugués. On souhaite étudier la différence entre ces deux
nombres quand ils tendent vers $+\infty$. La méthode sera semblable à
celle utilisée pour l'équation de Van der Pol.
\begin{align*}
t^+(v^+)=&\,t^-(v^-) \\
\frac{{z^+}'(v^+)}{z^+(v^+)} =&\, \frac{{z^-}'(v^-)}{z^-(v^-)} \\
			     =&\, \frac{{z^+}'(v^-)
\left(1+i\sqrt{2{\pi}}\frac{(v^2-1)e^{-v^2/2}}{{z^+}'(v^-)}\right)}
{z^+(v^-)\left(1+i\sqrt{2{\pi}}\frac{ve^{-v^2/2}}{{z^+}(v^-)}\right)}
\\
\intertext{Comme les termes $e^{-v^2/2}/z^{(')}$ sont
exponentionnellement petits quand $v\rightarrow\infty$, on peut écrire
des développements limités:} 
			     =&\, \frac{{z^+}'(v^-)}{z^+(v^-)}
\left(1+i\sqrt{2{\pi}}\left[\frac{(v^2-1)e^{-v^2/2}}{{z^+}'(v^-)}
-\frac{ve^{-v^2/2}}{{z^+}(v^-)}\right]+o(\ldots)\right)\\
t^+(v^+)-t^-(v^-)\sim &\, t^+(v^-) i\sqrt{2{\pi}}e^{-v^2/2}
\left[\frac{(v^2-1)}{{z^+}'(v^-)}
-\frac{v}{{z^+}(v^-)}\right]\\
\end{align*}
\begin{align*}
\intertext{D'où, avec les équivalents pour $z$ et $z'$ trouvés en
\eqref{eq:z(v)} et \eqref{eq:z'(v)}}
t^+(v^+)-t^-(v^-)\sim &\,  -i\sqrt{2{\pi}}e^{-v^2/2} v^4 \\
\intertext{Or, comme $v^+$ et $v^-$ sont proches, on peut aussi
écrire}
t^+(v^+)-t^-(v^-)\sim &\left(v^+-v^-\right) {t^+}'(v^+)\,. \\
\intertext{Donc finalement}
v^+-v^- \sim& \frac{i}{2}\sqrt{2{\pi}}e^{-v^2/2} v^6\,. \\  
\intertext{On en déduit facilement les équivalents correspondants pour
$(Y_0^+-Y_0^-)(X)$, sachant que}
Y_0^+-Y_0^-=&\, \frac{v^+-v^-}{X^2}\ \text{ et }\ v\sim\frac{2}{t}\sim
2X \,.\\
\intertext{On obtient}
(Y_0^+-Y_0^-)(X)\sim&\, 32i\sqrt{{2\pi}}X^4 
e^{-2X^2}\ \,, \eqref{Stokes:br}\\ 
\intertext{avec donc un coefficient de Stokes}
 & \mathcal{C}=32i\sqrt{2\pi} \,.
\end{align*}

\subsection{L'équation différentielle pour l'estimation de $a^+-a^-$}

Comme dans la première partie on cherche une équation différentielle 
en $d(x)=(z^+-z^-)(x)$; on pose aussi $b=a^+-a^-$.

À partir de l'équation \eqref{Brus0} ci-dessous
\begin{align*}
\eps \frac{dz^+(x)}{dx}&=\frac{2x}{(1+x)^2} 
\Bigl(1-\frac{1}{2z^+(1+x)^3}\Bigr)-\frac{a^+-1}{(1+x)^2} 
-\Bigl(z^+-\frac{1}{2(1+x)^3}\Bigr)\frac{2\eps}{1+x}\\
\intertext{et de la même équation écrite pour $(z^-,a^-)$, on arrive, en
faisant la différence, à}
\eps d'(x)&=\frac{x}{(1+x)^5}\Bigl(\frac{1}{z^-}-\frac{1}{z^+}\Bigr)
-\frac{b}{(1+x)^2}-d\frac{2\eps}{1+x}\\
\eps d'(x)&=d\left(\frac{x}{(1+x)^5z^-z^+}-\frac{2\eps}{1+x}\right)
- \frac{b}{(1+x)^2}\,.
\end{align*}
D'où, pour tout $x_0$ réel,
\begin{equation*}
d(x_0)=\frac{-b}{\eps} \int^{x_0}_{+\infty}
\frac{\exp\left({\frac{G_T(x_0)-G_T(t)}{\eps}}\right)}{(1+t)^2} dt
\end{equation*}
avec $G_T(x)=\int^x_0 \frac{t}{(1+t)^5z^+(t)z^-(t)} dt -2\eps \ln
(1+x)$, soit
\begin{equation}\label{eq:Bb}
d(x_0)=\frac{-be^{G(x_0)/\eps}}{\eps(1+x_0)^2} 
\int^{x_0}_{+\infty} e^{-G(t)/\eps} dt
\end{equation}
où $G(x)=\int^x_0 \frac{t}{(1+t)^5z^+(t)z^-(t)} dt$.

\subsection{Calcul de l'équivalent}

{On peut estimer chacun des termes du produit de \eqref{eq:Bb} (on
travaillera ici avec un $\eps$ réel positif).}

Avec la méthode du point col,
\begin{equation}\label{Bint}
\int^{x_0}_{+\infty} e^{-G(t)/\eps} dt \sim
-\sqrt{2\pi\eps}|z(0)|\ \ \text{ avec }\qquad|z(0)|\sim \Phi_0(0)
=\frac{1}{2}\,.
\end{equation}

Puis
\begin{equation}\label{Bd}
d(x_0)\sim \eps^{-3/2} \bigl(Y_0^+(X_l)-Y_0^-(X_l)\bigr)\,.
\end{equation}

Reste $e^{G(x_0)/\eps}$, où on prend $x_0=-1+X_l\sqrt\eps$. On coupe
l'intégrale correspondant à $G(x_0)$ en deux, en un point
$\lambda(\eps)$ dépendant peu de $\eps$, comme par exemple
$\lambda(\eps)=-1+{\eps}^{1/20}$; le but étant qu'à droite on puisse
utiliser l'estimation $z(t,\eps)=\Phi_0(t)+\eps \Phi_1(t)+
\eps^2 \Phi_2(t)+o({\eps}^{2-2/20})$ (théorème ~\ref{th:-loin}); et 
à gauche, gr\^ace au corollaire \ref{cor:pabo}, l'estimation
$z(t,\eps)=\f1\eps (Y_0(t/\eps)+\sqrt\eps Y_1(t/\eps)+\eps Y_2(t/\eps))
+o(\sqrt\eps^{1-3/20}$).

$$G(x_0)=\int_0^\lambda \frac{t}{(1+t)^5z^+(t)z^-(t)} dt +
\int_\lambda^{x_0} \frac{t}{(1+t)^5z^+(t)z^-(t)} dt\,.$$

Pour $\mathcal{I}_1=\int_0^\lambda \frac{t}{(1+t)^5z^+(t)z^-(t)} dt$,
on utilise donc pour $z$ la forme $z=\Phi_0+\eps \Phi_1+ o(\eps)$ (on
a bien $o(\eps)$ car $\Phi_2$ est bornée), qui donne
$$\frac{1}{z^+(t)z^-(t)}=\frac{1}{\Phi_0^2}-2\eps\frac{\Phi_1}{\Phi_0^3}
+o(\eps)\,,$$
donc
\begin{align*}
\mathcal{I}_1&=\!\int_0^\lambda\!\! 4t(1+t)dt-
6\eps\!\int_0^\lambda\!\!\frac{t}{1+t}(2+t)dt\\
\mathcal{I}_1&=
\frac{4}{3}\lambda^3+2\lambda^2-\eps\Bigl(3(1+\lambda)^2-3
-6\ln(1+\lambda)\Bigr)\label{BRbemol}\,.\tag{$R\mathcal{I}_1$}
\end{align*}

Pour $\mathcal{I}_2=\int_\lambda^{x_0} \frac{t}{(1+t)^5z^+(t)z^-(t)}
dt$, on se place au voisinage de $t=-1$, et on change de variable:
$t=-1+X\sqrt\eps$, $z=\eps^{-3/2}Y$.
\begin{align*}
\mathcal{I}_2 =&\sqrt\eps
\int_{\frac{\lambda+1}{\sqrt\eps}}^{\frac{x_0+1}{\sqrt\eps}} 
\frac{-1+\sqrt\eps X}{\sqrt\eps^5 X^5 \eps^{-3}Y^-Y^+}dX
 =\eps \int_{\ell}^{X_l}
\frac{-1+\sqrt\eps X}{X^5 Y^-(X)Y^+(X)}dX\,,\\
\intertext{où on pose bien sûr $\ell=\frac{\lambda+1}{\sqrt\eps}$. Le
calcul ressemble à celui effectué plus haut dans le cas de van der
Pol; on commence par développer les fonctions dans l'intégrale suivant
$\eps$, avant d'estimer chacun des termes obtenus, en vérifiant que le
dernier terme, et les suivants, sont bien négligeables par rapport aux
précédents.}
\intertext{Comme $Y(X)\sim \frac{1}{2X^3}$, l'intégrale converge, y
compris pour $x_0=-1$, et}
\mathcal{I}_2 =&\,\eps\!\int^{X_l}_\ell\!\frac{-1+\sqrt\eps X}{X^5 Y^-(X)Y^+(X)}dX\\
  =&\,\eps\!\int^{X_l}_\ell\!\frac{-1+\sqrt\eps X}{X^5 Y_0^-(X)Y_0^+(X)
\left(1+\sqrt\eps\frac{Y_1^+}{Y_0^+}+
\sqrt\eps\frac{Y_1^-}{Y_0^-}+\eps\frac{Y_2^+}{Y_0^+}
+\eps\frac{Y_2^-}{Y_0^-}+o(\eps)\right)}dX\\
  =&\,\eps\!\!\int^{X_l}_\ell\!\!\frac{-1+\sqrt\eps X}{X^5 Y_0^-(X)Y_0^+(X)}dX
-\eps^{3/2} \int^{X_l}_\ell\!\! \frac{-1+\sqrt\eps X}{X^5 Y_0^-(X)Y_0^+(X)}
\left(\frac{Y_1^+}{Y_0^+}+\frac{Y_1^-}{Y_0^-}\right)dX\\
 &\phantom{\eps} 
-\eps^2\int^{X_l}_\ell\frac{-1+\sqrt\eps X}{X^5 Y_0^-(X)Y_0^+(X)}
\left(\frac{Y_2^+}{Y_0^+}+\frac{Y_2^-}{Y_0^-}-
\left(\frac{Y_1^+}{Y_0^+}+\frac{Y_1^-}{Y_0^-}\right)^2
\right)dX +{o}(\eps)\\
  =&\underbrace{\int^{X_l}_\ell\frac{-\eps}{X^5Y_0^+Y_0^-}dX}_{(\text{A})}
+\underbrace{\int^{X_l}_\ell
\frac{\eps^{3/2}}{X^4Y_0^2}dX}_{(\text{B})}
+\underbrace{2\eps^{3/2}\int^{X_l}_\ell
\frac{1-\sqrt\eps X}{X^5Y_0^2}\frac{Y_1}{Y_0}dX}_{(\text{C})}\\
   &\phantom{\eps}+\underbrace{2\eps^2\int^{X_l}_\ell
\frac{1-\sqrt\eps X}{X^5Y_0^2}\left(\frac{Y_2}{Y_0}-
2\frac{Y_1^2}{Y_0^2}\right)dX}_{(\text{D})}+o(\eps)\\
\intertext{On s'occupe de chacun des termes de la somme}
\text{D}\sim&\, 2\eps^2\int^{X_l}_\ell\frac{1-\sqrt\eps X}{X^5Y_0^2}
\left(\frac{3}{2X^2}+\f{9}{8X^2}\right)dX
\sim 2\eps^2\int_\ell (1-\sqrt\eps X)
\frac{6+9/2}{X} dX\\
   =&o(\eps)\ \,\,
\left( \text{Rappel: $\ \ell=\frac{\lambda+1}{\sqrt\eps}$} \right)\\
\text{C} \sim&\, 2\eps^{3/2}\int^{X_l}_\ell\frac{1-\sqrt\eps X}{X^5Y_0^2}
\frac{Y_1}{Y_0}dX\sim2\eps^{3/2}\int_\ell(1-\sqrt\eps X)4X\frac{3}{4X}dX\\
   \sim&\,\eps\left(-6(\lambda+1)+3(\lambda+1)^2\right)\\
\text{B} \sim&\, \eps^{3/2} \int_\ell \frac{dX}{X^4\frac{1}{4X^6}
\left(1+\frac{3}{2X^2}+o\left(\frac{1}{X^3}\right)\right)}\\
   \sim&\, \eps^{3/2} \int_\ell\left(4X^2-6+o\left(1/X\right)\right)dX
\sim -\eps^{3/2}\left(\frac{4}{3}\ell^3-6\ell\right)\\
   \sim& -\frac{4}{3}(\lambda+1)^3+6\eps (\lambda+1)\\
\end{align*}
{Pour A, il faut être plus précis et distinguer $Y_0^+$ de 
$Y_0^-$; on repasse donc par l'équation différentielle \eqref{eq:BY0},}  
\begin{align*}
(Y_0^+-Y_0^-)'=&-\frac{2}{X}(Y_0^+-Y_0^-)
+\frac{1}{X^5}\frac{-Y_0^++Y_0^-}{Y_0^-Y_0^-}\\
\intertext{qui donne, en intégrant l'équation par la méthode de variation
de la constante, pour tous $\xi$, $\xi_0$ positifs:}
(Y_0^+-Y_0^-)(\xi)=&(Y_0^+-Y_0^-)(\xi_0)\exp\left(-2\ln\frac{\xi}{\xi_0}
-\int^\xi_{\xi_0}\frac{1}{t^5Y_0^-(t)Y_0^+(t)}dt\right)\\
 \int^{X_l}_{\ell}\frac{-1}{t^5Y_0^-Y_0^+}dt=&
\ln\frac{(Y_0^+-Y_0^-)(X_l)}{(Y_0^+-Y_0^-)(\ell)}+2\ln\frac{X_l}{\ell}\,.\\
\end{align*}
{D'après le résultat \eqref{Stokes:br} obtenu plus haut,}
\begin{align*}
\text{A} \sim& \,\eps\left(\ln\left(Y_0^+-Y_0^-\right)(X_l)
-\ln\left({\cal C} \ell^4 e^{-2\ell^2}\right)+2\ln\frac{X_l}{\ell}\right)\\
	 \sim& \,\eps\biggl(\ln\left(Y_0^+-Y_0^-\right)(X_l) -\ln{\cal C}
-4\ln(\lambda+1)+2\frac{(\lambda+1)^2}{\eps}+2\ln\eps\\
            &\phantom{\ln\left(Y_0^+-Y_0^-\right)(X_l)}
+2\ln(x_0+1)-2\ln(\lambda+1)\biggr)\,.\\ 
\end{align*}
Enfin,
\begin{multline}\tag{$R\mathcal{I}_2$}
\mathcal{I}_2\sim\text{A+B+C+D}\sim 
-\frac{4}{3}\lambda^3-2\lambda^2+\frac{2}{3}+\eps\biggl(
\ln\left(Y_0^+\!-\!Y_0^-\right)(X_l)+2\ln(x_0+1)\\
-6\ln(\lambda+1)+2\ln\eps-\ln{\cal C}+3(\lambda+1)^2\biggr)\,,
\end{multline}
\begin{align*}
\intertext{puis en reprenant le résultat \eqref{BRbemol}}
\exp\left(\frac{\mathcal{I}_1+\mathcal{I}_2}{\eps}\right) &\sim 
\exp\left(\frac{2}{3\eps}\right)\left(Y_0^+-Y_0^-\right)(X_l)(x_0+1)^2
\frac{\eps^2}{\cal C}e^3\,,\\
\intertext{donc avec \eqref{Bd} et \eqref{Bint}, l'équation
\eqref{eq:Bb} devient}
\eps^{-3/2}\left(Y_0^+\!-\!Y_0^-\right)\!(X_l)&
\sim b\frac{\sqrt{2\pi\eps}e^3
\left(Y_0^+-Y_0^-\right)(X_l)(x_0+1)^2\eps^2e^{2/3\eps}}{2{\cal C}\eps
(x_0+1)^2}\\
\eps^{-3/2} &\sim b\frac{e^3\sqrt{2\pi}}{2{\cal C}}\eps^{3/2}
e^{\frac{2}{3\eps}}\\
b&\sim \frac{2{\cal C}}{\sqrt{2\pi}e^3}
\frac{\exp\left(-\frac{2}{3\eps}\right)}{\eps^3}\,.
\end{align*}
Finalement, on arrive au résultat du
\begin{Th}\label{th:eqbrussel}
$$(a^+-a^-)(\eps)
\sim 64ie^{-3}\exp\left(-\frac{2}{3\eps}\right)\eps^{-3}$$
\end{Th}

\section{Conséquence pour le développement en série de $a(\eps)$}

\begin{Cor}
Les coefficients $a_n$ du développement asymptotique de $a^+$ et $a^-$
ont une croissance Gevrey, et vérifient, quand $n\rightarrow\infty$, 
$$a_n\sim\, 54n^2\left(\frac{3}{2}\right)^{n} n!$$
\end{Cor}

Étudions d'abord ce qui se passe pour la fonction $a^+(\eps)$ quand on
fait varier l'argument $\theta$ de $\eps$. Il est clair que
$a^+(\eps)$ et $v^+(x,\eps)$ dépendent holomorphiquement de $\eps$
dans des secteurs centrés en $0$ (mais pas en $0$
lui-m\^eme). Regardons en particulier ce qui se passe quand on fait
varier $\theta$ entre $0$ à $-2\pi$. En étudiant le relief
correspondant à $\theta$, on peut déduire le domaine d'existence de la
solution $v^+$ correspondante (qui varie continûment avec $\theta$ lui
aussi); cette évolution du domaine est représenté
figure~\ref{fig:tourne}. On remarque tout de suite que pour $\eps$
réel, le domaine de $z^+$ pour $\eps e^{-2i\pi}$ correspond au domaine
de $z^-$ pour la valeur $\eps$. Par unicité des solutions, on en
déduit que $a^+(\eps e^{-2i\pi})=a^-(\eps)$. C'est de cette propriété
que l'on se servira ci-dessous.

\begin{figure}[ht!]
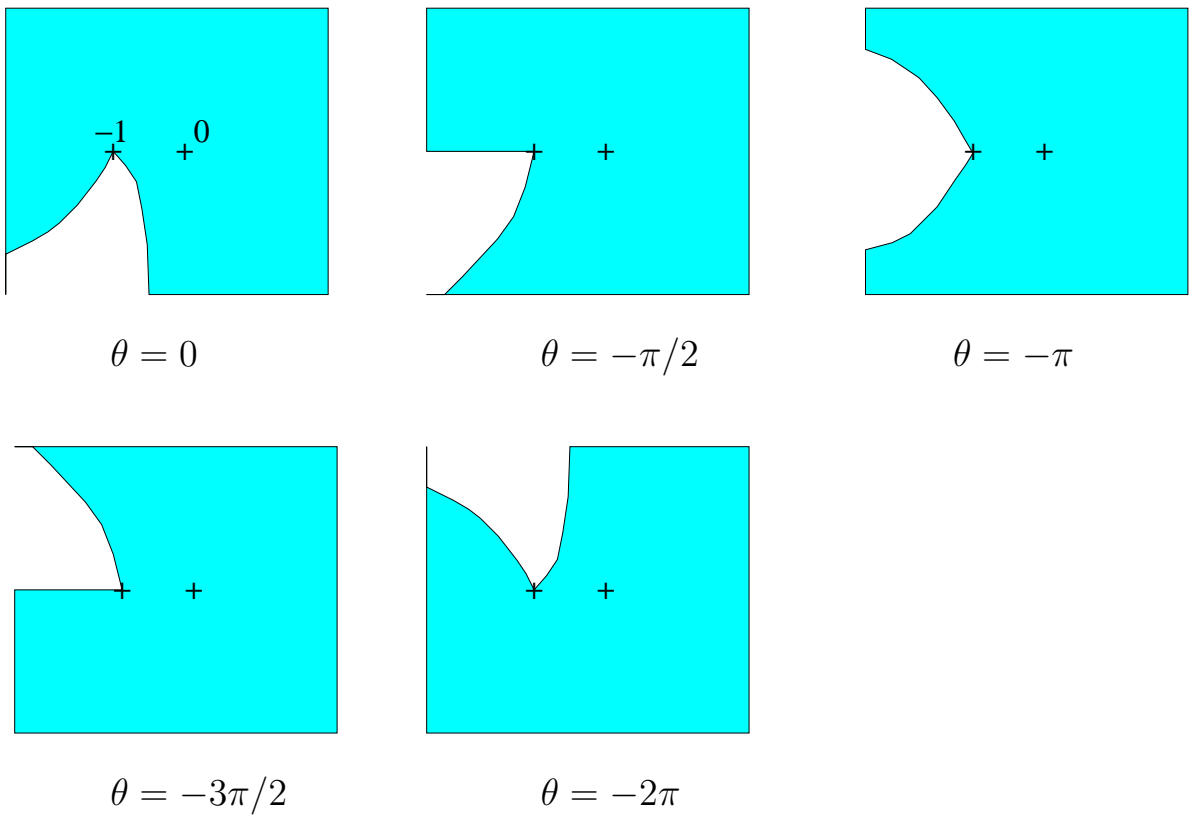

\caption{Domaine d'existence de $z^+$ pour différentes valeurs de 
$\theta=\arg(\eps)$.}
\centerline{\input tourne.pstex_t }
\label{fig:tourne}
\end{figure}

On peut faire exactement les mêmes déductions que pour l'équation de
Van der Pol, les deux solutions $y^+(x,\eps)$ et $y^-(x,\eps)$ étant
reliées de la même manière. On intègre à nouveau sur un chemin autour
de $\eps=0$ (celui de la figure~\ref{fig:integ}), pour trouver que
\begin{align*}
a_n\sim&\, \frac{1}{2i\pi}\int_0^{\eps_0} b(\eps) \eps^{-n-1} d\eps\,.\\
\intertext{En utilisant le théorème~\ref{th:eqbrussel},}
   \sim&\, \frac{1}{2i\pi}\int_0^{\eps_0} 64ie^{-3} 
  \exp\left(-\frac{2}{3\eps}\right)\eps^{-3}\eps^{-n-1} d\eps \\
   \sim&\, \frac{32e^{-3}}{\pi}\int_0^{\eps_0} 
\exp\left(-\frac{2}{3\eps}\right) \eps^{-n-4} d\eps \,.\\
\intertext{Puis on pose $t=2/3\eps$ et on transforme l'intégrale en
   une autre, équivalente:}
a_n\sim&\, \frac{32e^{-3}}{\pi} \int_0^{+\infty} 
\left(\frac{3}{2}\right)^{n+3} e^{-t}t^{n+2} dt \\
a_n\sim&\, \frac{32e^{-3}}{\pi} \left(\frac{3}{2}\right)^{n+3} (n+2)!
\intertext{Soit, sous une forme plus normale,}
a_n\sim&\, \frac{108e^{-3}}{\pi} n^2\left(\frac{3}{2}\right)^{n} n!
\end{align*}

\end{document}